\documentclass{amsart}
\usepackage[utf8]{inputenc}

\title{Weight Filtrations and Derived Motivic Measures}
\author{Anubhav Nanavaty}
\thanks{This work was supported in part by NSF Grant No. DMS-1944862}

\usepackage{amsmath}
\usepackage{graphicx}
\usepackage{scalerel}
\usepackage{stackengine,wasysym}
\usepackage{hyperref}

\usepackage{amssymb}
\usepackage{amsthm}
\usepackage{mathtools}
\usepackage{adjustbox}
\usepackage{dsfont}
\usepackage{MnSymbol}
\usepackage{mathrsfs}
\usepackage{stmaryrd}
\usepackage[all]{xy}
\usepackage[mathcal]{eucal}
\usepackage{verbatim}  
\usepackage{fullpage}  
\usepackage{hyperref}
\usepackage{setspace}
\usepackage{tikz-cd}
\usetikzlibrary{positioning}
\newcommand{\A}{\mathbb{A}}

\newcommand{\C}{\mathbb{C}}

\newcommand{\F}{\mathbb{F}}
\newcommand{\G}{\mathbb{G}}

\newcommand{\N}{\mathbb{N}}

\newcommand{\Q}{\mathbb{Q}}

\renewcommand{\S}{\mathbb{S}}

\newcommand{\V}{\mathbb{V}}

\newcommand{\Z}{\mathbb{Z}}

\newcommand{\Spec}{\mathrm{Spec}}
\newcommand{\Pro}{\mathrm{Pro}}
\newcommand{\Pre}{\mathrm{Pre}}
\newcommand{\Sch}{\mathrm{Sch}}
\newcommand{\Ho}{\mathrm{Ho}}

\newcommand{\Corr}{\mathrm{Corr}}
\newcommand{\s}{\mathrm{s}}
\newcommand{\cdh}{\mathrm{cdh}}
\newcommand{\sShv}{\mathrm{sShv}}
\newcommand{\GW}{\mathrm{GW}}
\newcommand{\Chow}{\mathrm{Chow}}
\newcommand{\Ch}{\mathrm{Ch}}
\newcommand{\eff}{\mathrm{eff}}
\newcommand{\perf}{\mathrm{perf}}
\newcommand{\Hom}{\mathrm{Hom}}

\newcommand{\End}{\mathrm{End}}
\renewcommand{\l}{\langle}
\renewcommand{\phi}{\varphi}
\renewcommand{\r}{\rangle}

\renewcommand{\o}{\overline}

\providecommand{\o}{\overline}
\usepackage{rotating}
\usepackage{amssymb}

\newtheorem{theorem}{Theorem}[section]

\newtheorem{proposition}[theorem]{Proposition}
\newtheorem{lemma}[theorem]{Lemma}
\newtheorem{corollary}[theorem]{Corollary}

\newtheorem{remark}[theorem]{Remark}
\newtheorem{convention}[theorem]{Convention}
\theoremstyle{definition}
\newtheorem{definition}[theorem]{Definition}

\newtheorem{hypothesis}[theorem]{Hypothesis}

\newtheorem{question}{Question}
\numberwithin{equation}{subsection}
\begin{document}

\maketitle
\begin{abstract}
Let $k$ be a field admitting resolution of singularities. We lift a number of motivic measures such as the Gillet--Soul\'{e} measure and the compactly supported $\A^1$-Euler characteristic to derived motivic measures in the sense of Campbell--Wolfson--Zakharevich, answering various questions in the literature. The obstruction to such lifts is that the Gillet--Soul\'{e} weight complex of a variety is built from data that is a priori functorial only after passing to a homotopy category. We remove this obstruction by showing that the collection of all weight complexes of $X$ assembles into a canonical \emph{weakly constant} pro-object in a Waldhausen category $\s\V_*^{\perf}$ of simplicial smooth projective varieties. On the way, we prove a statement of independent interest: under mild assumptions, the $K$-theory of a Waldhausen category is equivalent to the $K$-theory of its weakly constant pro-objects. This leads us to a new proof and generalization of the existence of the Gillet--Soul\'{e} weight filtration to both the unstable and stable motivic homotopy category. Lastly, we organize all of the various maps of spectra into a single homotopy commutative diagram out of $K(\mathcal{V}_k)$, the Zakharevich $K$-theory of varieties. 
\end{abstract}
\maketitle
\section{Introduction}

Many constructions that probe the geometry of algebraic varieties provide \emph{motivic measures}, or ring homomorphisms of the form:
\[
K_0(\mathcal{V}_k)\to R
\]
where $K_0(\mathcal{V}_k)$ is the Grothendieck ring of $k$-varieties, the free abelian group generated by isomorphism classes of separated $k$-schemes of finite type, modulo the relations: \[[Z]+[X\setminus Z]=[X]\]
for $Z\subset X$ a Zariski closed subvariety. The ring structure is given by $[X]\cdot [Y]:=[X\times Y]$. The target of the ring homomorphism, $R$, is usually $\pi_0$ of a commutative ring spectrum. Motivic measures provide useful ways to distinguish varieties, or prove facts about how they decompose. 

Throughout this paper, let $k$ be a field admitting resolution of singularities. Under this assumption, \cite{GS} construct the integral Chow motive for a variety $X$, giving rise to a motivic measure:
\[
K_0(\mathcal{V}_k)\to K_0(\Ch^{\mathrm{b}}(\Chow^\eff))
\]
where $\Ch^{\mathrm{b}}(\Chow^\eff)$ is the category of bounded chain complexes of effective Chow motives (see Section 2.1 for a full review of the construction). 

Similarly, the compactly supported $\A^1$-Euler characteristic provides a motivic measure sending
\[
[X]\mapsto \chi_{c,\A^1}(X)\in \GW(k)=\pi_0(\End(1_{\mathrm{SH}(k)}))
\]
where $\GW(k)$ is the Grothendieck-Witt ring over $k$, and $\End(1_{\mathrm{SH}(k)})$ is the endomorphism spectrum of the unit object in $\mathrm{SH}(k)$, the motivic stable homotopy category. This invariant can be computed for many examples, yielding interesting geometric information about the isomorphism class of $X$  (see \cite{MBOWZ} for an exposition). Notably, the target of both of these motivic measures is (non-trivially) $\pi_0$ of a commutative ring spectrum, and in \cite{zaka}, Zakharevich showed that the same is true of the source by constructing a spectrum $K(\mathcal{V}_k)$ whose path components recover the Grothendieck ring of varieties. It is therefore natural to ask whether these motivic measures lift to {\em derived motivic measures}, i.e. maps of spectra that recover the motivic measure on path components.  One reason for doing so is that such lifts allow one to probe the higher K-theory of varieties $K_i(\mathcal{V}_k)$, which to quote \cite{CWZ} ``remembers not only that certain varieties decompose into other varieties, which is what $K_0(\mathcal{V}_k)$ does, but also how." This approach is manifest in the work of \cite{CWZ}, which lifts the derived motivic measure associated to compactly supported \'{e}tale cohomology $[X]\mapsto [H^*_{\text{\'{e}t}}(X,\Q_\ell)]$ to show that for $p$ prime, there is a splitting $K(\mathcal{V}_{\F_p})\simeq \S\vee \tilde{K}(\mathcal{V}_{\F_p})$, where $\S$ is the sphere spectrum and $\tilde{K}_i(\mathcal{V}_{\F_p})$ is non-zero for arbitrarily large $i\in\N$. In \cite{BG}, it is further shown that for $k\subset \C$, there are infinitely many linearly independent non-torsion elements of arbitrarily high degree in $K(\mathcal{V}_{k})$ by using the derived motivic measure associated to the mixed Hodge structure constructed in \cite{BGN}. 

This paper grew out of an attempt to answer the following questions from the literature:
\begin{question}[7.5 of \cite{CWZ}]
Can one lift the Gillet--Soul\'{e} motivic measure?
\end{question}
\begin{question}[1.14 of \cite{MBOWZ}]
Can one lift the compactly supported $\A^1$-Euler characteristic?
\end{question}
Unfortunately, motivic measures are usually constructed non-functorially and at the level of homotopy categories, meaning that an enhancement to a derived motivic measure involves significant theory building and a deep understanding of the structure of the choices made in the original construction. For example, the measure referenced in Question 2 is defined only on smooth projective varieties, and extended to all of $K_0$ via compactifying and resolving singularities (both non-functorial processes). The main goal of this paper answer these questions affirmatively and in the process construct a plethora of other derived motivic measures that are interesting in their own right.
\begin{theorem}\label{mainthm}
Let $k$ be a field admitting resolution of singularities. There exists a homotopy commutative diagram of spectra
\[\begin{tikzcd}
    &&K(\sShv_{\cdh}^{\perf}(\Sch(k))_*\arrow[r,"\ref{prop:symspt}"]\arrow[d,"\ref{thm:doldkan}"]&K(\mathrm{SH}^{\perf}(k))\arrow[r,"\ref{cor:trace}"]&\End(1_{\mathrm{SH}(k)})\\
K(\mathcal{V}_k)\arrow[r,"\ref{cor:sVmotivicmeasure}"]&K(\s\V_*^{\perf}(k))\arrow[ur,"\ref{prop:sshv}"]\arrow[r,"\ref{prop:Mdg}"]\arrow[dr,"\ref{def:Chow}"]&K(\mathrm{C}^{\mathrm{gm}}(k,R))\arrow[d,"\ref{rem:thmheart}"]\\
    &&K(\Ch^{\perf}(\Chow^{\eff}(R
    )))
\end{tikzcd}\]
where $K(\mathcal{V}_k)$ denotes the Zakharevich spectrum of $k$-varieties, constructed in \cite{zaka}, and the other spectra are $K$-theory spectra of Waldhausen categories. At the level of $\pi_0$, the map $K_0(\mathcal{V}_k)\to K_0(\Ch^{\perf}(\Chow^{\eff}(R)))$ recovers the Gillet--Soul\'{e} measure, and the map $K_0(\mathcal{V}_k)\to \pi_0(\End(1_{SH(k)}))=\GW(k)$ lifts the compactly supported $\A^1$-Euler characteristic. We label each map with the associated theorem/proposition/corollary in the paper where we prove that such a map exists.\begin{footnote}
        {The maps labelled \ref{prop:symspt} and \ref{cor:trace} are constructed in \cite{Ro}, but we cite them in our paper for the convenience of the reader.}
    \end{footnote} 
\end{theorem}
The majority of the work done is to construct the map of Corollary \ref{cor:sVmotivicmeasure}.
The key insight behind Theorem \ref{mainthm} is to re-interpret and generalize the construction of the Gillet--Soul\'{e} weight complex. In particular, we construct a category of ``well-pointed" simplicial smooth projective varieties, $\s\V_*^\perf$ in Definition \ref{def:sVperf}, along with a subcategory of weak equivalences $\textbf{w}(\s\V_*^\perf)$ in Definition \ref{def:we}. Let $L(\mathcal{C},\mathcal{W})$ denote the $\infty$-categorical localization of a 1-category $\mathcal{C}$ with respect to a wide subcategory $\mathcal{W}\subset \mathcal{C}$ of weak equivalences, and let $\Ho(\mathcal{C}):=hL(\mathcal{C},\mathcal{W})$ denote the classical 1-categorical localization. By carefully studying the work done in \cite{GS}, we obtain an assignment $X\mapsto M(X)\in \Ho(\s\V_*^{\perf})$ which we call the Gillet--Soul\'{e} premotive (see Theorem \ref{prop:boundedweight}). We prove:
  \begin{theorem}[Theorem \ref{RefinedGS}]\label{mainthm2}
      Let $k$ be a field admitting resolution of singularities, and $\mathcal{V}_k$ denote the category of $k$-varieties. The assignment to a variety of its \emph{Gillet--Soul\'{e} premotive}
     \[
     X\mapsto M(X)\in \Ho(\s\V_*^{\perf})
     \]
    further extends to $\infty$-functors:
     \[
     f^!:(\mathcal{V}_k^{\mathrm{open}})^{op}\to L\Big(\s\V_*^{\perf}, \textbf{\emph{w}}(\s\V_*^{\perf})\Big)
     \]
     \[
     f_!:(\mathcal{V}_k^{\mathrm{closed}})\to L\Big(\s\V_*^{\perf}, \textbf{\emph{w}}(\s\V_*^{\perf})\Big)
     \]
     where $\textbf{\emph{w}}(\s\V_*^{\perf})$ is the subcategory of weak equivalences given by the Waldhausen structure on $\s\V_*^{\perf}$ in Definition \ref{def:we}. Here, $\mathcal{V}_k^{\mathrm{open}}$ is the subcategory of open immersions of $\mathcal{V}_k$, and $\mathcal{V}_k^{\mathrm{closed}}$ is the subcategory of closed embeddings of $\mathcal{V}_k$. $f_!$ and $f^!$ also agree on isomorphisms in $\mathcal{V}_k$.
 \end{theorem} 
 To prove this theorem, it's enough to work with relative categories and hammock localizations defined in \cite{DK} as models for $\infty$-categories and their localizations. This theorem is the main tool that allows us to construct the diagram in Theorem \ref{mainthm}, and results in a new proof of the Gillet--Soul\'{e} weight filtration being well-defined. In the original work, Gillet and Soul\'{e} \cite{GS} use key properties of the Gersten resolution to prove well-definedness, whereas in our work we rely on the `realization functors' out of $\s\V_*^{\perf}$ and basic geometric arguments from \cite{GS}.

 We now give a brief summary of the work done in each section of the paper. Section 2 reviews Chow motives and the Gillet--Soul\'{e} construction, along with some preliminaries on Waldhausen categories and pro-objects. In Section 3, the results of Section 2 allow us to put a certain Waldhausen structure on $\s\V_*^\perf(k)$ (see Corollary \ref{prop:ffwcsV}). In Section 4, we put a Waldhausen structure on the category of its `weakly-constant' pro-objects $\Pro^{wc}(\s\V_*^\perf(k))$ (see Definition \ref{def:prosV} for a proper definition) and show the equivalences both of their Dwyer-Kan hammock localizations and $K$-theory spectra via the following general theorem:
 \begin{theorem}[Theorems \ref{thm:DKequiv} and \ref{thm:proobj}]
 Let $\mathcal{C}$ be a Waldhausen category with FFWC such that all morphisms are weak cofibrations,\begin{footnote}
     {See Definition \ref{def:ffwc} for a precise definition of FFWC and weak cofibrations.}
 \end{footnote}and such that every undercategory $c\downarrow \mathcal{C}$ also inherits a Waldhausen with FFWC structure from $\mathcal{C}$. Let $\textbf{\emph{w}}\mathcal{C}$ denote the subcategory of weak equivalences. The fully faithful inclusion $\mathcal{C}\hookrightarrow \Pro^{wc}(\mathcal{C})$ induces an equivalence of $\infty$-categories:
 \[L\Big(\Pro^{wc}(\mathcal{C}),\textbf{\emph{w}}(\Pro^{wc}(\mathcal{C})\Big)\simeq L\Big(\mathcal{C},\textbf{\emph{w}}\mathcal{C}\Big)\]
along with an equivalence of spectra:
     \[K(\Pro^{wc}(\mathcal{C}))\simeq K(\mathcal{C})\]
 \end{theorem}
 In order to do this, we heavily use results from \cite{BM}. In Section 5, we prove the most technical lemmas of the paper. The most important lemma is Proposition \ref{prop:proobj}, which proves that for any variety $X$, the collection of all Gillet--Soul\'{e} pre-motives of $X$ forms a pro-object in a certain category, which we call the category of \emph{pro-weight complexes}. The work in Sections 4 and 5 demonstrates that the construction of Gillet--Soul\'{e} gives a more general notion of a weight filtration than in \cite{GS}, and so in Section 6 we begin by proving Theorem \ref{mainthm2}. We also construct the map $K(\mathcal{V}_k)\to K(\s\V_*^{\perf})$ in Corollary \ref{cor:sVmotivicmeasure}. Theorem \ref{mainthm2} and Corollary \ref{cor:sVmotivicmeasure} are the main results of the paper. In Sections 7-10, the remaining maps in Theorem \ref{mainthm} are constructed. In doing so, we show that the Gillet--Soul\'{e} weight filtration is well-defined in both the homotopy category of simplicial cdh sheaves and the stable homotopy category, generalizing their original construction.\\\\
\tableofcontents

\textbf{Acknowledgements.} The author would like to acknowledge Jesse Wolfson, his advisor, for suggesting this project and for countless insightful discussions and suggestions, and the anonymous referee on a previous submission that made extensive corrections and suggestions on a previous draft of the paper. The author acknowledges Denis-Charles Cisinski and Marc Hoyois for helpful discussions. Lastly, the author acknowledges Liam Keenan, Vladimir Sosnilo, Thomas Brazelton, Oliver Braunling and Alexander Haberman for a careful reading on earlier drafts.
\section{Preliminaries on The Chow Weight Complex, Pro-objects, and K-theory}
\subsection{Chow Motives and the Gillet--Soul\'{e} Construction}
We first briefly review parts of the Gillet--Soul\'{e} construction of the weight complex of a variety.
For the rest of the paper, assume $k$ is a field that admits resolution of singularities. We will drop $k$ from the notation whenever possible. 
\begin{definition}\label{def:var}
Let $\mathcal{V}$ denote the category of varieties, or separated schemes of finite type over $k$, where the morphisms are $k$-morphisms. Let $\mathcal{V}^{prop}$ be the wide subcategory of varieties where we restrict our morphism sets to proper morphisms. Write $\V$ as the category of smooth projective varieties. 
\end{definition}

\begin{definition}
  Write $\Corr$ as the additive category of correspondences modulo numerical equivalence. The objects are smooth projective varieties, and morphisms are defined as follows:
\[\Hom_{\Corr}(X,Y):=\bigoplus_{i\in I}A^{\dim(Y_i)}(X\times Y_i)\]
where $I$ indexes the connected components of $Y$, and $A^k(-)$ denotes the $k$-th Chow group. The composition map
\[\circ:\Hom_{\Corr}(X,Y)\times \Hom_{\Corr}(Y,Z)\to \Hom_{\Corr}(X,Z)
\]
where $X,Y,W$ are connected, is generated by sending two cycles $\Gamma\in Z^{\dim(Y)}(X,Y)$ and $\Lambda\in Z^{\dim(W)}(Y,Z)$, to $(\Gamma\times  Z)\cap (X\times \Lambda)$ in $A^{\dim(Z)}(X\times Z)$ via the pushforward map induced by the projection $X\times Y\times Z\to X\times Z$, and where `$\cap$' denotes the cycle-theoretic intersection.
\end{definition}
\begin{definition}\label{def:Chowmotives}
$\Chow^{\eff}$, or the category of effective Chow motives, is then defined to be the category with objects:
$(X,p)$, where $X$ is a smooth projective variety and $p\in A^{\dim(X)}(X\times X)$ is a projector, i.e. $p^2=p$. The morphism sets are defined as:
\[\Hom_{\Chow^{\eff}}\Big((X,p),(Y,q)\Big):=\{f\in \Hom_{\Corr}(X,Y):f\circ p=q\circ f\}\]
For a further discussion of effective Chow motives, see \cite{GS}.
\end{definition}

\begin{remark}\label{rmk:2.4}
For a morphism $\alpha\in \Hom_{\Chow^{\eff}}(X,Y)$, let $\alpha^{T}\in \Hom_{\Chow^{\eff}}(Y,X)$ denote its image under the isomorphism induced by the swap $X\times Y\cong Y\times X$. With this notation, there is a contravariant functor $I:\V^{op}\hookrightarrow \Chow^{\eff}$ which is the identity on objects and sends a morphism $(f:X\to Y)$ to $[\Gamma_f]^{T}\in \Hom_{\Chow^{\eff}}(Y,X)$, where $\Gamma_f$ denotes the graph of $f$ in $X\times Y$, and the brackets denote the image in the Chow ring. 
\end{remark}
\begin{definition}
Let $\s\V$ denote the category of simplicial objects in $\V$
\end{definition}
\begin{definition}\label{def:phi}
We can use the functor $I$ in Remark \ref{rmk:2.4} to define a new functor
\[\Phi:\s\V\to \Ch(\Chow^{\eff})
\]
Which sends a simplicial variety $X_\bullet$ to the \emph{Moore complex} of $I(X_\bullet)$, i.e. the complex:
\[\dots\to I(X_{n+1})\xrightarrow[]{\sum_{i=0}^{n+1}(-1)^id_i^{n+1}} I(X_{n})\xrightarrow[]{\sum_{i=0}^{n}(-1)^id_i^n} I(X_{n-1})\to\dots
\]
\end{definition}
We are now in a position to review the Gillet--Soul\'{e} construction in \cite{GS}, where they build, for each variety $X$, a complex $M_\bullet(X)$ in the homotopy category of bounded chain complexes of effective Chow motives, $\Ch^{\mathrm{b}}(\Chow^{\eff})$. We start by defining a topology on the category schemes.
\begin{proposition}[1.4.1 of \cite{GS}]\label{prop:envelope}
    An \emph{envelope} is a proper morphism of schemes $f:X\to Y$ such that for all field extensions $K/k$, the induced map on $K$ points $f:X(K)\to Y(K)$ is surjective. These covers specify a Grothendieck topology on $\mathcal{V}$. 
\end{proposition}
\begin{definition}
    A \emph{hyperenvelope} is a morphism of simplicial schemes $f:X_\bullet\to Y_\bullet$ that is a hypercover in the envelope topology on $\mathcal{V}$. In other words, for all field extensions $K/k$, the induced morphism on $K$-points $X_\bullet(K)\to Y_\bullet(K)$ is a trivial Kan fibration. 
\end{definition}
\begin{convention}
 If $h_X:\tilde{X}_\bullet\to X_\bullet$ is a hyperenvelope where $\tilde{X}_n$ is smooth and projective for all $n$, then we call $h_X$ a \emph{smooth projective hyperenvelope}.
\end{convention}
If resolution of singularities is available, we have the following theorems at our disposal from \cite{GS}:
\begin{theorem}[Lemma 2 of \cite{GS}]\label{thm:hypenvexists}
Given any proper simplicial variety $X_\bullet$ (i.e. where each $X_n$ is proper), there exists a smooth projective hyperenvelope $\tilde{X}_\bullet$ of $X$.
\end{theorem}

\begin{corollary}[Lemma 2 of \cite{GS}]\label{cor:2.10}
Given any proper morphism $X_\bullet\to Y_\bullet$ between proper simplicial varieties, there exists smooth projective hyperenvelopes $h_X:\tilde{X}_\bullet\to X_\bullet$ and $h_Y:\tilde{Y}_\bullet\to Y_\bullet$ respectively and a map $\tilde{f}:\tilde{X}_\bullet\to\tilde{Y}_\bullet$ such that the following diagram of simplicial varieties commutes:
\[\begin{tikzcd}
\tilde{X}_\bullet\arrow[d,"h_X"]\arrow[r,"\tilde{f}"]&\tilde{Y}_\bullet\arrow[d,"h_Y"]\\
X_\bullet\arrow[r,"f"]&Y_\bullet
\end{tikzcd}\]
Further, the induced morphism of simplicial varieties:
\[\tilde{X}_\bullet\to X\times_{Y_\bullet}\tilde{Y}_\bullet
\]
is also a hyperenvelope.
\end{corollary}
\begin{convention}
  Let $s\mathcal{V}$ denote the category of simplicial varieties. Recall that $\mathrm{Ar}(s\mathcal{V})$ denotes the arrow category of $s\mathcal{V}$, where the objects of $\mathrm{Ar}(s\mathcal{V})$ are morphisms of $s\mathcal{V}$ and the morphisms of $\mathrm{Ar}(s\mathcal{V})$ are commutative squares. As constructed in Corollary \ref{cor:2.10}, $\tilde{f}$ is called a \emph{smooth projective hyperenvelope} of $f$ in $\mathrm{Ar}(s\mathcal{V})$. 
\end{convention}
\begin{proof}[Proof of Corollary \ref{cor:2.10}]
First, consider the hyperenvelope $X_\bullet\times_{Y_\bullet}\tilde{Y}_\bullet\to X_\bullet$ (the pullback of a hypercover is a hypercover). Then, take a smooth projective hyperenvelope $\tilde{X}_\bullet$ of $X_\bullet\times_{Y_\bullet}\tilde{Y}_\bullet$, which exists by Theorem \ref{thm:hypenvexists}. Compositions of hypercovers remain hypercovers, so we have that $\tilde{X}_\bullet$ is a hyperenvelope of $X_\bullet$, and the induced morphism $\tilde{f}:=\tilde{X}_\bullet\to \tilde{Y}_\bullet$ is the desired hyperenvelope of $f$.  
\end{proof}
The proof of the following theorem is similar to the previous one.
\begin{theorem}
Given a morphism of arrows $\Phi:f\to g$ in $\mathrm{Ar}(s\mathcal{V})$, where $f:X_\bullet\to Y_\bullet$, $g:Z_\bullet\to W_\bullet$, and the vertical morphisms are proper simplicial morphisms of proper simplicial varieties, there exist smooth projective hyperenvelopes $\tilde{f},\tilde{g}$ of $f$ and $g$ respectively and a morphism of arrows $\tilde{\Phi}:\tilde{f}\to \tilde{g}$ such that the following diagram commutes \[
\begin{tikzcd}
\tilde{f}\arrow[d]\arrow[r]&\tilde{g}\arrow[d]\\
f\arrow[r]&g
\end{tikzcd}
\]
Further, we require that the induced morphism $\tilde{\Phi}\to f\times_g\tilde{g}$ is a smooth projective hyperenvelope. We will call such $\tilde{\Phi}$ \emph{smooth projective hyperenvelopes} of $\Phi$.
\end{theorem}
Using these theorems, we can now construct the integral Chow motive of $X$. The procedure is as follows:
\begin{enumerate}
    \item Choose a proper variety $\o{X}$ such that there is an open embedding $X\to\o{X}$. $\o{X}$ is called a \emph{compactification} of $X$\begin{footnote}
        {A more standard use of the term `compactification' is when the embedding $X\to \o{X}$ is made to be dense. We do not require this.}
    \end{footnote}. 
    \item Choose a smooth projective hyperenvelope $\tilde{j}_X$ of the proper inclusion $j_X:\o{X}-X\to \o{X}$, viewed as a morphism between constant simplicial varieties.
    \item  Define $\textbf{M}(X)$ to be the isomorphism class of $Cone(\tilde{j}_X)\in \Ho(\Ch^{\perf}(Chow))$. Much of the hard work is in this step, i.e. showing that $\textbf{M}(X)$ is well-defined up to homotopy equivalence and is homotopic to a bounded chain complex, despite all the choices made.
\end{enumerate}
For closed embeddings $f:Z\to X$, choose compactifications $\o{Z}$ of $Z$, $\o{X}$ of $X$ and find smooth projective hyperenvelopes over the following diagram:
\[\begin{tikzcd}
\o{Z}&\o{\Gamma_f}\arrow[r]\arrow[l]&\o{X}\\
\o{Z}-Z\arrow[u,"j_Z"]&\o{\Gamma_f}-\Gamma_f\arrow[r]\arrow[l]\arrow[u,"j_\Gamma"]&\o{X}-X\arrow[u,"j_X"]
\end{tikzcd}
\]
where $\o{\Gamma_f}$ is the closure of the graph $\Gamma_f$ viewed in $\o{X}\times\o{Z}$, and the horizontal arrows are simply projections. In Section 2.3, p.141 of \cite{GS}, it is shown that any hyperenvelope $\tilde{\pi_X}:\tilde{j}_\Gamma\to \tilde{j}_X$ of the projection $\pi_X:j_\Gamma\to j_X$ induces a homotopy equivalence $Cone(\tilde{j}_\Gamma)\to Cone(\tilde{j}_X)$ in chain complexes of Chow motives, giving us a morphism $\textbf{M}(Z)\to\textbf{M}(X)$. Finally, for any open inclusion $f:U\to X$, choose a compactification $\o{X}-X$, observe that it is also a compactification of $U$, and find a smooth projective hyperenvelope of the proper inclusion $\o{X}-X\to \o{X}-U$ to give a contravariant morphism $\textbf{M}(X)\to\textbf{M}(U)$. We observe at the outset that the approaches towards assigning maps to closed embeddings and open immersions differ significantly, and it is not immediately clear how to make either of these assignments functorial.\\\\
\subsection{Waldhausen $K$-theory}
We now review some preliminaries about Waldhausen categories with functorial factorizations of weak cofibrations (or FFWC) and recall how to endow DG categories with a Waldhausen structure with FFWC. In the case of chain complexes in an idempotent complete exact category, we use work of \cite{Buh} to show that this type of category will have a more concrete Waldhausen with FFWC structure, which will be useful for later on.
\begin{definition}\label{def:wald}
    A Waldhausen category $\mathcal{W}$ is a category with a choice of a zero object, $0$, equipped with a subcategory of weak equivalences $\textbf{w}(\mathcal{W})$ and a subcategory of cofibrations $\textbf{co}(\mathcal{W})$ such that:
    \begin{enumerate}
        \item Isomorphisms are both included in $\textbf{w}(\mathcal{W})$ and $\textbf{co}(\mathcal{W})$
        \item For each object $A\in ob(\mathcal{W})$ the map $0\to A$ is a cofibration
        \item If $A\to B$ is a cofibration, and $A\to C$ is any morphism, the pushout $B\cup_AC$ exists, and the induced map $C\to B\cup_AC$ is also a cofibration
        \item Given the following diagram:
        \[\begin{tikzcd}
        A\arrow[d,"\simeq"] &B\arrow[r,hook]\arrow[l]\arrow[d,"\simeq"] &C\arrow[d,"\simeq"]\\
        A'&B'\arrow[r,hook]\arrow[l] &C'
        \end{tikzcd}\]
        where the right horizontal morphisms are cofibrations and the vertical morphisms are weak equivalences, the induced map $A\cup_BC\to A'\cup_{B'}C'$ is a weak equivalence. 
    \end{enumerate}
\end{definition}

\begin{definition}
In a Waldhausen category $\mathcal{W}$, a \emph{weak cofibration} is a morphism $f:X\to Y$ such that in $\mathrm{Ar}(\mathcal{W})$, there is a zig-zag of morphisms $f\to g_1\leftarrow g_2\to\dots\leftarrow g_n\to g$, where the legs of each morphism of arrows are weak equivalences in $\mathcal{W}$ and $g$ is a cofibration in $\mathcal{W}$.
\end{definition}
\begin{definition}
    For $n\in\N$, let $[n]$ denote the category associated to the poset $\{0<1<\dots<n\}$. 
\end{definition}
\begin{definition}
    For any category $\mathcal{W}$, define $\mathrm{Fun}([n], \mathcal{W})$ to be the category of functors from $[n]$ to $\mathcal{W}$.
\end{definition}
Just taking $n=1$, a functor $[1]\to\mathcal{W}$ picks out a morphism of $\mathcal{W}$. Further, a morphism in $\mathrm{Fun}([1], \mathcal{W})$ denotes a commutative square in $\mathcal{W}$. 
\begin{definition}\label{def:ffwc}
 Let $\mathcal{W}$ be a Waldhausen category. Let $\mathrm{Fun}^{wc}([1], \mathcal{W})$ be the full subcategory of $\mathrm{Fun}([1], \mathcal{W})$ consisting of those functors whose image is a weak cofibration. Write $\mathrm{Fun}^{c,w}([2], \mathcal{W})$ for the full subcategory of $\mathrm{Fun}([2], \mathcal{W})$ consisting of those diagrams which are a cofibration followed by a weak equivalence. A functorial factorization of weak cofibrations (FFWC) is a
functor \[T: \mathrm{Fun}^{wc}([1], \mathcal{W})\to \mathrm{Fun}^{c,w}([2], \mathcal{W})\] such that \[(d^1)^{*}\circ T= id_{\mathrm{Fun}^{wc}([1],C)}\]
Here $d^1:[1]\to [2]$ takes $0$ to $0$ and $1$ to $2$. More transparently, on objects $T$ is given by:
\[
(X\xrightarrow{f} Y)\mapsto \begin{tikzcd}
    X\arrow[r,hook]\arrow[dr,"f"]&Tf\arrow[d,"\simeq"]\\
    &Y
\end{tikzcd}
\]
\end{definition}
\begin{definition}[2.6 of \cite{BM}]\label{def:fmcwc}
    If $\mathcal{W}$ is a Waldhausen category that admits FFWC with functor $T$, it further admits functorial mapping cylinders for weak cofibrations (FMCWC) if there is a natural transformation $Y\to Tf$ splitting the natural weak equivalence $Tf\xrightarrow{\simeq}Y$, giving us the diagram:
    \[
    \begin{tikzcd}
        &Y\arrow[d]\arrow[dr,equal]&\\
       X\arrow[r,hook]&Tf\arrow[r,"\simeq"]&Y 
    \end{tikzcd}
    \]
\end{definition}
In the paragraph directly after Definition 2.6 of \cite{BM} the following lemma is proved:
\begin{lemma}
    If $\mathcal{W}$ is Waldhausen with FFWC such that all morphisms are weak cofibrations, then $\mathcal{W}$ also admits FMCWC.
\end{lemma}
In \cite{BM}, FMCWC is an important property because of the following theorem, stated for below our purposes. It will be important for the proof of Theorem \ref{thm:proobj}:
\begin{theorem}[5.5 of \cite{BM}]\label{lcof}
    Let $\mathcal{W}$ denote a saturated Waldhausen category (i.e.\ one whose weak equivalences satisfy the two-out-of-three property; this assumption is in force throughout \cite{BM}) with FMCWC, and with every morphism being a weak cofibration. If $\mathrm{\textbf{w}}(\mathcal{W})$ denotes the category of weak equivalences, then $(\mathcal{W}, \mathrm{\textbf{w}}(\mathcal{W}))$ admits a homotopy calculus of left fractions.\begin{footnote}{
        See Definition 5.3 of \cite{BM} for a review of this concept.}
    \end{footnote}
\end{theorem}

In \cite{BGN}, the $K$-theory of a DG category $\mathcal{A}$ is defined as the Waldhausen $K$-theory of a category $\perf_{\mathrm{co}}(\mathcal{A}^{op}\mathrm{-dgm})$ of cofibrant perfect right DG $\mathcal{A}$-modules, which we now introduce. Its cofibrations and weak equivalences come from a model structure discussed in Section 3.1 of \cite{To}; the reason we care that this model structure is \emph{cofibrantly generated} is that cofibrant generation is precisely what produces the functorial factorizations (FFWC) on $\perf_{\mathrm{co}}(\mathcal{A}^{op}\mathrm{-dgm})$, see Proposition \ref{prop:dgffwc} below.
\begin{proposition}[Proposition 1 of \cite{BGN}]\label{prop:2.19} 
Let $\mathcal{A}$ be a category enriched in $R\mathrm{-dgm}$ (the category of differential graded $R$ modules for any commutative ring $R$) and $\mathcal{A}^{op}\mathrm{-dgm}$ denote the category of right DG $\mathcal{A}$-modules, i.e. the category of $R\mathrm{-dgm}$-enriched functors $\mathrm{Fun}(\mathcal{A},R\mathrm{-dgm})$. There is a cofibrantly generated model structure on $\mathcal{A}^{op}\mathrm{-dgm}$, such that a morphism $F\to G$ of right DG modules is a weak equivalence (resp. a fibration) if and only if for every
$X\in \mathcal{A}$ the induced morphism of complexes
$F(X)\to G(X)$
is a weak equivalence (respectively a fibration) of differential graded chain complexes over the commutative ring $R$, where fibrations and weak equivalences are defined via the projective model structure on the category of chain complexes in $R\mathrm{-dgm}$.
\end{proposition}

\begin{definition}
An object of $\mathcal{A}^{op}\mathrm{-dgm}$ is called \emph{perfect} if it is a compact object of the derived category of DG $\mathcal{A}$-modules (equivalently, if it is homotopically finitely presented in the sense of Section 3.1 of \cite{To}). We write $\perf_{\mathrm{co}}(\mathcal{A}^{op}\mathrm{-dgm})$ for the full subcategory of $\mathcal{A}^{op}\mathrm{-dgm}$ spanned by the cofibrant perfect objects. The cofibrations and weak equivalences of the model structure of Proposition \ref{prop:2.19} endow $\perf_{\mathrm{co}}(\mathcal{A}^{op}\mathrm{-dgm})$ with a Waldhausen structure.
\end{definition}
\begin{definition}\label{def:kthydg}
The Waldhausen $K$-theory of a DG category $\mathcal{A}$ is defined to be the Waldhausen $K$-theory of cofibrant and perfect right DG $\mathcal{A}$-modules, i.e. $\perf_{\mathrm{co}}(\mathcal{A}^{op}\mathrm{-dgm})$.
\end{definition}
\begin{definition}[Before Remark 1 of \cite{BGN}]
Let $H^0(\mathcal{A})$ denote the homotopy category of the DG category $\mathcal{A}$. $\mathcal{A}$ is called triangulated if the Yoneda embedding $\mathcal{A}\hookrightarrow \mathcal{A}^{op}\mathrm{-dgm}$ induces an equivalence of homotopy categories:
\[H^0(\mathcal{A})\hookrightarrow \Ho(\perf_{\mathrm{co}}(\mathcal{A}^{op}\mathrm{-dgm}))
\]
Where $\Ho(\perf_{\mathrm{co}}(\mathcal{A}^{op}\mathrm{-dgm}))$ is the homotopy category given by the model structure on $\perf_{\mathrm{co}}(\mathcal{A}^{op}\mathrm{-dgm})$. 
\end{definition}
Note that every object of $\mathcal{A}$ is now equated with its representable DG-presheaf in $\mathcal{A}^{op}\mathrm{-dgm}$, which will certainly be cofibrant and perfect. The following proposition is proved in Lemma 1 of \cite{BGN}:
\begin{proposition}\label{prop:dgffwc}
     $\perf_{\mathrm{co}}(\mathcal{A}^{op}\mathrm{-dgm})$ has FFWC, where every morphism is a weak cofibration.
\end{proposition}
We now recall the setup for the Waldhausen structure on chain complexes of an idempotent complete additive category, in particular that of effective Chow motives. We use terminology from \cite{Buh}; see loc. cit. for a more detailed discussion.
\begin{definition}
 Let $\Ch^{\perf}(\mathcal{A})$ denote the full subcategory of chain complexes $A_\bullet$ 
 that are compact objects of the derived category $D(\mathcal{A})$ as constructed in 10.4 of \cite{Buh}. Notably, all chain complexes quasi-isomorphic to bounded complexes are objects in this category. We call this the subcategory of perfect chain complexes of $\mathcal{A}$.
\end{definition}
The next proposition follows quickly from the work done in \cite{Buh}.
\begin{proposition}\label{prop:2.27}
If $\mathcal{A}$ is idempotent complete, $\Ch^{\perf}(\mathcal{A})$ has a Waldhausen structure where:
\begin{enumerate}
    \item Cofibrations are the \emph{level-wise split inclusions}. By definition, these are the chain maps $f_\bullet:A_\bullet\to C_\bullet$ such that for each $n$, the morphism $f_n:A_n\to C_n$ is an inclusion onto a direct summand in $\mathcal{A}$, i.e.\ $C_n\cong A_n\oplus B_n$ with $f_n$ the natural inclusion. In this situation we write $C_\bullet=A_\bullet\oplus B_\bullet$; see Remark \ref{rem:oplus} for the abuse of notation involved.
    \item Weak equivalences are quasi-isomorphisms.
\end{enumerate}
\end{proposition}
\begin{remark}\label{rem:oplus}
    Note that cofibrations are not necessarily split as chain maps. For notational convenience, when we write $A_\bullet\oplus B_\bullet$, we mean only that the underlying \emph{graded} object of $\mathcal{A}$ splits as the coproduct of the graded objects underlying $A_\bullet$ and $B_\bullet$; the differential need not be the direct sum differential. Concretely, the differential of $A_\bullet\oplus B_\bullet$ has the upper triangular matrix form $\begin{bmatrix} d_A&\delta\\ 0&d_B\end{bmatrix}$, where $d_A$ is the differential of the subcomplex $A_\bullet$, $d_B$ is the induced differential on the quotient, and $\delta$ is some degree $-1$ component $B_\bullet\to A_\bullet$. 
\end{remark}
\begin{proof}
We verify each of the conditions in Definition \ref{def:wald}. It is immediate that isomorphisms are cofibrations, and since isomorphisms are chain homotopy equivalences they are also quasi-isomorphisms by 11.2 of \cite{Keller}, verifying condition 1. Condition 2 is satisfied immediately. For condition 3, if we take the pushout along a cofibration $A_\bullet\to A_\bullet\oplus B_\bullet$ diagram in $\Ch(\mathcal{A})$:
\[
\begin{tikzcd}
A_\bullet\arrow[r]\arrow[d]& A_\bullet\oplus B_\bullet\arrow[d]\\\
C_\bullet\arrow[r]&C_\bullet\oplus B_\bullet
\end{tikzcd}
\]
Here the pushout $C_\bullet\oplus B_\bullet$ is computed degree-wise, $(C_\bullet\oplus B_\bullet)_n=C_n\oplus B_n$, with differential $\begin{bmatrix} d_C&u\circ\delta\\ 0&d_B\end{bmatrix}$, where $u:A_\bullet\to C_\bullet$ is the left vertical map and $\delta:B\to A$ is the corresponding component of the differential of $A_\bullet\oplus B_\bullet$ as in Remark \ref{rem:oplus}. Now, if $A_\bullet\oplus B_\bullet$ and $C_\bullet$ are objects of $\Ch^{\perf}(\mathcal{A})$, then since compact objects in $D(\mathcal{A})$ form a strictly full subcategory (i.e.\ a full subcategory closed under isomorphisms), $B_\bullet$ and therefore the pushout must also be compact in $D(\mathcal{A})$ and therefore in $\Ch^{\perf}(\mathcal{A})$, showing us that pushouts along cofibrations exist in $\Ch^{\perf}(\mathcal{A})$. Since $C_\bullet\to C_\bullet\oplus B_\bullet$ is an inclusion onto a direct summand of graded objects of $\mathcal{A}$, we have satisfied condition 3. For condition 4, consider any diagram
\[
\begin{tikzcd}
A_\bullet\arrow[d,"f"]&B_\bullet\arrow[r,hook]\arrow[d,"g"]\arrow[l,"j"]&B_\bullet\oplus C_\bullet\arrow[d,"h"]\\
A_\bullet'& B_\bullet'\arrow[l,"j'"]\arrow[r,hook]& B_\bullet'\oplus C_\bullet'
\end{tikzcd}
\]
where the right horizontal maps are cofibrations and the vertical maps are weak equivalences. We need to show that the induced morphism on pushouts $\psi:A_\bullet\oplus C_\bullet\to A'_\bullet\oplus C_\bullet'$ is a weak equivalence. We can more explicitly write $h=\begin{bmatrix} g & h_{12}\\
0&h_{22}
\end{bmatrix}$. We can then consider the diagram:
\[
\begin{tikzcd}
A_\bullet\oplus B_\bullet\oplus C_\bullet\arrow[r,"\phi"]\arrow[d,"f\oplus h"]&A_\bullet\oplus C_\bullet\arrow[d,"\psi"]\\
A_\bullet'\oplus B_\bullet'\oplus C_\bullet'\arrow[r,"\phi'"]&A_\bullet'\oplus C_\bullet'
\end{tikzcd}
\]
where $\phi=\begin{bmatrix} 1&j&0\\
0&0&1
\end{bmatrix}$, $\phi'=\begin{bmatrix} 1&j'&0\\
0&0&1
\end{bmatrix}$, and $\psi=\begin{bmatrix} f& j'\circ h_{12}\\
0&h_{22}
\end{bmatrix}$ is the induced map on pushouts. This diagram commutes: the equality $\psi\circ\phi=\phi'\circ(f\oplus h)$ follows by direct computation from the commutativity $f\circ j=j'\circ g$ of the left-hand square in the previous diagram. Note that $f\oplus h$ is a quasi-isomorphism, as $f$ and $h$ are. The morphism $\phi$ is degree-wise split surjective, and since $j$ is a chain map, its kernel is the subcomplex
\[\{(-j(b),b,0)\ :\ b\in B_\bullet\}\cong B_\bullet,
\]
and similarly for $\phi'$. The vertical morphisms of the diagram above therefore assemble into a morphism of degree-wise split short exact sequences of complexes
\[\begin{tikzcd}[column sep=1.6em]
0\arrow[r]&B_\bullet\arrow[r]\arrow[d,"g"]&A_\bullet\oplus B_\bullet\oplus C_\bullet\arrow[r,"\phi"]\arrow[d,"f\oplus h"]&A_\bullet\oplus C_\bullet\arrow[r]\arrow[d,"\psi"]&0\\
0\arrow[r]&B_\bullet'\arrow[r]&A_\bullet'\oplus B_\bullet'\oplus C_\bullet'\arrow[r,"\phi'"]&A_\bullet'\oplus C_\bullet'\arrow[r]&0
\end{tikzcd}
\]
Since they are degree-wise split, the two rows define distinguished triangles in the homotopy category of complexes over $\mathcal{A}$, and the vertical morphisms define a morphism between these triangles. The cofibers of the vertical morphisms therefore fit into a distinguished triangle
\[C_g\to C_{f\oplus h}\to C_\psi\to C_g[1].
\]
Since $g$ and $f\oplus h$ are quasi-isomorphisms, $C_g$ and $C_{f\oplus h}$ are acyclic, and since the acyclic complexes form a triangulated subcategory of the homotopy category (Corollary 10.5 of \cite{Buh}), $C_\psi$ is acyclic as well. Hence $\psi$ is a quasi-isomorphism, i.e. a weak equivalence, concluding the proof.
\end{proof}
\begin{proposition}\label{prop:ffwc}
The Waldhausen structure on $\Ch^{\perf}(\mathcal{A})$ in Proposition \ref{prop:2.27} has FFWC. 
\end{proposition}
\begin{proof}
The first observation to make is that any morphism in $\Ch^{\perf}(\mathcal{A})$ is a weak cofibration, i.e. $\mathrm{Fun}^{wc}([1],\Ch^{\perf}(\mathcal{A}))=\mathrm{Fun}([1],\Ch^{\perf}(\mathcal{A}))$. Indeed, for a morphism $f_\bullet:X_\bullet\to Y_\bullet$ in $\Ch^{\perf}(\mathcal{A})$, we define a functor $\phi:\mathrm{Fun}([1],\Ch^{\perf}(\mathcal{A}))\to \mathrm{Fun}([2],\Ch^{\perf}(\mathcal{A}))$ defined on objects by
\[\phi(f_\bullet)=(X_\bullet\to M_{f_\bullet}\to Y_\bullet)
\]
where $M_{f_{\bullet}}$ is the homological mapping cylinder of $f_\bullet$. Further, for any morphism $\Psi$ in $\mathrm{Fun}([1],\Ch^{\perf}(\mathcal{A}))$ as shown below:
\[\begin{tikzcd}
X_\bullet\arrow[r,"f_\bullet"]\arrow[d]&Y_\bullet\arrow[d]\\
X'_\bullet\arrow[r,"g_\bullet"]&Y_\bullet' 
\end{tikzcd}\]
$\phi(\Psi):\phi(f_\bullet)\to \phi(g_\bullet)$ is defined to be the induced morphism on mapping cylinders:
\[\begin{tikzcd}
X_\bullet\arrow[r]\arrow[d]&M_{f_\bullet}\arrow[r]\arrow[d]&Y_\bullet\arrow[d]\\
X'_\bullet\arrow[r,]&M_{g_\bullet}\arrow[r]&Y_\bullet' 
\end{tikzcd}\]
Since the compositions of the upper and lower horizontal morphisms are $f$ and $g$ respectively, we can conclude that $(d^1)^*\phi=id_{\mathrm{Fun}^{wc}([1],\Ch^{\perf}(\mathcal{A}))}$.
\end{proof} 

We recall the notion of an exact/weakly exact functor from \cite{W} and \cite{BM} for future use:
\begin{definition}\label{def:exact}
A functor $F:\mathcal{W}\to\mathcal{V}$ between two Waldhausen categories $\mathcal{W},\mathcal{V}$ is (resp. weakly) \emph{exact} if:
\begin{itemize}
    \item It preserves (resp. weak) cofibrations
    \item It preserves weak equivalences 
    \item It sends pushouts along cofibrations to (resp. weak) pushouts along cofibrations
\end{itemize}
A weak pushout is simply a commutative square such that there is a zig-zag of weak equivalences from it to a pushout square.
\end{definition}
\begin{proposition}{\cite{W},\cite{BM}}
 An exact functor $F:\mathcal{W}\to\mathcal{V}$ between two Waldhausen categories $\mathcal{W},\mathcal{V}$ lifts to a map of $K$-theory spectra $K(F):K(\mathcal{W})\to K(\mathcal{V})$. The same result holds if both $\mathcal{W},\mathcal{V}$ have FFWC and $F$ is weakly exact.
\end{proposition}
We now state a key theorem of \cite{BM} that we will use repeatedly in the paper (we state a slightly weaker version for our purposes):
\begin{theorem}[1.3 of \cite{BM}]\label{thm:BM}
    Let $F:\mathcal{C}\to\mathcal{D}$ be a weakly exact functor between two saturated Waldhausen (i.e. whose weak equivalences satisfy the 2-out-of-3 property) categories with FFWC. Suppose that $F$ induces an equivalence on homotopy categories, and that $F(f)$ is a weak equivalence iff $f$ is a weak equivalence. then 
    \[
    K(F):K(\mathcal{C})\to K(\mathcal{D})
    \]
    is an equivalence of spectra.
\end{theorem}
\subsection{Pro-objects}
 We review some more background on pro-objects, including the straightening lemma for morphisms of pro-objects (this can be found in more detail both the appendix of \cite{AM} and Section 3 of \cite{tateobj}). The main purpose of this review is to give ways to concretely describe morphisms of pro-objects. 
\begin{definition}
Given a filtered category $I$ and a small category $\mathcal{C}$, a pro-object of $\mathcal{C}$ is a functor:
\[X:I^{op}\to \mathcal{C}
\]
\end{definition}
\begin{definition}
The collection of pro-objects of $\mathcal{C}$ form a category $\Pro(\mathcal{C})$, where the objects are pro-objects $X:I^{op}\to\mathcal{C}$ and the morphisms between pro-objects, say $X:I^{op}\to \mathcal{C}$ and $Y:J^{op}\to\mathcal{C}$, are of the form:
\[\Hom_{\Pro(\mathcal{C})}(X,Y):=\varprojlim_{j\in J}\Big(\varinjlim_{i\in I}\Hom(X(i), Y(j))\Big)
\]
In other words, a morphism $\Phi$ is given by compatible collections of morphisms:
\[\{[f_j]\in \varinjlim_{i\in I}\Hom(X_{i},Y_j)\}_{j\in ob(J)}
\]
\end{definition}
Here, compatible means that for any morphism $\phi:j'\to j$ and $f:X_i\to Y_j\in [f_j]$, we have $Y_{\phi}\circ f\in [f_{j'}]$. We can describe morphisms between pro-objects more generally via the straightening lemmas, a full discussion of which can be found in appendix A.3 of \cite{AM} or in Section 3 of \cite{tateobj}. 
\begin{definition}\label{def:rep}
    Given a morphism $\{[f_j]\in \varinjlim_{i\in I}Hom(X_i,Y_j)\}_{j\in ob(J)}$ from $X\to Y$, it is possible to represent this morphism as a pair $(\phi,\{f_j\}_{j\in ob(J)})$, where $\phi: Ob(J)\to Ob(I)$ is a function and $f_j:X_{\phi(j)}\to Y_j$ is a representative of $[f_j]$. We call $(\phi,\{f_j\}_{j\in ob(J)})$ a \emph{representation of $\phi$ indexed by $Y$}
\end{definition}
\begin{definition}[1.5 of the Appendix of \cite{AM}]\label{def:cofinal}
Let $I$ be filtered. A functor $\Phi:I\to J$ is \emph{cofinal} if:
\begin{enumerate}
    \item For every $j\in ob(J)$ there is an $i\in ob(I)$ and a morphism $j\to \phi(i)$
    \item If $j\in ob(J)$ and $i\in ob(I)$ and if $f,g:j\to \phi(i)$ are two morphisms in $J$, then there exists $i'$ and a morphism $h:i\to i''$ such that $\phi(h)\circ f=\phi(h)\circ g$.
\end{enumerate}
Note that $\phi$ being cofinal implies that $J$ is also filtered, and that cofinality is transitive.
\end{definition}
\begin{remark}\label{rem:cofinalrestrict}
The relevance of cofinality is the following (see 1.5--1.8 of the Appendix of \cite{AM}): if $\Phi:I\to J$ is cofinal and $X:J^{op}\to\mathcal{C}$ is a pro-object, then the restriction $X\circ\Phi^{op}:I^{op}\to\mathcal{C}$ is isomorphic to $X$ as an object of $\Pro(\mathcal{C})$. We use this observation repeatedly in Section 5: it allows us to show that various concretely constructed indexing categories (such as the categories $A_f$ and $\Phi_f$ of Proposition \ref{prop:covar} and Theorem \ref{prop:contra}) represent morphisms between the pro-objects $P_X$.
\end{remark}
\begin{definition}[Before 3.1 of the Appendix of \cite{AM}]\label{def:downarrow}
Let $f:X\to Y$ be a morphism of pro-objects, where $X:I^{op}\to\mathcal{C}$ and $Y:J^{op}\to\mathcal{C}$. We say that a single morphism $\alpha_{ij}:X_i\to Y_j$ of $\mathcal{C}$ \emph{represents} $f$ if $\alpha_{ij}$ belongs to the class $[f_j]\in\varinjlim_{i\in I}\Hom(X_i,Y_j)$ determined by $f$; equivalently, if $\alpha_{ij}$ appears in some representation of $f$ in the sense of Definition \ref{def:rep}. Let $I\downarrow_{\mathcal{C}} J$ denote the category with objects the triples $\{(i,j,\alpha_{ij}:X_i\to Y_j)\}$ where $\alpha_{ij}$ represents $f$, and where a morphism $(i,j,\alpha_{ij})\to (i',j',\alpha_{i'j'})$ corresponds to a commutative square of the form:
\[
\begin{tikzcd}
X_i\arrow[d]\arrow[r,"\alpha_{ij}"]&Y_j\arrow[d]\\
X_{i'}\arrow[r,"\alpha_{i'j'}"]&Y_{j'}
\end{tikzcd}
\]
where the vertical morphisms are transition morphisms of the respective pro-objects. The purpose of this construction is to package a morphism of pro-objects as a pro-object in the arrow category $\mathrm{Ar}(\mathcal{C})$: the assignment $(i,j,\alpha_{ij})\mapsto \alpha_{ij}$ defines a functor $(I\downarrow_{\mathcal{C}} J)^{op}\to \mathrm{Ar}(\mathcal{C})$, i.e.\ a pro-object of arrows, which records $f$ together with all of its representatives.

More generally, let $\mathcal{J}$ be a finite category and let $D:\mathcal{J}\to \Pro(\mathcal{C})$ be a diagram of pro-objects, say with $D(x):I_x^{op}\to\mathcal{C}$ for each object $x$ of $\mathcal{J}$. We say that a functor $F:\mathcal{J}\to \mathcal{C}$ \emph{represents} $D$ if for every object $x$ of $\mathcal{J}$ we have $F(x)=D(x)(i_x)$ for some index $i_x\in I_x$, and for every morphism $e:x\to y$ of $\mathcal{J}$, the morphism $F(e):F(x)\to F(y)$ represents the morphism $D(e):D(x)\to D(y)$ in the sense above. 
\end{definition} 
The following result, a form of the \emph{straightening lemma} of \cite{AM} (see A.3.3 of loc.\ cit.), allows one to replace a diagram of pro-objects by a pro-object of diagrams. We will use it in this form in Sections 4 and 5:
\begin{proposition}[3.3 of the Appendix of \cite{AM}]\label{prop:straightening}
Let $\mathcal{J}$ denote a finite diagram with commutation relations, and suppose that $\mathcal{J}$ has no loops, i.e., that the beginning and end of a chain of arrows are always distinct. Consider a functor $D:\mathcal{J}\to \Pro(\mathcal{C})$ and the collection of diagrams representing $D:\mathcal{J}\to \Pro(\mathcal{C})$:
\[
S_{D}:=\{F:\mathcal{J}\to \mathcal{C}\ |\ \text{$F$ represents $D$}\}
\]
$S_{D}$ is cofiltered, and if $i\in I$ indexes the objects of $\mathcal{J}$ and $D_i$ denotes the image of $i$ in $\Pro(\mathcal{C})$, then the projection maps $S_{D}^{op}\to D_i^{op}$ are cofinal. So, in this situation, one can represent a diagram of pro-objects by a pro-object of diagrams.
\end{proposition}

\section{A Waldhausen Category of Simplicial Smooth Projective Varieties}
Recall that $\V$ is the category of smooth projective varieties and $\s\V$ denotes its simplicial objects. We first review some constructions in $\s\V$ which will be useful later on.

\begin{definition}\label{def:Deltan}
    Let $\Delta^n:\Delta^{op}\to \s\V$ denote the simplicial variety given on objects by:
    \[[\ell]\mapsto \coprod_{0\leq i_{0}\leq \dots \leq i_{\ell}\leq n}*\]
    where $*:=\Spec(k)$. We label the copy of $*$ associated to $i_{0}\leq \dots \leq i_{\ell}$ as $[i_{0}, \dots, i_{\ell}]$. The face and degeneracy maps are the same deletion and insertion maps of $\Delta^n$ as in the usual setting of simplicial sets.
\end{definition}
Note that the face and degeneracy maps are simply coproducts of identity maps, making them morphisms of varieties.
\begin{remark}
    If $\mathrm{FinSet}$ denotes the category of finite sets, and $\mathrm{FinSet}\hookrightarrow \V$ denotes the embedding sending a finite set $X\mapsto \coprod_{X}\Spec(k)$, then $\Delta^n$ in Definition \ref{def:Deltan} is just the image of the usual $\Delta^n$ via the induced functor on simplicial objects $s\mathrm{FinSet}\hookrightarrow \s\V$.
\end{remark}
\begin{definition}
    Let $\{0\}$ denote the simplicial subvariety of $\Delta^n$ defined by the elements of the form $[0\dots 0]$ in $\Delta^n$, and likewise, let $\{1\}$ denote the simplicial subvariety of $\Delta^n$ defined by elements of the form $[1\dots 1]$. Note that these are both constant simplicial sets. 
\end{definition}

\begin{definition}
   Given a morphism of simplicial varieties $f_\bullet: X_\bullet\to Y_\bullet$, we define the simplicial mapping cylinder to be the pushout:
    \[sM_{f_{\bullet}}:=(X_\bullet\times \Delta^1)\cup^{i_1,f_\bullet} (Y_\bullet)
    \]
    where $i_1: X_\bullet\hookrightarrow X_\bullet\times\Delta^1$ is the inclusion map onto $X_\bullet\times\{1\}$.
    \end{definition}
    For a morphism  $f_\bullet:X_\bullet\hookrightarrow Y_\bullet$ such that for all $n$, the map $f_n$ is an inclusion of the form $f_n:X_n\to X_n\coprod X_n'$, pushouts along $f_\bullet$ exist in $\s\V$. So, $sM_{f_\bullet}$ exists in $\s\V$.
    The following proposition is well-documented in the literature, as the proof is the same as in simplicial sets.
    \begin{proposition}\label{prop:sM}
        Any morphism $f_\bullet:X_\bullet\to Y_\bullet$ in $\s\V$ factors as:
        \[X_\bullet\xrightarrow{i_0} sM_{f_\bullet}\xrightarrow{\simeq} Y_\bullet \]
        where the first simplicial map is a level-wise inclusion of $X_\bullet$ onto the direct summand $X_\bullet\times\{0\}$, and the second map is a simplicial homotopy equivalence.
    \end{proposition} 

    \begin{definition}\label{def:cone}
        We define the simplicial mapping cone of a morphism $f_\bullet: X_\bullet\to Y_\bullet$ as the pushout:
        \[\begin{tikzcd}
        X_\bullet\arrow[d]\arrow[r,"i_0"]&sM_{f_\bullet}\arrow[d]\\
            *\arrow[r]& sC_{f_\bullet}
        \end{tikzcd}\]
        
    \end{definition}
    In order to put a Waldhausen structure on the category of simplicial varieties, we must first make them \emph{well-pointed}, which we now discuss. 
\begin{definition}
    We call the undercategory $*/\s\V$ the category of \emph{pointed} simplicial smooth projective varieties, i.e. whose objects are morphisms $(*\to X_\bullet)$ from the constant simplicial set $*$ to a simplicial variety $X_\bullet$, and morphisms between two objects, say from $(*\to X_\bullet)$ to $(*\to Y_\bullet)$, are diagrams of the form:
    \[
    \begin{tikzcd}
        &*\arrow[dr]\arrow[dl]&\\
        X_\bullet\arrow[rr,"f_\bullet"]&&Y_\bullet
    \end{tikzcd}
    \]
    i.e. simplicial maps $f_\bullet:X_\bullet\to Y_\bullet$ that respect *.
\end{definition}
\begin{definition}\label{def:sVperf}
    Let $\s\V_*$ denote the full subcategory of $*/\s\V$ consisting the pointed objects such that for each $n$, the pointed morphism is of the form:
    \[*\hookrightarrow X_n\coprod *\]
    i.e. the pointed morphism is an inclusion onto a direct summand $*\hookrightarrow X_\bullet\coprod *$, viewed as graded schemes (i.e. forgetting the degeneracy and face maps). We call $\s\V_*$ the category of \emph{well-pointed} simplicial smooth projective varieties. As a slight abuse of notation, we write $X_\bullet\coprod *$ to denote a well-pointed object going forward unless necessary.  
\end{definition}
 The following corollary justifies that we can replace $*/\s\V$ with $\s\V_*$:
\begin{corollary}
Write $\mathrm{w.e.}$ for the class of pointed morphisms whose underlying morphism of simplicial varieties is a simplicial homotopy equivalence; we emphasize that neither the homotopy inverse nor the homotopies are required to respect the base points. The inclusion $\s\V_*\hookrightarrow */\s\V$ induces an equivalence of relative categories
\[
(\s\V_*,\mathrm{w.e.})\simeq (*/\s\V,\mathrm{w.e.})
\]
in the sense that it is fully faithful, and essentially weakly surjective: every object in $*/\s\V$ is connected to an object of $\s\V_*$ by a pointed morphism whose underlying morphism of simplicial varieties is a simplicial homotopy equivalence.
\end{corollary}
\begin{proof}
    We just need to show that the inclusion is essentially weakly surjective. For an object $(f_\bullet:*\to X_\bullet)$ of $*/\s\V$, since the map $i_0:*\to sM_{f_\bullet}$ is an inclusion onto a direct summand by Proposition \ref{prop:sM}, the map defines a well-pointed object of $*/\s\V$. Further, Proposition \ref{prop:sM} states that the map $sM_{f_\bullet}\to X_\bullet$ is a simplicial homotopy equivalence of underlying simplicial varieties, and it is a pointed morphism: the base point $*\times\{0\}$ of $sM_{f_\bullet}$ is sent to the base point $f_\bullet(*)$ of $X_\bullet$. Note that the homotopy inverse need not be pointed.
\end{proof}
 In working with $\s\V_*$, we have also not lost any of the information in $\s\V$, as the following proposition is immediate:
\begin{proposition}
    There is faithful embedding:
    $\s\V\hookrightarrow \s\V_*$ given on objects by:
    \[X_\bullet\mapsto X_\bullet\coprod *
    \]
    and on morphisms, $f_\bullet:X_\bullet\to Y_\bullet$ gets sent to $f_\bullet\coprod id_*:X_\bullet\coprod *\to Y_\bullet\coprod *$.
\end{proposition}
We now put a Waldhausen structure on $\s\V_*$:
\begin{definition}\label{def:cof}
    In $\s\V_*$, let $\textbf{co}(\s\V_*)$ denote the subcategory with the same objects as $\s\V_*$ but whose morphisms $f_\bullet:X_\bullet\to Y_\bullet$ are level-wise isomorphisms onto a direct summand, i.e. for each $n\in\N$: \[f_n:X_n\to X_n\coprod Y_n\] 
    So, $f_\bullet$, viewed as just a morphism of graded schemes, is an inclusion onto a direct summand 
\end{definition}
\begin{remark}[Notational Remark]\label{rem:notation}
    From now on, when we write $A_\bullet\to A_\bullet\coprod B_\bullet$, we mean an inclusion on to the direct summand of graded schemes (and not simplicial schemes), i.e a cofibration in $\s\V_*$. 
\end{remark}
\begin{remark}\label{rem:pushouts}
Given any diagram:
\[A_\bullet\leftarrow B_\bullet\to B_\bullet\coprod C_\bullet
\]
where the right morphism is a cofibration, then the pushout \[A_\bullet\cup_{B_\bullet}(B_\bullet\coprod C_\bullet)\]
is well-defined, as it is done level-wise in $\s\V_*$. At the level of graded schemes, the pushout is simply $A_\bullet\coprod C_\bullet$. Further, the induced map $A_\bullet\to A_\bullet\coprod C_\bullet$ is again a cofibration in $\s\V_*$, as it is an inclusion onto a direct summand of graded schemes. 
\end{remark}
\begin{definition}\label{def:we1}
In $\s\V_*$, let $w(\s\V_*)$ denote the smallest subcategory containing the following morphisms:
\begin{enumerate}
    \item Pointed morphisms whose underlying morphism in $\s\V$ is a simplicial homotopy equivalence (as before, the homotopy inverse and the homotopies are not required to be pointed), 
    \item Morphisms of the form:
\[f_\bullet\coprod id_*:X_\bullet\coprod *\to Y_\bullet\coprod *\] such that $f_\bullet:X_\bullet\to Y_\bullet$ is a hyperenvelope,
\item Extensions of (1) and (2), i.e. given a diagram:
\[\begin{tikzcd}
    X_\bullet\arrow[d,"\simeq"]\arrow[r,"f_\bullet"]&Y_\bullet\arrow[d,"\simeq"]\\
    X'_\bullet\arrow[r,"f_\bullet'"]&Y_\bullet'
\end{tikzcd}\]
where the vertical morphisms (decorated with `$\simeq$' in the diagram) are compositions of morphisms of type (1) and (2), then we define the induced map on mapping cones $sC_{f_\bullet}\to sC_{f_\bullet'}$ to also be a member of $w(\s\V_*)$. We stress that the decoration `$\simeq$' here refers to these compositions; the subcategory of weak equivalences is only being defined in the present definition and in Definition \ref{def:we}.
\item If we have a diagram \[
    \begin{tikzcd}[row sep=2.2em, column sep=2.2em]
    X_\bullet \arrow[rr,"f_\bullet"] \arrow[dr] \arrow[dd,"h_X"] &&
    Y_\bullet \arrow[dd,"h_Y",swap] \arrow[dr] \\
    & X_\bullet'\arrow[rr,"f_\bullet' "]\arrow[dd,"h_{X'}",swap] &&
    Y_\bullet' \arrow[dd,"h_{Y'}",swap] \\
     X \arrow[rr] \arrow[dr] && Y\arrow[dr] \\
    & X' \arrow[rr] &&Y'
    \end{tikzcd}
\]
$h_X,h_Y,h_{X'},h_{Y'}$ are hyperenvelopes and the bottom square is an abstract blowup square of proper varieties over $k$, i.e. a Cartesian square where $X'\to Y'$ is a closed embedding and $Y\to Y'$ is a proper map, then we also define the induced map on mapping cones $sC_{f_\bullet}\to sC_{f_\bullet'}$ to also be a weak equivalence. 
\end{enumerate} 

\end{definition}
\begin{definition}\label{def:we}
Let $\textbf{w}(\s\V_*)$ denote the smallest subcategory containing $w(\s\V_*)$ such that:
\begin{enumerate}
    \item Given any diagram:
\[\begin{tikzcd}
A_\bullet\arrow[d,"\simeq"] &B_\bullet\arrow[r,hook]\arrow[l]\arrow[d,"\simeq"] & B_\bullet\coprod C_\bullet\arrow[d,"\simeq"]\\
A'_\bullet &B'_\bullet\arrow[r,hook]\arrow[l] &B_\bullet'\coprod C'_\bullet
\end{tikzcd}\]
as in axiom 4 of Definition \ref{def:wald} where the vertical morphisms are in $\textbf{w}(\s\V_*)$ and the horizontal morphisms are cofibrations as in Definition \ref{def:cof}, then the induced morphism on pushouts:
\[A_\bullet\coprod C_\bullet\to A'_\bullet\coprod C'_\bullet 
\]
is in $\textbf{w}(\s\V_*)$        
\item Given morphisms $g:X_\bullet\to Y_\bullet$ and $f:Y_\bullet\to Z_\bullet$ such that any two of the morphisms: $f,g,f\circ g$ are in $\textbf{w}(\s\V_*)$, then so is the third.
\end{enumerate}
\end{definition}
\begin{remark}
It should be noted that $\emph{\textbf{w}}(\s\V_*)$ can be explicitly constructed as follows:
Let $\{\mathcal{C}_i\subset \s\V_*\}_{i\in\Lambda}$ denote the collection of subcategories of $\s\V_*$ containing $w(\s\V_*)$ and satisfying properties (1) and (2) in Definition \ref{def:we}. Note that $\s\V_*$ is contained in this set, so it is non-empty. Then, we can define
\[\emph{\textbf{w}}(\s\V_*):=\bigcap_{i\in\Lambda}\mathcal{C}_i
\]
more precisely, the objects of $\emph{\textbf{w}}(\s\V_*)$ are the same as $\s\V_*$ (as $w(\s\V_*)$ is a wide subcategory) and the morphisms will are defined as the intersection of morphism sets in $\mathcal{C}_i$. One can readily prove that $\emph{\textbf{w}}(\s\V_*)$ in fact a category, and it will also satisfy both the axioms of Definition \ref{def:we}. by virtue of all of its morphisms living in each of the $\mathcal{C}_i$. By construction, it will be the minimal choice.
\end{remark}
\begin{proposition}\label{prop:3.17}
$\Big(\s\V_*,\emph{\textbf{w}}(\s\V_*),\textbf{co}(\s\V_*)\Big)$ is a Waldhausen category.
\end{proposition}
\begin{proof}
First, observe that $*$ is a zero object in $\s\V_*$, and we denote it $0$. We now explain why each of the axioms in Definition \ref{def:wald} is satisfied.
\begin{enumerate}
    \item Since isomorphisms are pointed morphisms whose underlying morphisms are simplicial homotopy equivalences, they are contained in $\textbf{w}(\s\V_*)$. Isomorphisms are also cofibrations by definition.
    \item Since the morphism from the zero object is a level-wise inclusion onto a direct summand, it is a cofibration.
    \item Axiom 3 is satisfied by Remark \ref{rem:pushouts}.
    \item Axiom 4 is automatically satisfied due to (1) of Definition \ref{def:we}.
    \end{enumerate}
\end{proof}
    
    \begin{proposition}\label{prop:3.18}
        If $f_\bullet:X_\bullet\coprod *\to Y_\bullet\coprod *$ is a morphism of well-pointed objects, then $sM_{f_\bullet}$ defines a well-pointed object such that the factorization
        \[X_\bullet\coprod *\xrightarrow{i_0} sM_{f_\bullet}\xrightarrow{\simeq} Y_\bullet \coprod *\]
        holds in $\s\V_*$, i.e. through morphisms of well-pointed objects. Here $i_0$ denotes the inclusion onto the direct summand $(X_\bullet\coprod *)\times\{0\}$ from Proposition \ref{prop:sM}, and the second map is a pointed morphism whose underlying morphism is a simplicial homotopy equivalence, i.e. a weak equivalence of type (1) in Definition \ref{def:we1}.
    \end{proposition}
    \begin{proof}
        Since $i_{0}:X_\bullet\coprod *\to sM_{f_\bullet}$ is an isomorphism onto a direct summand of graded simplicial schemes (forgetting the simplicial maps), it follows that the restriction to $*$:
        \[i_{0}|_*:*\to sM_{f_\bullet}\]
        defines a well-pointed object, promoting $i_{0}$ to a morphism of well-pointed objects, and a cofibration in $\s\V_*$. Since $f_\bullet$ respects the inclusions of $*$ into both $X_\bullet\coprod *$ and $Y_\bullet\coprod *$, the map: \[sM_{f_\bullet}\xrightarrow{\simeq}  Y_\bullet\coprod *\] respects the inclusions of $*$ into $sM_{f_\bullet}$ and $Y_\bullet\coprod *$ as well. It is therefore a pointed morphism whose underlying morphism of simplicial varieties is a simplicial homotopy equivalence (the homotopy inverse supplied by Proposition \ref{prop:sM} need not be pointed), and thus a weak equivalence in $\s\V_*$ of type (1) in Definition \ref{def:we1}.
    \end{proof}
\begin{proposition}\label{prop:sVffwcAll}
    $\s\V_*$ has FFWC.
\end{proposition}

\begin{proof}
The first observation to make is that any morphism in $\s\V_*$ is a weak cofibration, i.e. $\mathrm{Fun}^{wc}([1],\s\V_*)=\mathrm{Fun}([1],\s\V_*)$. Indeed, for a morphism $f_\bullet:X_\bullet\coprod *\to Y_\bullet\coprod *$ in $\s\V_*$, we have the commutative diagram in $\s\V_*$:
\begin{equation*}
\begin{tikzcd}
X_\bullet\coprod *\arrow[rr,"f_\bullet"]&&Y_\bullet\coprod *\\
X_\bullet\coprod *\arrow[r,hook]\arrow[u,equal] &sM_{f_\bullet}\arrow[r,"\simeq"]&Y_\bullet\coprod *\arrow[u,equal]
\end{tikzcd}
\end{equation*}
where in the bottom row, the first map is a cofibration and the second is a weak equivalence in $\s\V_*$ by Proposition \ref{prop:3.18}. The rest of the proof is the same as in Proposition \ref{prop:ffwc}, but with $sM_{f_\bullet}$ instead of $M_{f_\bullet}$.
\end{proof}
\begin{remark}\label{rem:cofcone}
If $f_\bullet:X_\bullet\to Y_\bullet$ is a cofibration, then the simplicial mapping cone of Definition \ref{def:cone} is weakly equivalent to the pushout $Y_\bullet\cup_{X_\bullet}*$. Indeed, applying the gluing axiom (1) of Definition \ref{def:we} to the diagram
\[\begin{tikzcd}
*\arrow[d,equal]&X_\bullet\arrow[l]\arrow[r,hook,"i_0"]\arrow[d,equal]&sM_{f_\bullet}\arrow[d,"\simeq"]\\
*&X_\bullet\arrow[l]\arrow[r,hook,"f_\bullet"]&Y_\bullet
\end{tikzcd}\]
whose right vertical morphism is the simplicial homotopy equivalence of Proposition \ref{prop:sM}, shows that the induced morphism on pushouts
\[sC_{f_\bullet}\to Y_\bullet\cup_{X_\bullet}*
\]
lies in $\textbf{w}(\s\V_*)$. We stress that this is a weak equivalence and not in general an isomorphism: $sC_{f_\bullet}$ is a quotient of the mapping cylinder $sM_{f_\bullet}$ rather than of $Y_\bullet$. Since the subcategories introduced in Definition \ref{def:sVb} are closed under isomorphism in $\Ho(\s\V_*)$, this distinction will not matter for our purposes.
\end{remark}
\begin{definition}\label{def:sVb}
    Let $\Ho(\s\V_*):=\s\V_*[\textbf{w}(\s\V_*)^{-1}]$ denote the homotopy category of $\s\V_*$ obtained by inverting weak equivalences. Let $\s\V_*^{\mathrm{b}}$ denote the smallest full subcategory of $\s\V_*$ which contains the finite coproducts and products of constant representable simplicial objects, and which satisfies two-out-of-three for simplicial cofiber sequences: whenever $f_\bullet:A_\bullet\to B_\bullet$ is a morphism of $\s\V_*$ and two of the three objects $A_\bullet$, $B_\bullet$, $sC_{f_\bullet}$ lie in $\s\V_*^{\mathrm{b}}$, then so does the third. In particular, $\s\V_*^{\mathrm{b}}$ is closed under finite iterations of simplicial mapping cones. Lastly, let $\s\V_*^{\perf}$ denote the smallest full subcategory of $\s\V_*$ which contains $\s\V_*^{\mathrm{b}}$, which contains every object isomorphic in $\Ho(\s\V_*)$ to one of its objects, and which again satisfies two-out-of-three for simplicial cofiber sequences. We call simplicial objects in $\s\V_*^{\perf}$ \emph{perfect} objects of $\s\V_*$.
\end{definition}
\begin{remark}\label{rem:desusp}
    The two-out-of-three condition is imposed by hand because $\Ho(\s\V_*)$ is not triangulated; in particular, we have no desuspension with which to recover the source of a morphism from its target and its simplicial mapping cone. In \cite{GS} the corresponding argument is carried out in $\Ch^{\perf}(\Chow^{\eff})$, where the analogous closure property is automatic since a thick triangulated subcategory satisfies two-out-of-three for distinguished triangles. We use this condition in the proof of Theorem \ref{prop:boundedweight}.
\end{remark}
\begin{corollary}\label{prop:ffwcsV}
The subcategories of cofibrations and weak equivalences in $\s\V_*$ from Proposition \ref{prop:3.17} restricted to $\s\V_*^{\perf}$ makes $\s\V_*^{\perf}$ Waldhausen with FFWC.
\end{corollary}
\begin{proof}
Notice that coproducts of weak equivalences are also weak equivalences by (1) of Definition \ref{def:we}, so pushouts along cofibrations will be contained in $\s\V_*^{\perf}$: given a cofibration $A_\bullet\to B_\bullet$ and a morphism $A_\bullet\to C_\bullet$ between perfect objects, write $P_\bullet:=C_\bullet\cup_{A_\bullet}B_\bullet$ for the pushout, so that $C_\bullet\to P_\bullet$ is again a cofibration. Pasting the two pushout squares
\[\begin{tikzcd}
A_\bullet\arrow[r]\arrow[d]&B_\bullet\arrow[d]\\
C_\bullet\arrow[r]\arrow[d]&P_\bullet\arrow[d]\\
*\arrow[r]&P_\bullet\cup_{C_\bullet}*
\end{tikzcd}\]
exhibits the outer rectangle as a pushout as well, giving a canonical isomorphism $P_\bullet\cup_{C_\bullet}*\cong B_\bullet\cup_{A_\bullet}*$. By Remark \ref{rem:cofcone} these two pushouts are weakly equivalent to $sC_{(C_\bullet\to P_\bullet)}$ and $sC_{(A_\bullet\to B_\bullet)}$ respectively, so $sC_{(C_\bullet\to P_\bullet)}$ and $sC_{(A_\bullet\to B_\bullet)}$ are isomorphic in $\Ho(\s\V_*)$. Since $A_\bullet$ and $B_\bullet$ are perfect, so is $sC_{(A_\bullet\to B_\bullet)}$, hence so is $sC_{(C_\bullet\to P_\bullet)}$, and two-out-of-three as in Definition \ref{def:sVb}, applied to the cofibration $C_\bullet\to P_\bullet$, shows that $P_\bullet$ is perfect. Containing pushouts along cofibrations also implies that $\s\V_*^{\perf}$ will contain all of its mapping cylinders, meaning that the proof of $\s\V_*^{\perf}$ having FFWC is the same as in the proof of Proposition \ref{prop:sVffwcAll}.
\end{proof}
A first observation is that one can use the Gillet and Soul\'{e} approach described in the previous section to construct a weight complex for a variety in $\s\V_*$.
\begin{definition}\label{def:premotive}
Let $X$ be an arbitrary variety, i.e. any object of $\mathcal{V}$; we emphasize that $X$ is not assumed to be smooth or projective. A Gillet--Soul\'{e} \emph{pre-motive} of $X$ in $\s\V_*$ is a simplicial mapping cone of the form $sC_{\tilde{j}_X}$, where $\tilde{j}_{X}$ is a smooth projective hyperenvelope of the proper inclusion $j_X:\o{X}-X\to \o{X}$, for $\o{X}$ a compactification of $X$. Here $j_X$ is viewed as a morphism of constant simplicial varieties; since it is a proper morphism of proper varieties, a smooth projective hyperenvelope of $j_X$ exists by Theorem \ref{thm:hypenvexists} and Corollary \ref{cor:2.10}.
\end{definition}
The following proposition is a re-wording of Theorem 2 of \cite{GS} for our purposes, and we adapt the proof to our setting:
\begin{proposition}\label{prop:3.23}
All Gillet--Soul\'{e} pre-motives of a variety $X$ in $\s\V_*$ are isomorphic in the homotopy category of $\s\V_*$.
\end{proposition}
\begin{proof}  
Consider any two pre-motives $sC_{\tilde{j}_{X,1}},sC_{\tilde{j}_{X,2}}$ of $X$. Recall that we need to choose compactifications $\o{X}_1,\o{X}_2$ and smooth projective hyperenvelopes $\tilde{j}_{X,i}:\widetilde{\o{X}_i-X}_\bullet\to \tilde{X}_{i,\bullet}$ to construct $sC_{\tilde{j}_{X,1}}$ and $sC_{\tilde{j}_{X,2}}$ respectively. We can consider the compactification $\o{X}:=\o{X}_1\times\o{X}_2$ of $X$ where the open immersion $X\to\o{X}$ is the diagonal map. This gives us the elementary blow-up squares for $i=1,2$:
\[\begin{tikzcd}
    \o{X}_f-X\arrow[r,"j_f"]\arrow[d]&\o{X}_f\arrow[d]\\
    \o{X}_i-X\arrow[r,"j_{X{,}i}"]&\o{X}_i
\end{tikzcd}\]
where $\o{X}_f$ is the closure of the graph of the diagonal embedding in $\o{X}_1\times\o{X}_2$. Finding a smooth projective hyperenvelope of $j_f$ (which exists by Corollary \ref{cor:2.10}) that fits into the diagram:
\[
\begin{tikzcd}
\widetilde{X_{f}-X}_\bullet\arrow[r,"\tilde{j}_f"]\arrow[d]&\tilde{X}_{f,\bullet}\arrow[d]\\
\widetilde{\o{X}_i-X}_\bullet\arrow[r,"\tilde{j}_{X,i}"] &\tilde{X}_{i,\bullet}
\end{tikzcd}
\]
we see that by (4) of Definition \ref{def:we1}, the induced maps on simplicial mapping cones $sC_{\tilde{j}_f}\to sC_{\tilde{j}_{X,i}}$ are weak equivalences in $\s\V_*$. This shows that we can find a zig-zag of weak equivalences between any two weight complexes in $\s\V_*$, concluding the proof.
\end{proof}
We now also reprove Theorem 2(i) of \cite{GS} for our setting:
\begin{theorem}\label{prop:boundedweight}
For any variety $X$, every weight complex of $X$ lies in $\s\V_*^{\perf}$. 
\end{theorem}
\begin{proof}
By Proposition \ref{prop:3.23}, it suffices to show that we can construct one weight complex that is in $\s\V_*^{\mathrm{b}}$. Given a variety $X$, we let $U\subset X$ be a smooth dense open subset and $Z$ denote its complement. The proof of Theorem 2.3 of \cite{GS} shows that we can choose weight complexes $sC_{\tilde{j}_U}, sC_{\tilde{j}_Z}, sC_{\tilde{j}_X}$ of $U,Z,X$ that give a pushout square where the top row is a cofibration:
 \[
 \begin{tikzcd}
sC_{\tilde{j}_Z}\arrow[d]\arrow[r]&sC_{\tilde{j}_X}\arrow[d]\\
     0\arrow[r]&sC_{\tilde{j}_U}
 \end{tikzcd}
 \]
Since the square is a pushout along a cofibration, $sC_{\tilde{j}_U}$ is weakly equivalent to the simplicial mapping cone of $sC_{\tilde{j}_Z}\to sC_{\tilde{j}_X}$ by Remark \ref{rem:cofcone}. So, assuming we know that the proposition holds true for $U$ and $Z$, two-out-of-three for simplicial cofiber sequences (Definition \ref{def:sVb}) applied to the two outer terms $sC_{\tilde{j}_Z}$ and $sC_{\tilde{j}_U}$ gives us that it holds for the middle term $sC_{\tilde{j}_X}$, i.e. for $X$.
By Noetherian induction, we can assume that $X$ is smooth and quasi-projective. We can then take a smooth projective compactification $\o{X}$ of $X$ such that $\o{X}-X=D^{red}$, where $D$ is a divisor with normal crossings in $\o{X}$. Now, we can choose $\o{X}$ to be its own weight complex if we view $\o{X}$ as a constant in $\s\V_*$, so the proposition holds for $\o{X}$. Since $\dim(\o{X}-X)<\dim(X)$, we have by the inductive hypothesis that the proposition holds for $\o{X}-X$ as well. Appealing again to the proof of Theorem 2.3 of \cite{GS}, we can find a pushout diagram whose top row is again a cofibration:
\[
 \begin{tikzcd}
sC_{\tilde{j}_X}\arrow[d]\arrow[r]&sC_{\tilde{j}_{\o{X}}}\arrow[d]\\
     0\arrow[r]&sC_{\tilde{j}_{\o{X}-X}}
 \end{tikzcd}
 \]
 for weight complexes of $X,\o{X}$, and $\o{X}-X$. As before, this square exhibits $sC_{\tilde{j}_{\o{X}-X}}$ as being weakly equivalent to the simplicial mapping cone of $sC_{\tilde{j}_X}\to sC_{\tilde{j}_{\o{X}}}$, so a second application of two-out-of-three, this time recovering the source from the target $sC_{\tilde{j}_{\o{X}}}$ and the cone $sC_{\tilde{j}_{\o{X}-X}}$, proves the statement for $X$. So, we can conclude the result for all $X$.  
\end{proof}
We are now in a position to construct the map of K-theory spectra from $\s\V_*^{\perf}$ to $\Ch^{\perf}(\Chow^{\eff})$. 
\begin{proposition}
$\Chow^{\eff}$ is an idempotent complete additive category, so by Proposition \ref{prop:2.27}, $\Ch^{\perf}(\Chow^{\eff})$ has a Waldhausen structure with FFWC.
\end{proposition}
\begin{proposition}
The functor $\Phi$ in Definition \ref{def:phi} induces a functor
\[\Phi_*:\s\V_*^{\perf}\to \Ch^{\perf}(Chow^{\eff})\] which is in fact an exact functor of Waldhausen categories. 
\end{proposition}
\begin{proof}
We first show that the restriction $\Phi_*:\s\V_*^{\perf}\to \Ch(Chow^{\eff})$ preserves finite colimits and is an exact functor of Waldhausen categories by verifying the axioms of exactness in the order listed in Definition \ref{def:exact}. Then, since $\Phi_*$ must preserve weak equivalences, it will follow that the image of $\Phi_*$ will land in $\Ch^{\perf}(Chow^{\eff})$. Note that the map $\V\to Chow^{\eff}$ in Remark \ref{rmk:2.4} preserves finite colimits and monomorphisms (and so inclusions onto direct summands), which implies that $\Phi_*$ preserves cofibrations. Further, a pushout along a monomorphism will remain a pushout along a monomorphism, i.e. pushouts along cofibrations are preserved. As for weak equivalences, note that $\Phi_*$ sends simplicial homotopy equivalences to chain homotopy equivalences, which are in turn quasi-isomorphisms by 10.7 of \cite{Buh}. Further, hyperenvelopes are sent to chain homotopy equivalences in $\Ch(Chow)$ by Proposition 2 of \cite{GS}. Next, observe for a morphism $f:X_\bullet\to Y_\bullet$, $\Phi_*(sC_f)$ is the un-normalized homological mapping cone, which is homotopy equivalent to usual homological mapping cone $C_{\Phi_*(f)}$, So, if we have a diagram:
\[\begin{tikzcd}
X_\bullet\arrow[r,"f"]\arrow[d]&Y_\bullet\arrow[d]\\
X_\bullet'\arrow[r,"g"]&Y_\bullet'
\end{tikzcd}
\]
in $\s\V_*$, where the vertical morphisms are compositions of chain homotopy equivalences and hyperenvelope maps, then this diagram will get sent to one in $\Ch(Chow^{\eff})$ where the vertical maps are chain homotopies, implying that the induced morphism on mapping cones $C_{\Phi_*(f)}\to C_{\Phi_*(g)}$ is also a chain homotopy. This gives us a diagram:
\[
\begin{tikzcd}
\Phi_*(sC_f)\arrow[r,"\simeq"]\arrow[d]&C_{\Phi_*(f)}\arrow[d,"\simeq"]\\
\Phi_*(sC_g)\arrow[r,"\simeq"]&C_{\Phi_*(g)}
\end{tikzcd}
\]
where the horizontal maps are chain homotopies as well, by the same proof of 2.5.3 of \cite{GoJar}. Since chain homotopy is an equivalence relation, we can conclude that $\Phi_*(sC_f)\to \Phi_*(sC_g)$ must also be a weak equivalence. Lastly, for any diagram as in (4) of Definition \ref{def:we1}: \[
    \begin{tikzcd}[row sep=2.2em, column sep=2.2em]
    X_\bullet \arrow[rr,"f_\bullet"] \arrow[dr] \arrow[dd,"h_X"] &&
    Y_\bullet \arrow[dd,"h_Y",swap] \arrow[dr] \\
    & X_\bullet'\arrow[rr,"f_\bullet' "]\arrow[dd,"h_{X'}",swap] &&
    Y_\bullet' \arrow[dd,"h_{Y'}",swap] \\
     X \arrow[rr] \arrow[dr] && Y\arrow[dr] \\
    & X' \arrow[rr] &&Y'
    \end{tikzcd}
\]
$h_X,h_Y,h_{X'},h_{Y'}$ are hyperenvelopes and the bottom square is an abstract blowup square of proper varieties over $k$, the induced map $sC_{f_\bullet}\to sC_{f_\bullet'}$ is sent to a chain homotopy in $\Ch(Chow)$ by the discussion in Section 2.3 of \cite{GS}. So, we see that all morphisms in $w(\s\V_*^{\perf})$ are sent to chain homotopies in $\Ch^{\perf}(Chow)$. Now, observe that the subcategory of $\s\V_*^{\perf}$ whose morphisms are sent to chain homotopy equivalences in $\Ch^{\perf}(Chow)$ will be a category containing $w(\s\V_*^{\perf})$ and satisfying properties (1) and (2) of Definition \ref{def:we} as pushouts along cofibrations are preserved and chain homotopy is an equivalence relation, $\textbf{w}(\s\V_*^{\perf})$ will be sent to the subcategory of chain homotopies in $\Ch^{\perf}(Chow)$, proving that weak equivalences are preserved.
\end{proof}   

\begin{definition}\label{def:Chow}
    We write the induced map of spectra in the previous proposition as:
    \[K(\Phi_*):K(\s\V_*^{\perf})\to K(\Ch^{\perf}(Chow))\]
\end{definition}
\section{Pro-Objects of Waldhausen Categories}
In order to get a map of spectra $K(\mathcal{V})\to K(\s\V_*^\perf)$, we will first need to construct a $K$-theory spectrum out of a suitable subcategory of pro-objects of $\s\V_*^\perf$; this spectrum will be homotopy equivalent to $K(\s\V_*^\perf)$. Further, we can show that their Dwyer-Kan hammock localizations are equivalent, where localizations are taken with respect to the subcategory of weak equivalences. These equivalences are instances of a more general phenomenon, which we outline in this section.
\begin{hypothesis}\label{def:conditionspro}
    For the remainder of this section, we let $(\mathcal{C},\textbf{w}(\mathcal{C}),\textbf{co}(\mathcal{C}))$ denote a \emph{saturated} Waldhausen category (i.e.\ one whose weak equivalences satisfy the two-out-of-three property; this is a standing assumption in \cite{BM}, whose results we use throughout) with FFWC (see Definition \ref{def:ffwc} for a definition) where every morphism in $\mathcal{C}$ is a weak cofibration, and for every object $c\in\mathcal{C}$, the associated undercategory has a Waldhausen structure: \[\Big((c\downarrow \mathcal{C},\textbf{w}(\mathcal{C})\cap (c\downarrow \mathcal{C}), \textbf{co}(\mathcal{C})\cap (c\downarrow \mathcal{C})\Big) \]
    which also has FFWC.
\end{hypothesis}

\begin{definition}\label{def:prosV}
Let $\Pro^{wc}(\mathcal{C})$ denote the subcategory of $\Pro(\mathcal{C})$ consisting of pro-objects $\{X:I\to \mathcal{C}\}$ such that all morphisms within $X$ are isomorphisms in the homotopy category $\Ho(\mathcal{C})$. We refer to these as \emph{weakly constant pro-objects}.
\end{definition}
We now use the Waldhausen structure on $\mathcal{C}$ to define one on $\Pro^{wc}(\mathcal{C})$. 
\begin{proposition}\label{prop:prowald}
$\Pro^{wc}(\mathcal{C})$ is a Waldhausen category defined such that a morphism $f:X\to Y$ is a weak equivalence if one can find a representation of $f$ indexed by $Y$ (see Definition \ref{def:rep}) by weak equivalences in $\mathcal{C}$, and where $\textbf{co}(\Pro^{wc}(\mathcal{C}))$ is the smallest subcategory of $\Pro^{wc}(\mathcal{C})$ containing every isomorphism of $\Pro^{wc}(\mathcal{C})$ together with every morphism admitting a representation indexed by its target by cofibrations of $\mathcal{C}$.
\end{proposition}
\begin{proof}
We verify the axioms of Definition \ref{def:wald}. For the first axiom, let $f:X\to Y$ be an isomorphism of pro-objects, where $X:I^{op}\to \mathcal{C}$ and $Y:J^{op}\to \mathcal{C}$. It is a cofibration because $\textbf{co}(\Pro^{wc}(\mathcal{C}))$ contains all isomorphisms by definition. We claim that $f$ is also a weak equivalence, i.e. that \emph{every} representation of $f$ consists of weak equivalences of $\mathcal{C}$. Pick any representation of $f$ indexed by $Y$, i.e. $\{f_j:X_{\phi(j)}\to Y_j\}_{j\in ob(J)}$, and any representation of $f^{-1}$ indexed by $X$, i.e. $\{f^{-1}_i:Y_{\psi(i)}\to X_i\}_{i\in ob(I)}$. Then
\[\{f_{i}^{-1}\circ f_{\psi(i)}:X_{\phi(\psi(i))}\to X_i\}_{i\in ob(I)}
\]
is a representation of $f^{-1}\circ f=\mathrm{id}_X$ indexed by $X$. Unwinding what this means: for each $i\in ob(I)$, the class of $f_i^{-1}\circ f_{\psi(i)}$ in $\varinjlim_{i'}\Hom(X_{i'},X_i)$ coincides with that of $\mathrm{id}_{X_i}$, so there exist an index $i'$ and transition morphisms $t:X_{i'}\to X_{\phi(\psi(i))}$ and $t':X_{i'}\to X_i$ of the pro-object $X$ with $f_i^{-1}\circ f_{\psi(i)}\circ t=t'$. Since $X$ lies in $\Pro^{wc}(\mathcal{C})$, its transition morphisms $t$ and $t'$ are weak equivalences, so two-out-of-three, available since $\mathcal{C}$ is saturated by Hypothesis \ref{def:conditionspro}, shows that $f_i^{-1}\circ f_{\psi(i)}$ is a weak equivalence. The symmetric argument applied to $f\circ f^{-1}=\mathrm{id}_Y$ shows that $f_j\circ f^{-1}_{\phi(j)}$ is a weak equivalence. Hence $[f_j]$ admits a left and a right inverse in $\Ho(\mathcal{C})$ and is therefore an isomorphism, so $f_j$ is a weak equivalence, again by saturation. This proves the first axiom. For the second axiom, note that the $0$ object of $\Pro^{wc}(\mathcal{C})$ is the constant $0$ object of $\mathcal{C}$, and any morphism from $0\to (X:I^{op}\to \mathcal{C})$ in $\Pro^{wc}(\mathcal{C})$ can only be written in the form $\{0\to X_i\}_{i\in ob(I)}$, which are all cofibrations in $\mathcal{C}$, showing that $0\to X$ is a cofibration in $\Pro^{wc}(\mathcal{C})$. For axiom 3, consider the diagram of pro-objects:
\[
\begin{tikzcd}
X\arrow[r,"f"]\arrow[d,"g"]&Y\\
Z
\end{tikzcd}
\]
where the horizontal morphism is a cofibration, and the indexing categories of $X,Y,Z$ are $I,J,K$ respectively. Assume first that we can find a compatible collection \[\{f_{j}:X_{\phi(j)}\to Y_{j}\}_{j\in ob(J)}\] cofibrations of $\mathcal{C}$ representing $f$, and we can pick any representation of $g$, $\{g_k:X_{\psi(k)}\to Z_k\}_{k\in \mathrm{ob}(K)}$. The induced morphism on the pushout is represented by the morphisms:
\[\{Z_k\to Z_k\cup_{X_{\psi(k)}}Y_j\}_{(k,j)\in ob(K\times J)}
\]
which are all cofibrations in $\mathcal{C}$. So, we have a pushout diagram:
\[
\begin{tikzcd}
X\arrow[r,"f"]\arrow[d,"g"]&Y\arrow[d]\\
Z\arrow[r]&Z\cup_XY
\end{tikzcd}
\]
where the bottom morphism is a cofibration in $\Pro^{wc}(\mathcal{C})$. Assume instead that $f$ is an isomorphism. Then the pushout $Z\cup_XY$ exists and the induced morphism $Z\to Z\cup_XY$ is an isomorphism, hence a cofibration. Now, a general cofibration is by definition a composite $f=f_n\circ\dots\circ f_1$ in which each $f_m$ is either an isomorphism or admits a representation by cofibrations of $\mathcal{C}$, so we induct on $n$. By the pasting law for pushouts, the pushout of $f$ along $g$ is computed by pushing out one $f_m$ at a time, and at each stage one of the two cases above applies, showing that the pushout exists and that the induced morphism is again a cofibration; composing these gives axiom 3. Axiom 4 follows similarly. 
\end{proof}
\begin{remark}\label{rem:wc}
Every pro-object $X:I^{op}\to \mathcal{C}$ in $\Pro^{wc}(\mathcal{C})$ is weakly equivalent, for any object $i$ of $I$, to the constant pro-object $X(i)$, as the identity morphism $\{id_{X(i)}:X(i)\to X(i)\}$ is a representation of a morphism $X\to X(i)$ indexed by $X(i)$. 
\end{remark}
\begin{proposition}\label{prop:proffwc}
The Waldhausen structure on $\Pro^{wc}(\mathcal{C})$ admits FFWC with every morphism being a weak cofibration.
\end{proposition}
\begin{proof}
Since $\mathcal{C}$ has FFWC where every morphism is a weak cofibration, there exists a factorization functor
$\phi:\mathrm{Fun}([1],\mathcal{C})\to \mathrm{Fun}^{c,w}([2],\mathcal{C})$ with properties outlined in Definition \ref{def:ffwc}. For shorthand notation, we write:
\[
\phi(f:X\to Y)=X\to T_f\to Y
\]
It is immediate that $\phi$ extends to a functor on pro-objects of each category:
\[\Pro(\phi):\Pro(\mathrm{Fun}([1],\mathcal{C}))\to \Pro(\mathrm{Fun}^{c,w}([2],\mathcal{C}))
\]
 For any morphism $X\to Y$ between weakly constant pro-objects $X:I^{op}\to \mathcal{C},\ Y:J^{op}\to \mathcal{C}$, we choose a representation $\{\phi,f_j:X_{\phi(j)}\to Y_j\}_{j\in J}$ indexed by $Y$, and define $\Pro(\phi)$ to send it to the composition $\{X_{\phi(j)}\to T_{f_j}\to Y_j\}_{j\in J}$, where the first morphisms are cofibrations, and the second morphisms are weak equivalences. Now, since $T:\mathrm{Ar}(\mathcal{C})\to \mathcal{C}$ defines a functor, it also sends pro-objects of morphisms of $\mathcal{C}$ to pro-objects of $\mathcal{C}$, making $\{T_{f_j}\}_{j\in J}$ a pro-object indexed by $J$. As a result, we have that $\{X_{\phi(j)}\to T_{f_j}\to Y_j\}_{j\in J}$ represents a composition of morphisms $X\to \Pro(T)_f \to Y$, where $\Pro(T)_f$ is the pro-object $\{T_{f_j}\}_{j\in J}$. The first morphism is represented by the inclusions $\{\phi,X_{\phi(j)}\to T_{f_j}\}_{j\in J}$, i.e. a collection of cofibrations indexed by $\{T_{f_j}\}$ and so it is a cofibration in $\Pro^{wc}(\mathcal{C})$. The second morphism is represented by the weak equivalences $\{id,T_{f_j}\to Y_j\}_{j\in J}$ indexed by $Y$, making it a weak equivalence in $\Pro^{wc}(\mathcal{C})$. So, $\Pro(\phi)$ restricts to a functor:
 \[
 \Pro(\phi):Fun([1],\Pro^{wc}(\mathcal{C}))\to Fun^{c,w}([2],\Pro^{wc}(\mathcal{C}))
 \] 
 The rest of the proof, i.e. that $\Pro(\phi)$ does in fact satisfy the necessary conditions to give us a functorial factorization of weak cofibrations, follows in a straightforward way from the fact that $\phi$ satisfies the required conditions. Since we defined $\Pro(\phi)$ on all morphisms of $\Pro^{wc}(\mathcal{C})$, we see that all morphisms are weak cofibrations. 
\end{proof}
For what follows, we recommend Section 2 of \cite{DK} for a discussion of the Dwyer-Kan hammock localization. We provide a quick definition here for the convenience of the reader:
\begin{definition}[2.1 of \cite{DK}]\label{def:DKloc}
    Let $(\mathcal{C},\mathcal{W})$ denote a category with weak equivalences. The Dwyer-Kan hammock localization, $L^H\mathcal{C}$, is a simplicial category (and hence a model for an $\infty$-category) defined as follows: for every $X,Y\in ob(\mathcal{C})$
    the $k$-simplices of $L^H\mathcal{C}(X,Y)$ are commutative diagrams of the form:
    \[
    \begin{tikzcd}
        &C_{0,1}\arrow[r]\arrow[ddl]\arrow[d]&C_{0,2}\arrow[d]\arrow[r]\arrow[l]&\arrow[l]\dots\arrow[r]&C_{0,n}\arrow[ddr]\arrow[l]\arrow[d]&\\
        &C_{1,1}\arrow[r]\arrow[dl]\arrow[d]&C_{1,2}\arrow[d]\arrow[r]\arrow[l]&\arrow[l]\dots\arrow[r]&C_{1,n}\arrow[dr]\arrow[l]\arrow[d]&\\
        X\arrow[uur]\arrow[dr]\arrow[ur]\arrow[r]&\vdots\arrow[r]\arrow[l]\arrow[d]&\vdots\arrow[r]\arrow[d]\arrow[l]&\arrow[l]\vdots\arrow[r]&\vdots\arrow[l]\arrow[d]&Y\arrow[uul]\arrow[ul]\arrow[l]\arrow[dl]\\
        &C_{k,1}\arrow[r]\arrow[ul]&C_{k,2}\arrow[r]\arrow[l]&\arrow[l]\dots\arrow[r]&C_{k,n}\arrow[ur]\arrow[l]&
    \end{tikzcd}
    \]
    where $n$ can vary, all vertical maps are in $\mathcal{W}$. In each column, all maps go in the same direction, and if they go left, then they are in $\mathcal{W}$. Further, no columns contain all identity maps, and maps in adjacent columns go in different directions. Face and degeneracy maps are the usual deletion and insertion maps applied to the \emph{rows} of the hammock, i.e. deleting the row $(C_{i,1},\dots,C_{i,n})$ or repeating it; no such operations are performed on columns.
\end{definition}

\begin{theorem}\label{thm:DKequiv}
    Let $\mathcal{C}$ be a Waldhausen category satisfying Hypothesis \ref{def:conditionspro}. The inclusion $\mathcal{C}\hookrightarrow \Pro^{wc}(\mathcal{C})$ induces an equivalence of simplicial categories:
    \[
    L^H(\mathcal{C})\to L^H(\Pro^{wc}(\mathcal{C}))
    \]
    where $L^H\mathcal{C}$ denotes the hammock localization of a Waldhausen category $\mathcal{C}$ with respect to its weak equivalences. 
\end{theorem}
\begin{proof}
    We aim to invoke Theorem 1.4 of \cite{BM}, which states that if \[\Ho(\mathcal{C})\hookrightarrow \Ho(\Pro^{wc}(\mathcal{C}))\] is equivalence of categories, and for each $X\in \mathrm{ob}(\mathcal{C})$, the induced map on undercategories:
    \[
    \Ho(X\downarrow \mathcal{C})\to \Ho(X\downarrow \Pro^{wc}(\mathcal{C}))
    \]
    is an equivalence, then $L^H(\mathcal{C})\to L^H(\Pro^{wc}(\mathcal{C}))$ is an equivalence of simplicial categories.\\
    We first show that the functor $\Ho(\mathcal{C})\hookrightarrow \Ho(\Pro^{wc}(\mathcal{C}))$ is fully faithful and essentially surjective. Essential surjectivity follows from Remark \ref{rem:wc}. For fully faithfulness, we invoke Lemma 3.4.1 of \cite{Krause}. Observe that $\Big(\Pro^{wc}(\mathcal{C}),\textbf{w}(\Pro^{wc}(\mathcal{C})\Big)$ admits a homotopy calculus of left fractions by Theorem \ref{lcof}. Further, $\mathcal{C}$ is a full subcategory of $\Pro^{wc}(\mathcal{C})$ such that 
\[
\textbf{w}(\mathcal{C})=\textbf{w}(\Pro^{wc}(\mathcal{C}))\cap \mathcal{C}
\]
 Further, for any weak equivalence \[\sigma:X\to Y^{\Pro}\] where $X$ is a constant pro-object (but $Y^{\Pro}$ may not be constant), we have a representation $\{X\to Y_j\}_{j\in ob(J)}$ where each $X\to Y_j$ is a weak equivalence in $\mathcal{C}$. Picking any $j\in ob(J)$, Remark \ref{rem:wc} gives a weak equivalence $\tau:Y^{\Pro}\xrightarrow{\simeq}Y_j$, and so we have
\[
\tau\circ \sigma:X\xrightarrow{\simeq}Y^{\Pro}\xrightarrow{\simeq} Y_j
\]
which is a weak equivalence in $\mathcal{C}$. These conditions allow us to invoke lemma 3.4.1 of \cite{Krause} which concludes that \[\Ho(\mathcal{C})\hookrightarrow \Ho(\Pro^{wc}(\mathcal{C}))\] is fully faithful. So, we have a fully faithful and essentially surjective inclusion, giving us an equivalence of homotopy categories.\\
We now show the second condition needed to invoke Theorem 1.4 of \cite{BM}. Observe that if we show for any $X\in \mathcal{C}$, that $X\downarrow \Pro^{wc}(\mathcal{C})$ admits a homotopy calculus of left fractions, we can use the exact same argument as in the last paragraph (i.e. use Lemma 3.4.1 of \cite{Krause}). Now, due to Propositions \ref{prop:prowald} and \ref{prop:proffwc} along with the assumption on the undercategories of $\mathcal{C}$ in Hypothesis \ref{def:conditionspro}, $X\downarrow \Pro^{wc}(\mathcal{C})$ admits a Waldhausen with FFWC structure (where every morphism is a weak cofibration) inherited from $X\downarrow \mathcal{C}$. So, by Theorem \ref{lcof}, $X\downarrow \Pro^{wc}(\mathcal{C})$ admits a homotopy calculus of left fractions and we can conclude. 
\end{proof}

We can now state and prove that the resulting K-theory spectrum of $\Pro^{wc}(\mathcal{C})$ is equivalent to that of $\mathcal{C}$.
\begin{theorem}\label{thm:proobj}
Let $\mathcal{C}$ be a Waldhausen category satisfying Hypothesis \ref{def:conditionspro}. The inclusion $i:\mathcal{C}\hookrightarrow \Pro^{wc}(\mathcal{C})$ induces an equivalence:
\[K(i):K(\mathcal{C})\xrightarrow{\simeq} K(\Pro^{wc}(\mathcal{C}))
\]
\end{theorem}
\begin{proof}
We verify that the functor $i$ satisfies the conditions of Theorem \ref{thm:BM} to conclude the result. First note that by definition of their respective Waldhausen structures, the inclusion is an exact functor of Waldhausen categories. A cofibration of $\mathcal{C}$ is sent to a morphism of constant pro-objects represented by that same cofibration, which is one of the generators of $\textbf{co}(\Pro^{wc}(\mathcal{C}))$ in Proposition \ref{prop:prowald}. Further, $i(f)$ is a weak equivalence if and only if $f$ is a weak equivalence. Both categories are saturated, as Theorem \ref{thm:BM} requires: $\textbf{w}(\mathcal{C})$ is saturated by Hypothesis \ref{def:conditionspro}, and $\textbf{w}(\Pro^{wc}(\mathcal{C}))$ inherits the two-out-of-three property from $\textbf{w}(\mathcal{C})$ by choosing compatible representations of the morphisms involved.
The fact that the inclusion induces an equivalence of homotopy categories follows from Theorem \ref{thm:DKequiv}, so Theorem \ref{thm:BM} gives us that $K(i):K(\mathcal{C})\to K(\Pro^{wc}(\mathcal{C}))$ is in fact an equivalence.
\end{proof}
\begin{corollary}
    $\s\V_*^{\perf}$ satisfies Hypothesis \ref{def:conditionspro}, implying that the inclusion $\s\V_*^{\perf}\hookrightarrow \Pro^{wc}(\s\V_*^{\perf})$ induces the equivalences of their hammock localizations with respect to weak equivalences:
\[
L^H(\s\V_*^{\perf})\to L^H(\Pro^{wc}(\s\V_*^{\perf}))
\]
along with maps of their $K$-theory spectra
\[
K(\s\V_*^{\perf})\xrightarrow{\simeq}K(\Pro^{wc}(\s\V_*^{\perf}))
\]
\end{corollary}
\begin{proof}
The fact that $\s\V_*^\perf$ is Waldhausen with FFWC and that every morphism is a weak cofibration is proven in Corollary \ref{prop:ffwcsV}. Moreover, $\textbf{w}(\s\V_*^{\perf})$ is saturated: property (2) of Definition \ref{def:we} is precisely the two-out-of-three property. We must simply check that for all $X\in \mathrm{ob}(\s\V_*^\perf)$, the induced Waldhausen structure on $X\downarrow \s\V^{\perf}_*$ has FFWC. This is true due to the properties of the simplicial mapping cylinder: For any morphism $f:Y\to Z$ under $X$, the simplicial mapping cylinder factorization gives us a diagram:
\[
\begin{tikzcd}
    &X\arrow[d]\arrow[dr]\arrow[dl]&\\
    Y\arrow[r,hook]&sM_{f}\arrow[r,"\simeq"]&Z
\end{tikzcd}
\]
In other words, $sM_f$ is in the undercategory of $X$ for any such $f:Y\to Z$. So in this case, the functorial factorization on $\s\V_*^{\perf}$ induces one on the undercategory $X\downarrow \s\V_*^{\perf}$.
\end{proof}
\section{The Category of Pro-Weight Complexes}

We re-frame the above approach using pro-objects in order to lift the Gillet--Soul\'{e} measure. The goal of this section is to construct the category of \emph{pro-weight complexes}, or $PH(\mathcal{V})$ to encode all of the choices that are made in constructing the integral Chow motive. 

\begin{definition}\label{def:Vhyp}
Let $\mathcal{V}^{hyp}$ denote the category where the objects are triples:
\[(X,\o{X},\tilde{j}_X)
\]
where $X$ is a variety over $k$, i.e. an object of $\mathcal{V}$, and $\o{X}=\o{X}'\coprod *$ is a disjoint union of $\o{X}'$, a compactification of $X$, with a copy of $\Spec(k)$\begin{footnote}
    {The disjoint base-point is used heavily in the proof of Theorem \ref{prop:contra}.}
\end{footnote}. The map $\tilde{j}_X:\widetilde{\o{X}-X}_\bullet\to \tilde{X}_\bullet$ is a smooth projective hyperenvelope of the map $j_X:\o{X}-X\to \o{X}$. Morphisms are of the form $(\o{f},\tilde{f}):(X,\o{X},\tilde{j}_X)\to (Y,\o{Y},\tilde{j}_Y)$, where $\o{f}:\o{X}\to\o{Y}$ is a proper morphism such that $\o{f}(\o{X}-X)\subset \o{Y}-Y$, and $\tilde{f}:\tilde{j}_X\to\tilde{j}_Y$ is a hyperenvelope of $\o{f}|_{\o{X}-X}$, that is we have the following diagram: 
\[\begin{tikzcd}
\tilde{j}_X\arrow[r,"\tilde{f}"]\arrow[d]&\tilde{j}_Y\arrow[d]\\
j_X\arrow[r,"\o{f}|_{\o{X}-X}"]&j_Y
\end{tikzcd}\] 
where the bottom row consists of $j_X,j_Y$ viewed as morphisms of constant simplicial varieties (not necessarily smooth/projective).
\end{definition}
\begin{remark}\label{rem:covar}
For a proper morphism $f:X\to Y$, we can pick compactifications $\o{X}'$ and $\o{Y}$. Defining $\o{X}$ to be the closure of the graph of $f$ in $\o{X}'\times\o{Y}$, we observe that $\o{X}$ is a proper variety and a compactification of $X$. Let $\o{f}:\o{X}\to\o{Y}$ denote the restriction to $\o{X}$ of the second projection $\o{X}'\times\o{Y}\to\o{Y}$; by construction, $\o{f}|_X=f$. Label $j_X:\o{X}-X\to \o{X}$ and $j_Y:\o{Y}-Y\to \o{Y}$. Since the map $f:X\to Y$ is closed and $\o{f}|_X=f$, we have $\o{f}(\o{X}-X)\subset\o{Y}-Y$, so there is a restriction morphism of arrows $\o{f}:j_X\to j_Y$. Taking a smooth projective hyperenvelope $\tilde{f}:\tilde{j}_X\to\tilde{j}_Y$ of the map $\o{f}$ via Corollary \ref{cor:2.10}, we get a morphism $(\o{f},\tilde{f}):(X,\o{X},\tilde{j}_X)\to (Y,\o{Y},\tilde{j}_Y)$ in $\mathcal{V}^{hyp}$. In particular, every proper morphism $f:X\to Y$ of $\mathcal{V}$ lifts to a morphism of $\mathcal{V}^{hyp}$, i.e. the forgetful functor $\mathcal{V}^{hyp}\to\mathcal{V}^{prop}$ is surjective on objects and full. (We will not need the stronger statement that this functor is fibered in filtered categories, and we do not claim it.) 
\end{remark}
\begin{remark}\label{rem:contra}
Recall the way that Gillet and Soul\'{e} construct a map of Chow motives for an open inclusion. If $U\to X$ is an open immersion, and we choose a compactification $\o{X}$ of $X$, then since $\o{X}$ is a compactification of $U$ as well, we get a square:
\[\begin{tikzcd}
\o{X}-X\arrow[r,"j_X"]\arrow[d]&\o{X}\arrow[d,equal]\\
\o{X}-U\arrow[r,"j_U"]&\o{X}
\end{tikzcd}
\]
Calling this morphism $i_{f}:j_X\to j_U$, and taking a smooth projective hyperenvelope $\tilde{i}_f:\tilde{j}_X\to \tilde{j}_U$, we again see that $(i_f,\tilde{i}_f):(X,\o{X},\tilde{j}_X)\to (U,\o{X},\tilde{j}_U)$ gives us a morphism in $\mathcal{V}^{hyp}$.
\end{remark}
These two types of morphisms will be used heavily in the next section; first they will be used to construct morphisms of pro-objects of $\mathcal{V}^{hyp}$:
\begin{definition}\label{def:PX}
Given $X\in \mathcal{V}$, we define the \emph{pro-weight complex} of $X$ to be the subcategory of $\mathcal{V}^{hyp}$ consisting of objects:
\[P_X:=\{(X,\o{X},\tilde{j}_X)\in\mathcal{V}^{hyp}:\o{X}=\o{X}'\coprod *\text{ where }\o{X}'\text{ is a compactification of $X$}\}
\]
and morphisms in $P_X$ are of the form $(\o{f},\tilde{f}):(X,\o{X},\tilde{j}_X)\to (X,\o{X}',\tilde{j}_X')$ such that if $i:X\to\o{X}$, $i':X\to\o{X}'$ denote the open immersions, $\o{f}\circ i=i'$. In other words, $\o{f}$ leaves $X$ unchanged.
\end{definition}
\begin{remark}\label{rem:basepoint}
As in Definition \ref{def:Vhyp}, the second entry of an object of $P_X$ is a compactification of $X$ together with a disjoint copy of $\Spec(k)$, which serves as a base point; the disjoint base point is what allows us to regard the associated simplicial mapping cones as objects of $\s\V_*$, and it is used in an essential way in the proof of Theorem \ref{prop:contra}. For the sake of readability we will from now on suppress the base point from the notation, writing $\o{X}$ for a compactification of $X$ and leaving the disjoint copy of $\Spec(k)$ implicit. This is harmless: the base point is a disjoint summand, so all of the constructions below, namely closures, products, pullbacks, equalizers and hyperenvelopes, may be carried out on the compactifications with the base point removed and the base point adjoined afterwards, as we make explicit in part (2) of the proof of Proposition \ref{prop:proobj}.
\end{remark}
 We now show that $P_X$ is in fact a pro-object: 
\begin{proposition}\label{prop:proobj}
For all $X\in\mathcal{V}$, $P_X\in \Pro(\mathcal{V}^{hyp})$.
\end{proposition}
\begin{proof}
We need to show the following things:
\begin{enumerate}
    \item $P_X$ is non-empty, i.e. the weight complex of a variety can always be constructed, as reviewed in the last paragraph of Section 2.1
    \item For any two objects $(X,\o{X}_1,\tilde{j}_{X_1}); (X,\o{X}_2,\tilde{j}_{X_2})$, there exists $(X,\o{Y},\tilde{j}_Y)$ and morphisms:
    \[\begin{tikzcd}
    &(X,\o{X}_1,\tilde{j}_{X_1})\\
    (X,\o{Y},\tilde{j}_Y)\arrow[ur]\arrow[dr]&\\
    &(X,\o{X}_2,\tilde{j}_{X_2})
    \end{tikzcd}\]
    In order to show this, recall from Definition \ref{def:Vhyp} and Remark \ref{rem:basepoint} that each $\o{X}_i$ is of the form $\o{X}_i'\coprod *$. Throughout this part of the proof we make the base points explicit: we work with the compactifications $\o{X}_i'$, that is, with the disjoint base points removed, and adjoin a base point at the end. Let $\o{Y}'$ denote the closure of the image of $X$ under the composite
    \[X\xrightarrow{\Delta} X\times X\to \o{X}_1'\times\o{X}_2',
    \]
    where the first morphism is the diagonal and the second is the product of the open embeddings $X\to \o{X}_i'$. The diagonal is a locally closed embedding, and the product of the two open embeddings is an open embedding, so the composite is a locally closed embedding. Hence $X$ is open in its closure $\o{Y}'$, so that $\o{Y}'$ is a compactification of $X$. We then set $\o{Y}:=\o{Y}'\coprod *$; since the base point is a disjoint summand, the closures, pullbacks and hyperenvelopes taken below are unaffected by its presence. Let us denote:
    \[j_{X_1,X_2}:\o{Y}-X\to \o{Y}
    \]
    and label $\pi_i$, for $i=1,2$ to be the arrows of proper maps:
    \[\begin{tikzcd}
        \o{Y}-X\arrow[r,"j_{X_1,X_2}"]\arrow[d,"\pi_i"]&\o{Y}\arrow[d,"\pi_i"]\\
        \o{X}_i-X\arrow[r,"j_{X_i}"]&\o{X}_i
    \end{tikzcd}\]
    obtained by projection. We define the hyperenvelopes $j_i:=\tilde{j}_{X_i}\times_{j_{X_i}}j_{X_1,X_2}$, i.e. to be pullbacks of $j_{X_1,X_2}$ along $\pi_i$. Now, we can take a smooth projective hyperenvelope of the pullback of $j_1\times_{j_{X_1,X_2}}j_2$. Calling this $\tilde{j}_{X_1,X_2}$, we get the following diagram of hyperenvelopes
    \begin{equation}
    \begin{tikzcd}
    \tilde{j}_{X_1,X_2}\arrow[r]\arrow[d]&j_1\arrow[r]\arrow[d]&\tilde{j}_{X_1}\arrow[d]\\
    j_2\arrow[d]\arrow[r]&j_{X_1,X_2}\arrow[d,"\pi_2"]\arrow[r,"\pi_1"]&j_{X_1}\\
    \tilde{j}_{X_2}\arrow[r]&j_{X_2}
    \end{tikzcd}
    \end{equation}
    So, $\tilde{j}_{X_1,X_2}$ is a smooth projective hyperenvelope of $j_{X_1,X_2}$. If we now label $\tilde{\pi}_i:\tilde{j}_{X_1,X_2}\to \tilde{j}_{X_i}$, we can conclude that $(X,\o{Y},\tilde{j}_{X_1,X_2})$ is the desired object we need, as it is in $P_X$ and equipped with morphisms \[(\pi_i,\tilde{\pi}_i):(X,\o{Y},\tilde{j}_{X_1,X_2})\to (X,\o{X}_i,\tilde{j}_{X_i})\]
    \item We now need to show that given any two morphisms $(\o{f}_1,\tilde{f}_1),(\o{f}_2,\tilde{f}_2):(X,\o{X},\tilde{j}_X)\to (X,\o{Y},\tilde{j}_Y)$, there exists an object $(X,\o{Z},\tilde{j}_Z)$ and a morphism $(\o{g},\tilde{g}):(X,\o{Z},\tilde{j}_Z)\to (X,\o{X},\tilde{j}_X)$ such that $\o{f_1}\circ \o{g}=\o{f_2}\circ \o{g}$ and $\tilde{f_1}\circ \tilde{g}=\tilde{f_2}\circ \tilde{g}$ . \\\\
    We emphasize that, as noted in the footnote to Section 2.1, we do not require the embedding $X\to \o{X}$ to be dense. Consequently two morphisms $\o{X}\to \o{Y}$ agreeing on $X$ need not be equal, which is why an equalizer is needed here rather than an appeal to density. \\\\
    Now, given two morphisms of the form:
    \[\o{f}_1,\o{f_2}:\o{X}\to \o{Y}
    \]
    where $\o{X},\o{Y}$ are compactifications of $X$ and $Y$, then $eq(\o{f}_1,\o{f}_2)$ is a proper subvariety of $\o{X}$ by \cite[\href{https://stacks.math.columbia.edu/tag/01KM}{Tag 01KM}]{stacks-project} (which relies on $\o{Y}$ being separated), and if $i:X\to \o{X}$ denotes the embedding, then $\o{f}_1\circ i=\o{f}_2\circ i$ implies that $X$ is also an open subvariety of $eq(\o{f}_1,\o{f}_2)$. Set \[j_{eq}:=eq(\o{f}_1,\o{f}_2)-X\to eq(\o{f}_1,\o{f}_2)\]
    Our goal will be to show that the induced map on equalizers: $eq(\tilde{f}_1,\tilde{f}_2)\to eq(\o{f}_1,\o{f}_2)$ is a hyperenvelope. We stress that this equalizer is in general neither smooth nor projective; we claim only that it is a hyperenvelope of $j_{eq}$, and it will be replaced by a smooth projective hyperenvelope at the end of the proof. As discussed in 1.4.1 \cite{GS}, it suffices to check that for all fields $K$, the map on $K$ points $eq(\tilde{f}_1,\tilde{f}_2)(K)\to eq(\o{f}_1,\o{f}_2)(K)$ is a trivial Kan fibration, i.e. it is surjective on $0$ simplices and satisfies the boundary-filling condition. Recall that $\tilde{f}_i$ is a morphism of arrows, i.e. a commutative square. For each $i=1,2$, we write $\tilde{f}_i^1:\tilde{X}_\bullet\to \tilde{Y}_\bullet$ for the component of this square on the targets of the two arrows, and $\tilde{f}_i^2:\widetilde{\o{X}-X}_\bullet\to \widetilde{\o{Y}-Y}_\bullet$ for its component on the sources; we call these the first and second leg of $\tilde{f}_i$ respectively. We see that for the commutative cube:
    \[
    \begin{tikzcd}[row sep=1.5em, column sep = 1.5em]
    \tilde{eq}(\tilde{f}_1^1,\tilde{f}_2^1) \arrow[rr] \arrow[dr] \arrow[dd,swap] &&
    \tilde{Y}_\bullet \arrow[dd] \arrow[dr,"\Delta"] \\
    & \tilde{X}_\bullet\arrow[dd] \arrow[rr,"(\tilde{f}_1^1{,}\tilde{f}_2^1)"] &&
    \tilde{Y}\times\tilde{Y} \arrow[dd] \\
    eq(\o{f}_1,\o{f}_2) \arrow[rr] \arrow[dr, "f_2"] && \o{Y}\arrow[dr,"\Delta"] \\
    & \o{X} \arrow[rr,"(\o{f}_1{,}\o{f_2})"] &&\o{Y}\times\o{Y}
    \end{tikzcd}
\]
the upper and lower squares are pullbacks by definition of the equalizer and the vertical legs (except for $\tilde{eq}(\tilde{f}_1^1,\tilde{f}_2^1)\to eq(\o{f}_1,\o{f}_2)$ are hyperenvelopes. Passing to $K$ points, we see that given a diagram of simplicial sets of the form:
\[
\begin{tikzcd}
\partial\Delta[n]\arrow[r]\arrow[d,hook]& eq(\tilde{f}_1^1,\tilde{f}_2^1) (K)\arrow[d]\\
\Delta[n]\arrow[r]&eq(\o{f}_1,\o{f}_2)(K)
\end{tikzcd}
\]
We use that $\tilde{X}_\bullet(K)\to \o{X}(K)$ and $\tilde{Y}_\bullet(K)\to \o{Y}(K)$ are trivial Kan fibrations, giving us lifts:
\[
    \begin{tikzcd}[row sep=1.5em, column sep = 1.5em]
    \partial\Delta[n] \arrow[rr] \arrow[dr] \arrow[dd,hook,swap] &&
    \tilde{Y}_\bullet \arrow[dd] \arrow[dr,"\Delta"] \\
    & \tilde{X}_\bullet \arrow[rr,"(\tilde{f}_1^1{,}\tilde{f}_2^1)"]\arrow[dd] &&
    \tilde{Y}\times\tilde{Y} \arrow[dd] \\
     \Delta[n]\arrow[ur,dotted]\arrow[uurr,dotted] \arrow[rr] \arrow[dr, "f_2"] && \o{Y}\arrow[dr,"\Delta"] \\
    & \o{X} \arrow[rr,"(\o{f}_1{,}\o{f_2})"] &&\o{Y}\times\o{Y}
    \end{tikzcd}
\]
which then must factor through a lift of the form
\[
\begin{tikzcd}
\partial\Delta[n]\arrow[r]\arrow[d,hook]& eq(\tilde{f}_1^1,\tilde{f}_2^1) (K)\arrow[d]\\
\Delta[n]\arrow[r]\arrow[ur,dotted]&eq(\o{f}_1,\o{f}_2)(K)
\end{tikzcd}
\]
by universal property of the pullback, proving that $eq(\tilde{f}_1^1,\tilde{f}_2^1)(K)\to eq(\o{f}_1,\o{f}_2)(K)$ is a trivial Kan fibration for all $K$, and therefore a hyperenvelope. A similar argument shows that $eq(\tilde{f}_1^2,\tilde{f}_2^2)\to eq(\o{f}_1|_{\o{X}-X},\o{f}_2|_{\o{X}-X})$ is a hyperenvelope. So, \[\tilde{eq}:=eq(\tilde{f}_1^2,\tilde{f}_2^2)\to eq(\tilde{f}_1^1,\tilde{f}_2^1)\]
will be a hyperenvelope of $j_{eq}$. Taking a smooth projective hyperenvelope of $\tilde{eq}$, calling it $\tilde{j}_{eq}$, we see that
\[(X,eq(\o{f}_1,\o{f}_2),\tilde{j}_{eq})
\]
is an object of $P_X$ and the induced morphism to $(X,\o{X},\tilde{j}_X)$ equalizes $(\o{f}_1,\tilde{f}_1)$ and $(\o{f}_2,\tilde{f}_2)$ as desired. 
    \end{enumerate}
\end{proof}
\begin{definition}
We define \emph{the category of pro-weight complexes}, $PH(\mathcal{V})$, to be the full subcategory of $\Pro(\mathcal{V}^{hyp})$ consisting of the objects $P_X$, for $X\in\mathcal{V}$.
\end{definition}
Now, we further understand what morphisms between these pro-objects look like in $\Pro(\mathcal{V}^{hyp})$.
\begin{proposition}\label{prop:covar}
Given a proper morphism $f:X\to Y$, there exists a morphism of pro-objects $P_f:P_X\to P_Y$, represented by the pro-object of morphisms:
\[A_f:=\Bigg\{(\o{f},\tilde{f}):(X,\o{X},\tilde{j}_X)\to (Y,\o{Y},\tilde{j}_Y)\Bigg|\ (X,\o{X},\tilde{j}_X) \in P_X,\ (Y,\o{Y},\tilde{j}_Y)\in P_Y,\ \o{f}|_X=f \Bigg\}\]
Further, the assignment is functorial, i.e. we have a functor $P:\mathcal{V}^{prop}\to PH(\mathcal{V})$, where $\mathcal{V}^{prop}$ as defined in Definition \ref{def:var} is the category of varieties, but where the morphisms are restricted to only proper morphisms. 
\end{proposition}
\begin{proof}
We must show 3 things:
\begin{enumerate}
    \item $A_f$ is a pro-object in the category of morphisms of $\mathcal{V}^{hyp}$,
    \item The projection maps $A_f^{op}\to P_X^{op}$ and $A_f^{op}\to P_Y^{op}$ are cofinal,
    \item Let $g:X\to Y$, $f:Y\to Z$ be proper morphisms. Recall that
    \[S_{P_f\circ P_g}:=\Bigg\{(\o{f},\tilde{f})\circ (\o{g},\tilde{g}):(X,\o{X},\tilde{j}_X)\to (Z,\o{Z},\tilde{j}_Z)\ \Bigg|\ (\o{f},\tilde{f})\in A_f,\ (\o{g},\tilde{g})\in A_g\text{ composable}\Bigg\}
    \]
    where $(\o{f},\tilde{f})$ and $(\o{g},\tilde{g})$ are composable when the target of $(\o{g},\tilde{g})$ agrees with the source of $(\o{f},\tilde{f})$, both being objects of $P_Y$.
    which represents the morphism $P_f\circ P_g$ via the straightening lemma, Proposition \ref{prop:straightening}. We show that there is a cofinal functor $A_{f\circ g}\to S_{P_f\circ P_g}$ so that $A_{f\circ g}$ represents $P_f\circ P_g$ (not just $P_{f\circ g}$). This will prove functoriality.
\end{enumerate}
\begin{enumerate}
    \item[Proof of (1):] We follow the same logic as in the proof of Proposition \ref{prop:proobj}. Throughout this proof we write $\o{X}_1,\o{X}_2$ for two compactifications of $X$, and likewise $\o{Y}_1,\o{Y}_2$ for two compactifications of $Y$, following the notation of the proof of Proposition \ref{prop:proobj}; the constructions carried out there will be applied to them without further comment. First, observe that $A_f$ is non-empty by Remark \ref{rem:covar}. Given objects of $A_f$, i.e. two morphisms:
\[(\o{f},\tilde{f}):(X,\o{X}_1,\tilde{j}_{X_1})\to (Y,\o{Y}_1,\tilde{j}_{Y_1})\] and \[(\o{g},\tilde{g}):(X,\o{X}_2,\tilde{j}_{X_2})\to (Y,\o{Y}_2,\tilde{j}_{Y_2})\] such that $\o{g}|_{X}=\o{f}|_X=f$ we see that there is an induced map:
\[\o{f}\times\o{g}:j_{X_1,X_2}\to j_{Y_1,Y_2}
\]
where $j_{X_1,X_2}$ and $j_{Y_1,Y_2}$ are constructed as in part (2) of the proof of Proposition \ref{prop:proobj}, applied to the pairs of compactifications $\o{X}_1,\o{X}_2$ and $\o{Y}_1,\o{Y}_2$ respectively, giving us the diagram of arrows:
\[\begin{tikzcd}
j_{X_1}\arrow[r,"\o{f}"] &j_{Y_1}\\
j_{X_1,X_2}\arrow[u,"\pi_1"]\arrow[d,"\pi_2"]\arrow[r,"\o{f}\times\o{g}"] &j_{Y_1,Y_2}\arrow[u,"\pi_1"]\arrow[d,"\pi_2"]\\
j_{X_2}\arrow[r,"\o{g}"] &j_{Y_2}\\
\end{tikzcd}
\]
pulling back $\tilde{f}:\tilde{j}_{X_1}\to \tilde{j}_{Y_1}$ along the square $\pi_1$ to get $\tilde{f}\times_{\o{g}}(\o{f}\times\o{g})$ and similarly pulling back $\tilde{g}:\tilde{j}_{X_2}\to \tilde{j}_{Y_2}$ along the square $\pi_2$ to get $\tilde{g}\times_{\o{g}}(\o{f}\times\o{g})$ we see that both of these morphisms are hyperenvelopes of $j_{X_1,X_2}\to j_{Y_1,Y_2}$. Take a smooth projective hyperenvelope over $\tilde{f}\times_{\o{g}}(\o{f}\times\o{g})$ and $\tilde{g}\times_{\o{g}}(\o{f}\times\o{g})$ and call it $\widetilde{f\times g}$. Then, we have the following morphisms
\[\begin{tikzcd}
&(\o{f},\tilde{f})\\
(\o{f}\times\o{g},\widetilde{f\times g})\arrow[ur]\arrow[dr]\\
&(\o{g},\tilde{g})
\end{tikzcd}\]
in $A_f$. Next, observe that for a pair of arrows $(\o{f},\tilde{f}),(\o{g},\tilde{g}):(X,\o{X},\tilde{j}_X)\to (Y,\o{Y},\tilde{j}_Y)$ the same proof as in part (3) of the proof of Proposition \ref{prop:proobj} will show that there will be a map \[(X,eq(\o{f},\o{g}),\tilde{j}_{eq})\to (X,\o{X},\tilde{j}_X)\] equalizing both $(\o{f},\tilde{f})$ and $(\o{g},\tilde{g})$. We can thus conclude that $A_f$ is a pro-object of morphisms of $\mathcal{V}^{hyp}$. 
\item[Proof of (2):] We prove that each of the maps is cofinal by verifying the conditions of Definition \ref{def:cofinal} in order. Cofinality of $A_f^{op}\to P_Y^{op}$ follows from the fact that the functor $A_f\to P_Y$ is surjective on both objects (verifying the first condition) and morphisms (verifying the second condition). To see why this is true on objects, observe that for every $(Y,\o{Y}_1,\tilde{j}_{Y_1})$, we can construct $(X,\o{X}_1,\tilde{j}_{X_1})$ and a map $(X,\o{X}_1,\tilde{j}_{X_1})\to (Y,\o{Y}_1,\tilde{j}_{Y_1})$ in $A_f$ by Remark \ref{rem:covar}. Further, given any morphism $(Y,\o{Y}_1,\tilde{j}_{Y_1})\to (Y,\o{Y}_2,\tilde{j}_{Y_2})$, we see that if we compose it with a morphism $(X,\o{X}_1,\tilde{j}_{X_1})\to (Y,\o{Y}_1,\tilde{j}_{Y_1})$ (which we know already exists), the composition $(X,\o{X}_1,\tilde{j}_{X_1})\to (Y,\o{Y}_2,\tilde{j}_{Y_2})$ will also lie in $A_f$, meaning that $A_{f}\to P_Y$ is also surjective on morphisms.\\
In order to verify the conditions for cofinality of $A_f^{op}\to P_X^{op}$, more work is required. For the first condition, taking any object $(X,\o{X}_1,\tilde{j}_{X_1})$ in $P_X$, we see that for any morphism $(X,\o{X}_2,\tilde{j}_{X_2})\to (Y,\o{Y},\tilde{j}_Y)$ in $P_Y$, we have a morphism:\[(X,\o{X}_1\times\o{X}_2,\tilde{j}_{X_1,X_2})\to (X,\o{X}_2,\tilde{j}_{X_2})\to (Y,\o{Y},\tilde{j}_Y)\]
in $A_f$, and further, we have a morphism $(X,\o{X}_1\times\o{X}_2,\tilde{j}_{X_1,X_2})\to (X,\o{X}_1,\tilde{j}_{X_1})$ in $P_X$ proving the first condition for cofinality. As for the second, given a morphism $(X,\o{X}_1,\tilde{j}_{X_1})\to (Y,\o{Y},\tilde{j}_Y)$ in $A_f$ and two morphisms of $P_X$:
\[(\o{g},\tilde{g}),(\o{h},\tilde{h}):(X,\o{X}_1,\tilde{j}_{X_1})\to (X,\o{X}_2,\tilde{j}_{X_2})
\]
By part (3) of the proof of Proposition \ref{prop:proobj}, there is a map \[i:(X,eq(\o{g},\o{h}),\tilde{j}_{eq})\to (X,\o{X}_1,\tilde{j}_{X_1})\] equalizing $(\o{g},\tilde{g})$ and $(\o{h},\tilde{h})$, but also if we compose $i$ with $(X,\o{X}_1,\tilde{j}_{X_1})\to (Y,\o{Y},\tilde{j}_Y)$, we see that the induced morphism $(X,eq(\o{g},\o{h}),\tilde{j}_{eq})\to (Y,\o{Y},\tilde{j}_Y)$ will also lie in $A_f$. So, we have the commutative triangle:
\[\begin{tikzcd}
(X,eq(\o{g},\o{h}),\tilde{j}_{eq})\arrow[dr]\arrow[d,"i"]&\\
(X,\o{X}_1,\tilde{j}_{X_1})\arrow[r]&(Y,\o{Y},\tilde{j}_Y)
\end{tikzcd}
\]
which will get sent to the morphism $i$ via the projection $A_f\to P_X$, which shows that there is a morphism in $A_f$ whose image in $P_X$ equalizes $(\o{g},\tilde{g})$ and $(\o{h},\tilde{h})$, proving the last condition of cofinality of $A_f^{op}\to P_X^{op}$. So, we can conclude that if $f$ is proper, then $A_f$ defines a morphism $P_X\to P_Y$. 
\item[Proof of (3):] 
We first observe that objects of $A_{f\circ g}$ are morphisms of the form $(X,\o{X},\tilde{j}_X)\to (Z,\o{Z},\tilde{j}_Z)$ that \emph{do not} necessarily factor through any $(Y,\o{Y},\tilde{j}_Y)$. However, since the elements of $S_{P_f\circ P_g}$ are by construction composites of composable pairs, $S_{P_f\circ P_g}$ is contained in the subcategory of $A_{f\circ g}$ of morphisms that factor through an object of $P_Y$. We consider the map $A_{f\circ g}\hookrightarrow S_{P_f\circ P_g}$. In order to show cofinality, consider any morphism $(\o{h},\tilde{h}):(X,\o{X},\tilde{j}_X)\to (Z,\o{Z},\tilde{j}_Z)$ of $A_{f\circ g}$, i.e. such that $\o{h}|_{X}=f\circ g$, we pick a compactification $\o{Y}$ of $Y$ and denote $\o{Y}_1$ to be the closure of the image of $Y$ in $\o{Y}\times\o{Z}$, and similarly, let $\o{X}_1$ denote the closure of $X$ in $\o{X}\times\o{Y}$. Define $\o{g}_1:=\o{X}_1\to \o{Y}_1$ and $\o{f}:=\o{Y}_1\to \o{Z}$. If $j_Y:=\o{Y}_1-Y\to \o{Y}_1$, then $\o{f}$ also gives us an arrow $j_Y\to j_Z$. If we take a pullback of $\tilde{j}_Z$ along $\o{f}:j_Y\to j_Z$, we get a hyperenvelope $j_Y\times_{j_Z}\tilde{j}_Z$ of $j_Y$. Taking a smooth projective hyperenvelope $\tilde{j}_Y$ of $j_Y\times_{j_Z}\tilde{j}_Z$, we have a morphism:
\[(Y,\o{Y}_1,\tilde{j}_Y)\to (Z,\o{Z},\tilde{j}_Z)
\]
 If we define $j_{X,1}:=\o{X}_1-X\to \o{X}_1$, then pulling back $\tilde{j}_Y$ along the map $j_{X,1}\to j_Y$ and taking a smooth projective hyperenvelope $\tilde{j}_{X,1}$, we get a morphism $(X,\o{X}_1,\tilde{j}_{X,1})\to (Y,\o{Y}_1,\tilde{j}_Y)$. By cofinality of  $A_{f\circ g}$ in $P_X$ (proved in pt. 1 of this proof) we can find morphisms $(X,\o{X}_2,\tilde{j}_{X,2})\to (X,\o{X}_1,\tilde{j}_{X,1})$ and $(X,\o{X}_2,\tilde{j}_{X,2})\to (X,\o{X},\tilde{j}_{X})$ fitting into the diagram:
\[
\begin{tikzcd}
(X,\o{X}_2,\tilde{j}_{X,2})\arrow[r]\arrow[dr]&(X,\o{X}_1,\tilde{j}_{X,1})\arrow[r,"(\o{f}{,}\tilde{f})"]&(Y,\o{Y}_1,\tilde{j}_Y)\arrow[r,"(\o{g}{,}\tilde{g})"]&(Z,\o{Z},\tilde{j}_Z)\\
&(X,\o{X},\tilde{j}_X)\arrow[urr, "(\o{h}{,}\tilde{h})",swap]&&
\end{tikzcd}
\]
 and the morphism in the first row is in $S_{P_f\circ P_g}$. This proves the first condition for a cofiltering. For the second condition, we must consider a commutative diagram of the form:
\[
\begin{tikzcd}
(X,\o{X}_1,\tilde{j}_{X_1})\arrow[d]\arrow[d,shift right=2]\arrow[r,"(\o{g}{,}\tilde{g})"]&(Y,\o{Y}_1,\tilde{j}_{Y_1})\arrow[r,"(\o{f}{,}\tilde{f})"]&(Z,\o{Z}_1,\tilde{j}_{Z_1})\arrow[d]\arrow[d,shift right=2]\\
(X,\o{X}_2,\tilde{j}_{X_2})\arrow[rr,"(\o{g}{,}\tilde{g})"]&&(Z,\o{Z}_2,\tilde{j}_{Z_2})
\end{tikzcd}
\]
Since $P_Z$ and $P_X$ are cofiltered, we have morphisms $(Z,\o{Z}_{eq},\tilde{j}_{eq,Z})\to (Z,\o{Z}_1,\tilde{j}_{Z_1})$ and $(X,\o{X}_{eq},\tilde{j}_{eq,X})\to (X,\o{X}_1,\tilde{j}_{X_1})$ equalizing the right-hand side vertical morphisms. Next, since $A_{f}^{op}\to P_Z^{op}$ is cofinal, there is a map $(Y,\o{Y}_2,\tilde{j}_{Y_2})\to (Z,\o{Z}_{eq},\tilde{j}_{eq,Z})$ fitting into the diagram:
\[
\begin{tikzcd}
(Y,\o{Y}_2,\tilde{j}_{Y_2})\arrow[r]\arrow[d]&(Z,\o{Z}_{eq},\tilde{j}_{eq,Z})\arrow[d]\\
(Y,\o{Y}_1,\tilde{j}_{Y_1})\arrow[r,"(\o{f}{,}\tilde{f})"]&(Z,\o{Z}_1,\tilde{j}_{Z_1})
\end{tikzcd}
\]
Similarly, cofinality of $A_{g}^{op}\to P_Y^{op}$ gives us a commutative square
\[
\begin{tikzcd}
(X,\o{X}_2,\tilde{j}_{X_2})\arrow[r]\arrow[d]&(Y,\o{Y}_2,\tilde{j}_{Y_2})\arrow[d]\\
(X,\o{X}_1,\tilde{j}_{X_1})\arrow[r,"(\o{g}{,}\tilde{g})"]&(Y,\o{Y}_1,\tilde{j}_{Y_1})
\end{tikzcd}
\]
Lastly, using that $P_X$ is cofiltered (using both conditions of cofiltering) we can find $(X,\o{X}_3,\tilde{j}_{X_3})$ that fits into the diagram:
\[\begin{tikzcd}
&(X,\o{X}_3,\tilde{j}_{X_3})\arrow[dr]\arrow[dl]&\\
(X,\o{X}_{eq},\tilde{j}_{eq,X})\arrow[dr]&&(X,\o{X}_2,\tilde{j}_{X_2})\arrow[dl]\\
&(X,\o{X}_1,\tilde{j}_{X_1})&
\end{tikzcd}
\]
This gives us the commutative diagram:
\[
\begin{tikzcd}
(X,\o{X}_3,\tilde{j}_{X_3})\arrow[r]\arrow[d]&(Y,\o{Y}_2,\tilde{j}_{Y_2})\arrow[r]\arrow[d]&(Z,\o{Z}_{eq},\tilde{j}_{eq,Z})\arrow[d]\\
(X,\o{X}_1,\tilde{j}_{X_1})\arrow[d]\arrow[d,shift right=2]\arrow[r,"(\o{g}{,}\tilde{g})"]&(Y,\o{Y}_1,\tilde{j}_{Y_1})\arrow[r,"(\o{f}{,}\tilde{f})"]&(Z,\o{Z}_1,\tilde{j}_{Z_1})\arrow[d]\arrow[d,shift right=2]\\
(X,\o{X}_2,\tilde{j}_{X_2})\arrow[rr,"(\o{g}{,}\tilde{g})"]&&(Z,\o{Z}_2,\tilde{j}_{Z_2})
\end{tikzcd}
\]
So, $A_{f\circ g}^{op}\to S_{P_f\circ P_g}^{op}$ is cofinal, finishing the proof that $P:\mathcal{V}^{prop}\to PH(\mathcal{V})$ is a functor.
\end{enumerate}

\end{proof}
Now, we move towards defining a morphism $P_X\to P_U$ given an open inclusion $U\hookrightarrow X$:

\begin{definition}
Recall (see Remark \ref{rem:contra} for more details) that for an open immersion $f:U\to X$, we have a \emph{contravariant} morphism between weight complexes $M(X)\to M(U)$. In the notation of Remark \ref{rem:contra}, let
\[\Phi_{f}:=\{(id,\tilde{i}_{f}):(X,\o{X},\tilde{j}_X)\to (U,\o{X},\tilde{j}_U)\}_{(\o{X},\tilde{j}_X)\in (P_X)_U}
\]
be a subcategory of morphisms in $\mathcal{V}^{hyp}$. Here, $\tilde{i}_f$ is a hyperenvelope of the inclusion $j_X\to j_U$ induced by the identity map on $\o{X}$. The morphisms of $\Phi_f$ are of the form $(\o{g},\tilde{g}):(id,\tilde{i}_{f})\to (id,\tilde{i}_{f}')$ in the following sense: $\o{g}$ is the square
\[\begin{tikzcd}
j_X\arrow[r,"i_f"]\arrow[d]\arrow[d,"\o{g}_1"]&j_U\arrow[d,"\o{g}_2"]\\
j_X'\arrow[r,"i_f"]&j_U'
\end{tikzcd}
\]
and $\tilde{g}$ is the square:
\[\begin{tikzcd}
\tilde{j_X}\arrow[r,"\tilde{i}_{f}"]\arrow[d,"\tilde{g}_1"]&\tilde{j_U}\arrow[d,"\tilde{g}_2"]\\
\tilde{j_X'}\arrow[r,"\tilde{i}_{f}"]&\tilde{j_U'}
\end{tikzcd}
\]
where $(\o{g}_1,\tilde{g}_1)$ is a morphism of $P_X$ and $(\o{g}_2,\tilde{g}_2)$ is a morphism of $P_U$. 
\end{definition}
\begin{theorem}\label{prop:contra}
For any open immersion $f:U\to X$, $\Phi_f$ is cofiltered and represents a morphism $P_X\to P_U$. If $\textbf{co}(\mathcal{V})$ denotes the category of open immersions of $\mathcal{V}$, then $\Phi$ defines a contravariant functor:
\[\textbf{co}(\mathcal{V})^{op}\to PH(\mathcal{V})
\]
\end{theorem}
\begin{proof}
Just like in the proof of Proposition \ref{prop:covar}, we must show the following 3 things:
\begin{enumerate}
    \item $\Phi_f$ is a pro-object of morphisms of $\mathcal{V}^{hyp}$
    \item The projection maps $\Phi_f^{op}\to P_X^{op}$ and $\Phi_f^{op}\to P_U^{op}$ are cofinal
    \item For a composition of embeddings $U\xrightarrow{f}V\xrightarrow{g} X$, we label the category of morphisms representing $\Phi_f\circ\Phi_g$ as:
    \[S_{\Phi_f\circ\Phi_g}
    \]
    We must show that the inclusion $\Phi_{g\circ f}\to S_{\Phi_f\circ\Phi_g}$ is cofinal.
\end{enumerate}
For (1)-(3), we follow the same strategy as the proof of Proposition \ref{prop:covar}. In (1) we verify each condition of a pro-object as in the proof of Proposition \ref{prop:covar}, and for (2) and (3) we prove both cofinality conditions for each functor in the order of Definition \ref{def:cofinal}.\\
\begin{enumerate}
    \item[Proof of (1):] We see that by Remark \ref{rem:contra}, $\Phi_f$ is non-empty. Next, given any two morphisms $i_k:=(i_{f,k},\tilde{i}_{f,k})\in \Phi_f$ for $k=1,2$, we need to find another morphism $i':=(i_{f},\tilde{i}_{f}')$ and arrows $i'\to i_k$ in $\Phi_f$ for $k=1,2$. So, if $\tilde{i}_{f,k}$ are the diagrams:
\[\begin{tikzcd}
\widetilde{\o{X}_k-X}_\bullet\arrow[r]\arrow[d]&\tilde{X}_{k,\bullet}\arrow[d]\\
(\widetilde{\o{X}_k-U})_{k,\bullet}\arrow[r]&\tilde{X}_{k,\bullet}
\end{tikzcd}
\]
which are smooth projective hyperenvelopes of $i_{f,k}$:
\[\begin{tikzcd}
\o{X}_k-X\arrow[r]\arrow[d]&\o{X}_k\arrow[d,equal]\\
\o{X}_k-U\arrow[r]&\o{X}_k
\end{tikzcd}
\]
Recalling that $\o{X}_1\times\o{X}_2$ is a compactification of $X$, writing $\o{X}':=\o{X}_1\times\o{X}_2$, $j_X':=\o{X}'-X\to \o{X}'$, $j_U':=\o{X}'-U\to \o{X}'$, and $j_U:=\o{X}'-U\to \o{X}'$ we have the following diagrams for each $k=1,2$:
\[
\begin{tikzcd}
j_X'\arrow[r]\arrow[d]&j_U'\arrow[d]\\
j_{X,k}\arrow[r]&j_{U,k}
\end{tikzcd}
\]
induced by the projections $\o{X}'\to \o{X}_k$. For each $k$, we can find smooth projective hyperenvelopes of the above square, i.e. the diagrams
\[
\begin{tikzcd}
\tilde{j_{X',k}}\arrow[r]\arrow[d]&\tilde{j_{U',k}}\arrow[d]\\
\tilde{j_{X,k}}\arrow[r]&\tilde{j_{U,k}}
\end{tikzcd}
\]
Take a smooth projective hyperenvelope $\tilde{j}_U'$ over both $\tilde{j_{U',k}}$. Next, take a smooth projective hyperenvelope $\tilde{j_X}'$ over the hyperenvelopes (for $k=1,2$) $\tilde{j}_U'\times_{j_U'}\tilde{j_{X',k}}$. Labelling $\tilde{i}_f'=\tilde{j_X}'\to \tilde{j_U}'$, we see that $(X,\o{X},\tilde{j}_X')\in (P_X)_U$ and \[i_f':=(id,\tilde{i}_f'):(X,\o{X},\tilde{j}_X')\to (U,\o{X},\tilde{j}_U')\] is an element of $\Phi_f$. Further, we have morphisms $i_f'\to i_{k}$ in $\Phi_f$. \\\\
We also need to show that for any parallel arrows $(\o{\phi},\tilde{\phi}),(\o{\phi}',\tilde{\phi'}):i_1=(i_{f,1},\tilde{i}_{f_1})\to i_2=(i_{f,2},\tilde{i}_{f_2})$ in $\Phi_f$, there exists a morphism $i'\in \Phi_f$ and an arrow $(\o{\psi},\tilde{\psi}):i'\to i_1$ such that $(id,\tilde{\phi'})\circ (id,\tilde{\psi})=(id,\tilde{\phi})\circ (id,\tilde{\psi})$. Firstly, $\o{\phi},\o{\phi'}$ give us the diagram:
\[
\begin{tikzcd}
j_X\arrow[r,"i_{f,1}"]\arrow[d,"\o{\phi}_1"]\arrow[d,shift right=2,"\o{\phi}_1'",swap]&j_U\arrow[d,"\o{\phi}_2",swap]\arrow[d,shift left=2,"\o{\phi}_2'"]\\
j_X\arrow[r,"i_{f,2}"]&j_U'
\end{tikzcd}
\]
and a similar diagram exists for $\tilde{\phi},\tilde{\phi'}$. Recalling the constructions of $(X,eq(\o{\phi_1},\o{\phi_1}'),\tilde{j}_{eq,X})\in P_X$ and $(U,eq(\o{\phi_2},\o{\phi_2}'),\tilde{j}_{eq,U})\in P_U$ in part 3 of Proposition 5.5, we first observe that on $\o{X}, eq(\o{\phi_1},\o{\phi_1}')=eq(\o{\phi_2},\o{\phi_2}')$, so we get an arrow of morphisms
\[\begin{tikzcd}
eq(\o{\phi_1},\o{\phi_1}')-X\arrow[d]\arrow[r]&eq(\o{\phi_1},\o{\phi_1}')\arrow[d]\\
eq(\o{\phi_2},\o{\phi_2}')-U\arrow[r]&eq(\o{\phi_2},\o{\phi_2}')
\end{tikzcd}
\]
Write the top morphism as $j_{eq,X}$ and the bottom morphism as $j_{eq,U}$. Taking $\tilde{j}_{eq,X}'$ to be a smooth projective hyperenvelope of $j_{eq,X}$ over both $j_{eq,X}\times_{j_{eq,U}}\tilde{j}_{eq,U}$ and $\tilde{j}_{eq,X}$, we not only see that $(X,eq(\o{\phi_1},\o{\phi_1}'),\tilde{j}_{eq}')\in P_X$, but we also have the commutative diagrams
\[
\begin{tikzcd}
j_{eq,X}\arrow[d]\arrow[r]&j_{eq,U}\arrow[d]\\
j_X\arrow[r,"i_{f,1}"]\arrow[d,"\o{\phi}_1"]\arrow[d,shift right=2,"\o{\phi}_1'",swap]&j_U\arrow[d,"\o{\phi}_2",swap]\arrow[d,shift left=2,"\o{\phi}_2'"]\\
j_X\arrow[r,"i_{f,2}"]&j_U'
\end{tikzcd}
\]
and
\[
\begin{tikzcd}
\tilde{j}_{eq,X}'\arrow[d]\arrow[r,"\tilde{i}_{f,eq}"]&\tilde{j}_{eq,U}\arrow[d]\\
\tilde{j}_X\arrow[r,"\tilde{i}_{f,1}"]\arrow[d,"\tilde{\phi}_1"]\arrow[d,shift right=2,"\o{\phi}_1'",swap]&\tilde{j}_U\arrow[d,"\tilde{\phi}_2",swap]\arrow[d,shift left=2,"\tilde{\phi}_2'"]\\
\tilde{j}_X\arrow[r,"\tilde{i}_{f,2}"]&\tilde{j}_U'
\end{tikzcd}
\]
So the morphism \[(id,\tilde{i}_{f,eq}):(X,eq(\o{\phi_1},\o{\phi_1}'),\tilde{j}_{eq,X}')\to (U,eq(\o{\phi_1},\o{\phi_1}'),\tilde{j}_{eq,U}')\] in $\Phi_f$ satisfies the desired properties.
    \item[Proof of (2):] In order to show that the projection $\Phi_f^{op}\to P_X^{op}$ is cofinal, the first condition is to show that for any object $(X,\o{X},\tilde{j}_X)$ of $P_X$, there is a morphism $(X',\o{X}',\tilde{j}_X')\to (X,\o{X},\tilde{j}_X)$ such that $(X',\o{X}',\tilde{j}_X')$ is the target of a morphism in $\Phi_f^{op}$. We see that by Remark \ref{rem:contra} the functor is surjective on objects, which satisfies the first condition of cofinality. For the second condition, suppose we are given the diagram
\begin{equation}\label{eq:parallel}
\begin{tikzcd}
(X,\o{X},\tilde{j}_X)\arrow[d]\arrow[d,shift right=2]\arrow[r] &(U,\o{X},\tilde{j}_U)\\
(X,\o{X}',\tilde{j}_X')
\end{tikzcd}
\end{equation}
where the morphism $(X,\o{X},\tilde{j}_X)\to (U,\o{X},\tilde{j}_U)$ is in $\Phi_f$. Then, since $P_X$ is cofiltered, there exists a morphism $(X,\o{X}_{eq},\tilde{j}_{eq})\to (X,\o{X},\tilde{j}_X)$ equalizing the vertical maps. Further, $\o{X}_{eq}$ is also a compactification of $U$, and if we write $j_{eq,U}:\o{X}_{eq}-U\to \o{X}_{eq}$ and $j_U:\o{X}-U\to \o{X}$, we get a morphism $j_{eq}\to j_U$. Pulling back $\tilde{j}_U$ along this map and taking a smooth projective hyperenvelope $\tilde{j}_{eq,U}$, we have the diagram
\[
\begin{tikzcd}
(X,\o{X}_{eq},\tilde{j}_{eq})\arrow[r]\arrow[d]&(U,\o{X}_{eq},\tilde{j}_{eq,U})\arrow[d]\\
(X,\o{X},\tilde{j}_X)\arrow[d]\arrow[d,shift right=2]\arrow[r] &(U,\o{X},\tilde{j}_U)\\
(X,\o{X}',\tilde{j}_X')
\end{tikzcd}
\]
In other words, there is a morphism in $\Phi_f$ whose image in $P_X$ equalizes the vertical morphisms in \eqref{eq:parallel}, which is the second condition of cofinality. 

Next we must show that $\Phi_f^{op}\to P_U^{op}$ is cofinal, which is harder. Given an object $(U,\o{U}\coprod *,\tilde{j}_U)$ of $P_U$, if $i:U\to\o{U}$ denotes the open immersion, write $\o{im(i)}$ as the closure of $U$ in $\o{U}$. $U$ is Zariski dense in $\o{im(i)}$. We also have that $U$ is Zariski dense in the closure of its image in $X$, so taking $V$ to be the complement of its Zariski closure, we have that $U\coprod V$ is Zariski dense in $X$. Given compactifications $V\to \o{V}$ and $X\to\o{X}$ such that the images are dense, we consider the closure of the image of the graph morphism:
\[\Gamma_f: U\coprod V\to (\o{im(i)}\coprod \o{V})\times \o{X}
\]
denoted $\o{\Gamma_f}\subset (\o{im(i)}\coprod \o{V})\times \o{X}$. Observe that $\o{\Gamma_f}$ is a compactification of $U$, and the morphism $\o{\Gamma_f}\to \o{X}$ is proper and birational. By 2.11 of \cite{Conrad} there exists $U\coprod V$-admissible blowups of $\o{\Gamma_f}$ and $\o{X}$ that are isomorphic, i.e. two blowup maps $\tilde{X}\to \o{\Gamma_f}$ and $\tilde{X}\to \o{X}$ with open immersion $U\coprod V\to \tilde{X}$ compatible with the blow-up map. In the proof 2.11 and as explained in Lemma 2.7 of loc. cit., we can actually make the blowup $\tilde{X}\to \o{X}$ to be $X$-admissible as $X$ both contains $U\coprod V$ and is also a Zariski dense open subset of $\o{X}$. Further, we have a map $\o{im(i)}\coprod \o{V}\to\o{im(i)}\coprod *$ sending $\o{V}$ to $*$. So, we have a composition of proper maps:
\[\tilde{X}\to \o{\Gamma_f}\to (\o{im(i)}\coprod \o{V})\to (\o{im(i)}\coprod *)\to (\o{U}\coprod *)\]
that will send $\tilde{X}-U$ to $\o{U}-U$
and when restricted to $U$ is the identity. We know we can construct a smooth projective hyperenvelope $\tilde{j}_U'$ that gives us a morphism $(U,\tilde{X},\tilde{j}_U')\to (U,\o{U},\tilde{j}_U)$. Next, we can find a hyperenvelope $\tilde{j}_X$ of $\tilde{X}-X\to \tilde{X}$ and a morphism $\tilde{j}_X\to \tilde{j}_U'$ that allows us to construct the diagram of morphisms:
\begin{equation}\label{eq:domcompact}
\begin{tikzcd}
    (U,\tilde{X},\tilde{j}_U')\arrow[d]&(X,\tilde{X},\tilde{j}_X)\arrow[l]\\
    (U,\o{U}\coprod *,\tilde{j}_U)
\end{tikzcd}
\end{equation}
which proves the first cofinality condition. For the second condition, we are given a pair of morphisms $(\o{f},\tilde{f}),(\o{g},\tilde{g}):(U,\o{X},\tilde{j}_U)\to (U,\o{U},\tilde{j}_U')$ such that there exists a morphism $(X,\o{X},\tilde{j}_X)\to (U,\o{X},\tilde{j}_U)$. Let $\o{im(X)}\subset \o{X}$ denote the closure of $X$ inside $\o{X}$ and $\o{im(U)}\subset eq(\o{f},\o{g})$ denote the closure of $U$ inside the equalizer (which is a compactification of $U$). Just as in \eqref{eq:domcompact}, we can find a compactification of $X$ and $U$ dominating $\o{im(U)}$ and $\o{im(X)}$ respecting the open embedding $U\hookrightarrow X$, call it $\tilde{X}$. Since both $\o{im(U)}$ and $\o{im(X)}$ are closed subsets of $\o{X}$, we have the commutative diagram:
\begin{equation}
\begin{tikzcd}
\tilde{X}\arrow[r,equal]\arrow[d]&\tilde{X}\arrow[d]\\
    eq(\o{f},\o{g})\arrow[d]&\o{im(X)}\arrow[d]\\
    \o{X}\arrow[r,equal]& \o{X}
\end{tikzcd}
\end{equation}
which respect the open immersion $U\hookrightarrow X$. If $(U,eq(\o{f},\o{g}),\tilde{j}_{eq})$ denotes the object that equalizes $(\o{f},\tilde{f})$ and $(\o{g},\tilde{g})$, then we can find hyperenvelopes $\tilde{j}_X'$ of $\tilde{X}-X\to \o{X}$ and $\tilde{j}_U'$ of $\tilde{X}-U\to \tilde{X}$ that give us the diagram:
\[\begin{tikzcd}
\tilde{j}_U'\arrow[d]&\tilde{j}_X'\arrow[dd]\arrow[l]\\
    \tilde{j}_{eq}\arrow[d]&\\
    \tilde{j}_U&\tilde{j}_X\arrow[l]
\end{tikzcd}\]
i.e. image of the morphism $(X,\tilde{X},\tilde{j}_X')\to (X,\o{X},\tilde{j}_X)$ in $P_U$ will equalize $(\o{f},\tilde{f})$ and $(\o{g},\tilde{g})$, proving cofinality.
\item[Proof of (3):] The proof that $\Phi_{g\circ f}\hookrightarrow S_{\Phi_f\circ \Phi_g}$ is the same as (3) in the proof of Proposition \ref{prop:covar}.
\end{enumerate} 

\end{proof}
Lastly, we would like to relate $PH(\mathcal{V})$ to $\Pro(s\mathcal{V}_*^{\perf})$. This is done via the simplicial mapping cone:
\begin{definition}\label{prefunctor}

The simplicial cone defines a functor, which we call the \emph{pre-motive functor} 
\[M:PH(\mathcal{V})\to \Pro^{wc}(\s\V_*^{\perf})
\]
where on objects $M_X:=M(P_X)=\{sCone(\tilde{j}_{X}): (X,\o{X},\tilde{j}_X)\in P_X\}$. For morphism $\Phi:P_X\to P_Y$, we define $M(\Phi):=\{sCone(\tilde{f}): (\o{f},\tilde{f})\text{ represents $\Phi$}\}$.
\end{definition}

\section{Lifting the Measure to $\s\V_*^{\perf}$}
We are now in a position to prove the existence of a well-defined weight filtration in $L^H(\s\V_*^\perf)$ equipped with covariant morphisms for closed embeddings and contravariant morphisms for open immersions, and to build the map on K-theory spectra via the construction of a weak W-exact functor in the sense of \cite{Cam}. Recall that in \cite{Cam}, $\mathcal{V}$ is shown to be an \emph{SW}-category, or a category with three distinguished classes of morphisms:
\begin{enumerate}
    \item \emph{Cofibrations}, consisting of closed embeddings. Let $\textbf{co}(\mathcal{V})$ denote the corresponding category.
    \item \emph{Complements}, consisting of open immersions. Let $\textbf{comp}(\mathcal{V})$ denote the corresponding category.
    \item \emph{Weak equivalences}, consisting of isomorphisms. Let $\textbf{w}(\mathcal{V})$ denote the corresponding category.
\end{enumerate}

We use these subcategories, along with work from Section 4 and 5 to state the following theorem:
\begin{theorem}\label{RefinedGS}
    The pre-motive functor $X\mapsto M_X$ induces the following infinity functors:
    \[\textbf{\emph{w}}(\mathcal{V})\to L^H(\s\V_*^\perf)\]
    \[
    \textbf{\emph{co}}(\mathcal{V})\to L^H(\s\V_*^\perf)
    \]
    \[
    \textbf{\emph{comp}}(\mathcal{V})^{op}\to L^H(\s\V_*^\perf)
    \]
    that on $\pi_0$ sends $X\mapsto M(X)\in \Ho(\s\V_*^\perf)$. In other words, the Gillet--Soul\'{e} premotive construction $M(X)$ is well-defined in $L^H(\s\V_*^\perf)$, and comes equipped with covariant morphisms associated to closed embeddings along with contravariant morphisms associated to open immersions.
\end{theorem}
\begin{proof}
    Composing Proposition \ref{prop:covar} with Definition \ref{prefunctor}, we have a functor of 1-categories:
    \[
    \textbf{co}(\mathcal{V})\to \Pro^{wc}(\s\V_*^{\perf})
    \]
    The fact that $M_X$ in fact lands in weakly constant pro-objects of $\s\V_*^{perf}$ follows from Theorem \ref{prop:boundedweight}. Dwyer-Kan hammock localization gives a map of $\infty$-categories \[\Pro^{wc}(\s\V_*^{\perf})\to L^H(\Pro^{wc}(\s\V_*^{\perf}))\] By Theorem \ref{thm:DKequiv}, we have a map $\Pro^{wc}(\s\V_*^{\perf})\to L^H(\s\V_*^\perf)$, and so a map $\textbf{co}(\mathcal{V})\to L^H(\s\V_*^\perf)$. Restricting to $\textbf{w}(\mathcal{V})\subset \textbf{co}(\mathcal{V})$, we acquire the first $\infty$-functor in the theorem. The last $\infty$-functor in the theorem is constructed in the same way as $\textbf{co}(\mathcal{V})\to L^H(\s\V_*^\perf)$ but using the functor in Theorem \ref{prop:contra}. The fact that they all agree on 0-simplices follows from the fact that they all agree on objects in $\Pro^{wc}(\s\V_*^\perf)$. 
\end{proof}
Instead of distinguished exact sequences in an exact category or quotients in a Waldhausen category, $\mathcal{V}$ has said to have \emph{subtraction sequences} of the form $Z\to X\leftarrow U$, where the first morphism is a cofibration, the second is a complement, and $X-Z=U$. Associated to an SW-category is a K-theory spectrum, first constructed in [Z] and then in \cite{Cam}. In order to construct maps from the K-theory spectrum of an SW category to the K-theory spectrum of a Waldhausen category, the notion of a W-exact functor is introduced in \cite{Cam}. We recall the definition here. 

\begin{definition}
Given an SW-category $\mathcal{C}$ and a Waldhausen category $\mathcal{W}$ with FFWC, a weak \emph{$W$-exact functor} is a triple of functors $(F_!,F^!,F_w)$ such that:
\begin{enumerate}
    \item $F_!:\textbf{co}(\mathcal{C})\to \mathcal{W}$
    \item $F^!:\textbf{comp}(\mathcal{C})^{op}\to \mathcal{W}$
    \item $F_w:\textbf{w}(\mathcal{C})\to \textbf{w}(\mathcal{W})$
    \item For an object $X\in ob(\mathcal{C})$, $F_!(X)=F^!(X)=F_w(X)$
    \item For any pullback diagram in $\mathcal{C}$ (on the left) where $i,i'$ are complements and $j,j'$ are cofibrations, the corresponding diagram commutes in $\mathcal{W}$ (on the right):
    \[\begin{tikzcd}
    X\arrow[r,"j"]\arrow[d,"i"]&Y\arrow[d,"i'"]\\
    Z\arrow[r,"j'"]&W
    \end{tikzcd}\qquad \begin{tikzcd}
    F(X)\arrow[r,"F_!(j)"]&F(Y)\\
    F(Z)\arrow[r,"F_!(j')"]\arrow[u,"F^!(i)"]&F(W)\arrow[u,"F^!(i')",swap]
    \end{tikzcd}\]
   
    \item For a subtraction sequence $Z\xrightarrow[]{f} X\xleftarrow{i} U$, the corresponding diagram in $\mathcal{W}$:
    \[\begin{tikzcd}
    F(Z)\arrow[r,"F_!(f)"]\arrow[d]&F(X)\arrow[d,"F^!(i)"]\\
    0\arrow[r]&F(U)
    \end{tikzcd}
    \]
    is weakly equivalent to a pushout diagram, i.e. it is connected to a pushout square by a zig-zag of weak equivalences of diagrams (cf. Definition \ref{def:exact}).
    \item For any commutative diagram in $\mathcal{C}$ (on the left) where $i,i'$ are weak equivalences and $j,j'$ are cofibrations, the corresponding diagram commutes in $\mathcal{W}$ (on the right):
    \[\begin{tikzcd}
    X\arrow[r,"j"]\arrow[d,"i",swap]&Y\arrow[d,"i'"]\\
    Z\arrow[r,"j'"]&W
    \end{tikzcd}\qquad \begin{tikzcd}
    F(X)\arrow[r,"F_!(j)"]\arrow[d,"F_w(i)",swap]&F(Y)\arrow[d,"F_w(i')"]\\
    F(Z)\arrow[r,"F_!(j')"]&F(W)
    \end{tikzcd}\]
    and a similar statement holds for complement maps.
\end{enumerate}
\end{definition}
The following proposition is directly from \cite{Cam}:
\begin{theorem}\label{thm:Cam}
Let $\mathcal{C}$ be an SW category, $\mathcal{W}$ be a Waldhausen category with FFWC, and $(F_!,F^!,F_w):\mathcal{C}\to\mathcal{W}$ be a weak $W$-exact functor. Then, $(F_!,F^!,F_w)$ determines a map of K-theory spectra:
\[K(\mathcal{C})\to K(\mathcal{W})
\]
which on $\pi_0$ sends the class $[X]$ of an object $X\in \mathcal{C}$ to the class $[F_!(X)]\in K_0(\mathcal{W})$.
\end{theorem}
We now use this definition and theorem to get the desired lift.
\begin{theorem}\label{thm:liftmeasure}
There exists a weak $W$-exact functor $(F_!,F^!,F_w):\mathcal{V}\to \Pro^{wc}(\s\V_*^{\perf})$ sending on objects:
\[X\mapsto M_X:=\{sCone(\tilde{j}_{X}): (X,\o{X},\tilde{j}_X)\in P_X\}
\]
Consequently, this assignment gives us a map of K-theory spectra, which we will label
\[\Pro(V):K(\mathcal{V})\to K(\Pro^{wc}(\s\V_*^{\perf}))
\]
\end{theorem}
\begin{definition}
Recalling the axioms in 5.2 of \cite{Cam}, we define for variety $X\in ob(\mathcal{V})$, \[F^!(X)=F_w(X)=F_!(X)=M_X\] Now, we define them on morphisms.
\begin{enumerate}
    \item $F_!:\textbf{co}(\mathcal{V})\to \Pro^{wc}(\s\V_*^{\perf})$ is defined by sending a cofibration $f:Z\to X$ to the following morphism of pro-objects:
    \[F_!(f):= M(P_f): M_Z\to M_X
    \]
    where $P_f$ is the morphism constructed in Proposition \ref{prop:covar}.
    \item $F^!:\textbf{comp}(\mathcal{V})^{op}\to \Big(\Pro^{wc}(\s\V_*^{\perf})\Big)$ is defined by sending an open immersion $f:U\to X$ to 
    \begin{align*}
        F^!(f):=M(\Phi_f):M_X\to M_U
    \end{align*}
    where $\Phi_f$ is the morphism constructed in Theorem \ref{prop:contra}.
    \item $F_w:\textbf{w}(\mathcal{V})\to \Pro^{wc}(\s\V_*^{\perf})$ is just defined as the restriction of $F_!$ to $\textbf{w}(\mathcal{V})$.
\end{enumerate}
\end{definition}
\begin{proof}[Proof of Theorem \ref{thm:liftmeasure}]
     The first 4 conditions of $W$-exactness can be immediately verified. Showing property 5 requires the most work. Given a diagram:
    \[
    \begin{tikzcd}
Z\arrow[r,hook,"f"]\arrow[d,"g"]&X\arrow[d,"g'"]\\
        W\arrow[r,hook,"f'"]&Y
    \end{tikzcd}
    \]
    where $f,f'$ are closed embeddings and $g,g'$ are open immersions, we show that we have a commutative diagram 
    \begin{equation}\label{eq:PHsquare}
    \begin{tikzcd}
P_Z \arrow[r,"P_f"]&P_X\\
P_{W}\arrow[u,"\Phi_g"]\arrow[r,"P_{f'}"]&P_{Y}\arrow[u,"\Phi_{g'}",swap]
    \end{tikzcd}
    \end{equation}
    in $PH(\mathcal{V})$. This can be done by observing that the subcategory of arrows of $PH(\mathcal{V})$ consists of diagrams (not necessarily commutative) in $PH(\mathcal{V})$ of the form:
    \begin{equation}\label{eq:PHrep}
    \begin{tikzcd}
        (Z,\o{Z},\tilde{j}_Z)\arrow[r,"(\o{f}{,}\tilde{f})"]&(X,\o{X},\tilde{j}_X)\\
        (W,\o{Z},\tilde{j}_{W})\arrow[r,"(\o{f}'{,}\tilde{f}')"]\arrow[u,"(id{,}\tilde{g})"]&(Y,\o{X},\tilde{j}_{Y})\arrow[u,"(id{,}\tilde{g}')",swap]
    \end{tikzcd}
    \end{equation}
Where
\begin{equation}\label{eq:comp1}(W,\o{Z},\tilde{j}_{W})\to (Z,\o{Z},\tilde{j}_Z)\to (X,\o{X},\tilde{j}_X)\end{equation}
represents $P_f\circ \Phi_g$ and 
\begin{equation}\label{eq:comp2}(W,\o{Z},\tilde{j}_{W})\to (Y,\o{X},\tilde{j}_{Y})\to (X,\o{X},\tilde{j}_X)\end{equation}
represents $\Phi_{g'}\circ P_{f'}$. This is a pro-object of morphisms representing \eqref{eq:PHsquare} (see \cite{AM} A.3.4 for the proof). Now, A.3.5 of \cite{AM} tells us that we just need to find for each diagram of the form \eqref{eq:PHrep} a morphism $(W,\o{W}',\tilde{j}_W')\to (W,\o{Z},\tilde{j}_W)$ equalizing \eqref{eq:comp1} and \eqref{eq:comp2}. Observe that the morphisms $j_W\to j_X$ in \eqref{eq:comp1} and \eqref{eq:comp2} are already equal, giving us the commutative diagram:
\[
\begin{tikzcd}
    \tilde{j}_W\arrow[d]\arrow[r,shift left=1]\arrow[r]&\tilde{j}_X\arrow[d]\\
    j_W\arrow[r]&j_X
\end{tikzcd}
\]
taking the equalizer of the two morphisms from $\tilde{j}_W\to\tilde{j}_X$, we have that by the same proof as in part (3) of the proof of Proposition \ref{prop:proobj}, the equalizer of these two morphisms which we label $j_{eq}$, is a hyperenvelope of $j_W$. Taking a smooth projective hyperenvelope $\tilde{j}_{eq}$ of $j_{eq}$, we get a morphism of the form $(W,\o{W},\tilde{j}_{eq})\to (W,\o{W},\tilde{j}_W)$ equalizing \eqref{eq:comp1} and \eqref{eq:comp2} as desired, proving that \eqref{eq:PHsquare} is a commutative diagram.\\\\
For condition 6, given a subtraction sequence $Z\to X\leftarrow U$, the proof of Theorem 2.4 of \cite{GS} constructs a morphism $(Z,\o{Z},\tilde{j}_Z)\to (X,\o{X},\tilde{j}_X)\to (U,\o{X},\tilde{j}_U)$ such that the induced morphism on cones:
\[sCone(\tilde{j}_Z)\to sCone (\tilde{j}_X)\to sCone(\tilde{j}_U)\] is a pushout square in $\s\V_*$. Since we have a weak equivalence of pro-objects:
\[\begin{tikzcd}
    M_Z\arrow[r]\arrow[d,"\simeq"]&M_X\arrow[d,"
    \simeq"]\arrow[r]&M_U\arrow[d,"\simeq"]\\
    sCone(\tilde{j}_Z)\arrow[r]& sCone (\tilde{j}_X)\arrow[r]& sCone(\tilde{j}_U)
\end{tikzcd}\]
we can conclude that the upper row is weakly equivalent to a pushout square. \\\\
For condition 7, $F_w=F_!$ so the proof of condition 7 holds by functoriality. 
\end{proof}
\begin{corollary}\label{cor:sVmotivicmeasure}
    There is a map of spectra (the derived motivic measure associated to $\s\V_*^{\perf}$)
     \[K(\mathcal{V})\to K(\s\V_*^{\perf})\]
     which at the level of $K_0$ sends $[X]\mapsto [sCone(\tilde{j}_X)]$, where $\tilde{j}_X$ is any smooth projective hyperenvelope of $j_X:\o{X}-X\to \o{X}$, and where $\o{X}$ is any compactification of $X$. In particular, if $X$ is smooth projective, then $[X]$ is sent to the class of $X$ viewed as a constant simplicial variety. 
\end{corollary}
\begin{proof}
    By Theorem \ref{thm:proobj}, the inclusion of constant pro-objects $\s\V_*^{\perf}\to \Pro^{wc}(\s\V_*^{\perf})$ induces a homotopy equivalence $K(\s\V_*^{\perf})\simeq K(\Pro^{wc}(\s\V_*^{\perf}))$. Composing the map $K(\mathcal{V})\to K(\Pro^{wc}(\s\V_*^{\perf}))$ of Theorem \ref{thm:liftmeasure} with the inverse of this equivalence yields the desired map $K(\mathcal{V})\to K(\s\V_*^{\perf})$. The description at the level of $K_0$ follows from the description on objects of the weak $W$-exact functor in Theorem \ref{thm:liftmeasure} together with Theorem \ref{thm:Cam}. For the last claim, note that a smooth projective variety $X$ is its own compactification, so we may take $j_{X}:\O\to X$, and then $sCone(j_{X})=X$ viewed as a constant simplicial object of $\s\V_*^{\perf}$. 
\end{proof}
\begin{corollary}\label{cor:derivedGS}
Composing the map of Corollary \ref{cor:sVmotivicmeasure} with the map on K-theory spectra of Definition \ref{def:Chow}, we obtain a map:
\[K(\mathcal{V})\to K(\Ch^{\perf}(Chow))
\]
which on the level of $K_0$ sends the class $[X]$ of a variety to the class $[\Phi_*(sCone(\tilde{j}_X))]$ of any of its weight complexes. By Definition \ref{def:premotive}, this is precisely the class of the Chow motive of $X$ in the sense of Gillet--Soul\'{e}, and so we call this map the \emph{derived Gillet--Soul\'{e} measure}.
\end{corollary}
\section{The Derived Measure To Simplicial cdh Sheaves}
We must recall some basic definitions used in the $\A^1$-homotopy theory of Morel and Voevodsky. Consider the cdh topology on $\Sch_k$, the category of Noetherian schemes over $k$. This is defined in 1.2.2 of \cite{Vo2}. Recall that envelopes were introduced in Proposition \ref{prop:envelope}; by 2.17 of \cite{Vo}, envelopes are coverings in the cdh topology. Another important fact about cdh sheaves is the following:
\begin{proposition}[1.2.6 of \cite{Vo2}]\label{prop:cdhblowup}
An \emph{abstract blow-up square} is a pullback square of the form:
\begin{equation}\begin{tikzcd}
B\arrow[d]\arrow[r]&Y\arrow[d,"p"]\\
A\arrow[r,"e"]&X
\end{tikzcd}
\end{equation}
where $p$ is proper, $e$ is a closed embedding and $p^{-1}(X-e(A))\to X-e(A)$ is an isomorphism. Every cdh sheaf $F$ sends abstract blow-up squares to pullback diagrams; that is, the natural map
\[F(X)\to F(A)\times_{F(B)}F(Y)
\]
is an isomorphism. 
\end{proposition}
The following proposition follows from 6.17 of \cite{HuKe}:
\begin{proposition}\label{prop:subcanonical}
    The $cdh$ topology on $\Sch_k$ restricted to $Sm_k$, the subcategory of smooth varieties over $k$, is sub-canonical, i.e. representable presheaves of smooth schemes are cdh-sheaves.  
\end{proposition}

\begin{lemma}\label{lem:Yonedapushouts}
Consider the Yoneda embedding $r:\s\V_*^{\perf}\to \sShv_{\cdh}(\Sch)_*$ sending a simplicial smooth projective variety $X_\bullet$ to its representable $r(X_\bullet)=(rX)_\bullet$ (which is a simplicial cdh sheaf by Proposition \ref{prop:subcanonical}). Then, $r$ preserves coproducts.
\end{lemma}
\begin{proof}
Let $A,B\in \s\V_*^{\perf}$ denote two simplicial objects concentrated in degree $0$. For any $cdh$-sheaf $\mathcal{F}$, we have the natural isomorphism in the presheaf category $\Pre(\Sch)_*$ via Yoneda:
\begin{align*}Hom_{\Pre(\Sch)_*}(r(A\coprod B),\mathcal{F})&=\mathcal{F}(A\coprod B)\\
&\cong\mathcal{F}(A)\times \mathcal{F}(B)\quad\quad\text{(the decomposition $A\coprod B$ is an abstract blow-up square, so Proposition \ref{prop:cdhblowup} applies)}\\
&\cong Hom_{\Pre(\Sch)_*}(rA,\mathcal{F})\times Hom_{\Pre(\Sch)_*}(rB,\mathcal{F})\\
&\cong Hom_{\Pre(\Sch)_*}(rA\coprod rB,\mathcal{F})
\end{align*}
The sheafification functor $\#:\Pre(\Sch)_*\to Shv_{\cdh}(\Sch)_*$ preserves finite products and colimits, and is left adjoint to the inclusion functor $Shv_{\cdh}(\Sch)_*\hookrightarrow \Pre(\Sch)_*$, further giving us the natural isomorphisms:
\[Hom_{\Pre(\Sch)_*}(r(A\coprod B),\mathcal{F})\cong Hom_{Shv_{\cdh}(\Sch)_*}(r(A\coprod B),\mathcal{F})\\
\]
and
\[
Hom_{\Pre(\Sch)_*}(rA\coprod rB,\mathcal{F})\cong Hom_{Shv_{\cdh}(\Sch)_*}((rA)\coprod (rB),\mathcal{F})
\]
In other words, we have the natural isomorphism for any sheaf $\mathcal{F}$:
\[Hom_{Shv(\Sch)_*}\Big((rA)\coprod (rB),\mathcal{F}\Big)\cong Hom_{Shv(\Sch)_*}\Big(r(A\coprod B),\mathcal{F}\Big)
\]
We can therefore conclude that in $\sShv_{\cdh}(\Sch)_*$ \[r(A\coprod B)=r(A)\coprod r(B)\]
i.e. the functor preserves coproducts on constant simplicial objects. Since coproducts of simplicial objects are taken degree-wise, we have that the functor on all of $\s\V_*^{\perf}$ must also preserve coproducts.
\end{proof}
 
In order to turn $Shv_{Nis}(\Sch)_*$ into a Waldhausen category, we recall the following proposition, Theorem 5 of \cite{Jar}, attributed to Joyal:
\begin{proposition}\label{prop:joyal}
Let $f:\mathcal{F}\to\mathcal{G}$ be a morphism of simplicial sheaves. $\sShv_{\cdh}(\Sch)_*$ has a cofibrantly generated proper model structure where:
\begin{itemize}
    \item $f$ is a weak equivalence if it induces an isomorphism on homotopy sheaves
    \item $f$ is a cofibration if it is a degree-wise monomorphism
    \item $f$ is a fibration if it has the right lifting property with respect to any cofibration that is a weak equivalence
\end{itemize}
\end{proposition}
\begin{definition}
We define a simplicial sheaf $\mathcal{F}_\bullet\in \sShv_{\cdh}(\Sch)_*$ to be perfect if there exists $\tilde{G}_\bullet$ that can be obtained from taking simplicial mapping cones and coproducts or products of constant representable simplicial sheaves such that $\mathcal{F}_\bullet$ and $\tilde{G}_\bullet$ are isomorphic in the homotopy category $\Ho(\sShv_{\cdh}(\Sch)_*)$.
\end{definition}
\begin{definition}
Let $\perf_{\mathrm{co}}(\sShv_{\cdh}(\Sch)_*)$ denote the subcategory of perfect and cofibrant objects in $\sShv_{\cdh}(\Sch)_*$. 
\end{definition}
\begin{proposition}
    The classes of cofibrations and weak equivalences in $\perf_{\mathrm{co}}(\sShv_{\cdh}(\Sch)_*)$ give it a Waldhausen structure with FFWC 
\end{proposition}
\begin{proof}
The proof that $\perf_{\mathrm{co}}(\sShv_{\cdh}(\Sch)_*)$ is Waldhausen follows immediately from the model category axioms, and the FFWC portion of the proof is the same as in Lemma 1 of \cite{BGN}, which relies on the fact that the model structure is cofibrantly generated.
\end{proof}

\begin{proposition}\label{prop:sshv}
The Yoneda embedding $r:\s\V_*^{\perf}\hookrightarrow \perf_{\mathrm{co}}(\sShv_{\cdh}(\Sch)_*)$ is an exact functor of Waldhausen categories, where $\s\V_*^{\perf}$ has the Waldhausen structure given in Definition \ref{def:we}.
\end{proposition}
\begin{proof}
 By Lemma \ref{lem:Yonedapushouts}, the functor preserves level-wise inclusions onto direct summands of coproducts, so since any cofibration $f_\bullet$ in $\s\V_*^{\perf}$ is of the form $f_n:X_n\hookrightarrow X_n\coprod Y_n$, we see that $rf_\bullet$ will also be such that $rf_n:rX_n\hookrightarrow rX_n\coprod rY_n$, which is a level-wise monomorphism in $\sShv_{\cdh}(\Sch)_*$. So, $r$ preserves cofibrations. A similar argument shows that level-wise pushouts along cofibrations in $\s\V_*^{\perf}$ are preserved. In the case of weak equivalences, isomorphisms and simplicial homotopy equivalences in $\s\V_*^{\perf}$ will be isomorphisms and homotopy equivalences respectively in $\s\Pre(\Sch)_*$ just by functoriality. Lemma 1 of \cite{Jar2} shows that all hypercovers are local weak equivalences in $\s\Pre(\Sch)_*$ given the same model structure as in Proposition \ref{prop:joyal}, and so they will remain local weak equivalences after sheafification. Now, consider a diagram of the form:
 \[\begin{tikzcd}
 A_\bullet\arrow[r,"f"]\arrow[d,"h_A"]&B_\bullet\arrow[r]\arrow[d,"h_B"] &sCone[f]\arrow[d]\\
 A_\bullet'\arrow[r,"g"]&B_\bullet'\arrow[r]& sCone[g]
 \end{tikzcd}
\]
where $h_A,h_B$ are hyperenvelopes or simplicial homotopy equivalences. We already know that $rh_A,rh_B$ are weak equivalences of simplicial sheaves, and therefore so must be the induced morphism on cones. Lastly, consider a cubical diagram \[
    \begin{tikzcd}[row sep=2.2em, column sep=2.2em]
    X_\bullet \arrow[rr,"f_\bullet"] \arrow[dr] \arrow[dd,"h_X"] &&
    Y_\bullet \arrow[dd,"h_Y",swap] \arrow[dr] \\
    & X_\bullet'\arrow[rr,"f_\bullet' "]\arrow[dd,"h_{X'}",swap] &&
    Y_\bullet' \arrow[dd,"h_{Y'}",swap] \\
     X \arrow[rr] \arrow[dr] && Y\arrow[dr] \\
    & X' \arrow[rr] &&Y'
    \end{tikzcd}
\]
where the bottom square is an abstract blow-up square of proper varieties and the vertical morphisms are hyperenvelopes. The map $Y\to Y'$ is a cdh cover, and so the induced map on the sheafification $rY\to rY'$ is an effective epimorphism by  \cite[\href{https://stacks.math.columbia.edu/tag/00WT}{Tag 00WT}]{stacks-project} and \cite[\href{https://stacks.math.columbia.edu/tag/086K}{Tag 086K}]{stacks-project}. This means that $rY\times_{rY'}rY\rightrightarrows rY'$ is a coequalizer, so we have a diagram of sheaves:
\[
\begin{tikzcd}
    rX\arrow[rr]\arrow[rd,hook]\arrow[dd]&&rY\arrow[dd]\\
    &rY\times_{rY'}rY\arrow[dr]\arrow[dl]\arrow[ur]\arrow[dr,shift left=2]&\\
    rX'\arrow[rr]&&rY'
\end{tikzcd}
\]
where every path from $rX\to rY'$ commutes, implying that the outer square is a pushout square. Since the map $rX\to rY$ is a level-wise monomorphism of sheaves (as Yoneda is limit preserving) the square is a homotopy pushout square. By 13.5.9 of \cite{Hirschorn} the upper square
\[
\begin{tikzcd}
    X_\bullet\arrow[r,"f_\bullet"]\arrow[d]&Y_\bullet\arrow[d]\\
    X_\bullet'\arrow[r,"f_\bullet'"]&Y'_\bullet
\end{tikzcd}
\]
is a homotopy cofiber square, and so the induced map on homotopy cofibers is a weak equivalence. The homotopy cofiber in $\sShv_{\cdh}(\Sch)_*$ of a morphism between representables of smooth projective varieties will be the same as the image of the simplicial mapping cone of a morphism in $\s\V_*^{\perf}$, meaning that $r(sC_{f_\bullet})\to r(sC_{f_\bullet'})$ will be a weak equivalence. So, Yoneda sends morphisms in $w(\s\V_*^{\perf})$ to weak equivalences in $\sShv_{\cdh}(\Sch)_*$. Now, consider the subcategory of $\s\V_*^{\perf}$ whose morphisms land in weak equivalences in $\sShv_{\cdh}(\Sch)_*$. Since Yoneda preserves pushouts along cofibrations and weak equivalences are saturated in $\sShv_{\cdh}(\Sch)_*$, this subcategory satisfies the two axioms of Definition \ref{def:we}, implying it contains $\textbf{w}(\s\V_*)$. So, we can actually conclude that Yoneda preserves all weak equivalences. Since Yoneda additionally preserves mapping cones and coproducts, it will by definition send perfect objects of $\s\V_*$ to perfect objects of $\sShv_{\cdh}(\Sch)_*$. Lastly, the image of $\s\V_*^{\perf}$ via Yoneda lands in cofibrant objects, as all objects are cofibrant in the model structure outlined in Proposition \ref{prop:joyal}.
\end{proof}

\begin{corollary}\label{cor:cdhmeasure}
Writing $K(\sShv_{\cdh}(\Sch)_*)$ as the Waldhausen $K$-theory of perfect and cofibrant objects of $\sShv_{\cdh}(\Sch)_*$, we have a motivic measure
\[K(\mathcal{V}_k)\to K(\sShv_{\cdh}(\Sch)_*)
\]
On the level of $K_0$, this map will be $[X]\mapsto [r(\tilde{j}_X)]$, where $[r(\tilde{j}_X)]$ is the isomorphism class of the representable presheaf associated to any weight complex of $X$.  
\end{corollary}
\begin{corollary}
    At the level of $K_0(\sShv_{\cdh}(\Sch)_*)$, for any proper variety $\o{X}$, the map in Corollary \ref{cor:cdhmeasure} sends:
    \[[\o{X}]\mapsto [(r\o{X})^{\#}]\] where $(r\o{X})^{\#}$ is the sheafification of the representable presheaf of $\o{X}$. It therefore follows by Corollary \ref{cor:cdhmeasure} that for any variety $X$, the isomorphism class in $\sShv_{\cdh}(\Sch)_*$ of the weight complex can be computed as
    \[[r(sC_{\tilde{j}_X})]=[(r\o{X})^{\#}]+[r(\o{X}-X)^{\#}]\]
    where $\o{X}$ is any compactification of $X$.
\end{corollary}
\begin{proof}
    When building the weight complex of a proper variety $\o{X}$, it serves as its own compactification, and so any smooth projective hyperenvelope $\tilde{X}_\bullet$ of $\o{X}$ is a weight complex of $\o{X}$ in $\s\V_*^{\perf}$. Further, the map $\tilde{X}_\bullet\to \o{X}$ is a cdh hypercover in $\Sch$ and also in $\sShv_{\cdh}(\Sch)_*$, implying that in $K_0(\sShv_{\cdh}(\Sch)_*)$, we have $[r\tilde{X}_\bullet]=[r(\o{X})^{\#}]$. This proves that the motivic measure in Corollary \ref{cor:cdhmeasure} sends the isomorphism class of a proper variety $[\o{X}]$ to the isomorphism class of its corresponding cdh sheaf $[r(\o{X})^{\#}]$. The rest of the corollary follows.
\end{proof}
\section{Recovering the Derived Voevodsky Measure}

Recall that \cite{BGN} constructs a derived motivic measure to  $K(\mathrm{C}^{\mathrm{gm}}(R))$, where \[\mathrm{C}^{\mathrm{gm}}(R)=\perf_{\mathrm{co}}((\mathcal{D}\mathcal{M}^{\eff}_{gm,Nis})^{op}-Rmod)\] gives us the $K$-theory of the DG category of perfect effective geometric Voevodsky motives with coefficients in a commutative ring $R$, i.e. $\mathcal{D}\mathcal{M}^{\eff}_{gm,Nis}(R)$ constructed in \cite{BV}. We first review the notion of both the usual and compactly supported Voevodsky motive in loc. cit.
\begin{definition}
For a commutative ring $R$, if $X$ and $Y$ are varieties over a field $k$, define:
\begin{align*}
\textit{Cor}_R(X,Y):=&R\l \Gamma\subset Y\times X,\  \Gamma\text{ is a reduced, irreducible subscheme},\\
&\pi_2:\Gamma\to X\text{ finite, dominant over a connected component}\r
\end{align*}
Here, $R\l S \r$ denotes the free $R$ module generated by the set $S$. 
Let $Sm$ denote the category of smooth varieties and $Cor_R(Sm)$ denote the category with the same objects, but the morphisms are the elements of $\textit{Cor}_R(X,Y)$. Note that there is an associative composition in $Cor_R(Sm)$ that is defined the same way as in $Chow$. 
\end{definition}
\begin{proposition}
There is an inclusion functor $\V\to Cor_R(Sm)$, which is the identity on objects and where $f:X\to Y$ in $\V$ will be sent to the graph $\Gamma_f\in Cor_R(X,Y)$.
\end{proposition}
\begin{definition}
Let $Cor_R(Sm)^{op}-\mathrm{R-mod}$ be the category of pre-sheaves of $R$-modules on $Cor_R(Sm)$, and $Cor_R(Sm)^{op}\mathrm{-dgm}$ similarly denote the category of pre-sheaves of DG $R$-modules on $Cor_R(Sm)$. Observe that $Cor_R(Sm)^{op}-\mathrm{mod}$ is a full subcategory of $Cor_R(Sm)^{op}\mathrm{-dgm}$.
\end{definition}
\begin{definition}
For $X\in ob(Cor_R(Sm))$, we let $R_{tr}[X]:=Hom_{Cor_R(Sm)}(-,X)$ denote its representable DG-presheaf, i.e. in $Cor(Sm)^{op}-\mathrm{R\mathrm{-dgm}}$ (it is actually an object of $Cor_R(Sm)^{op}-\mathrm{R-mod}$).
\end{definition}

\begin{definition}
 We also define
\begin{align*}
QDom_R(X,Y):=&R\l \Gamma\subset Y\times X,\  \Gamma\text{ is a reduced, irreducible subscheme},\\
&\pi_2:\Gamma\to X\text{ quasi-finite, dominant over a connected component}\r
\end{align*}
\end{definition}

\begin{proposition}
There is a functor \[R_{tr}^c:\mathcal{V}^{prop}\to Cor_R(Sm)^{op}-mod\]
Here, $\mathcal{V}^{prop}$ is the category where the objects are k-varieties and the morphisms are only proper ones. For a variety $X$, $R_{tr}^c[X]$ is a pre-sheaf sending, on objects
\[Y\mapsto  QDom_R(Y,X)
\]
On morphisms $R_{tr}^c[X]$ sends a finite correspondence $\Gamma\in Cor_R(Y,Z)$ to the map $ QDom_R(Z,X) \to QDom_R(Y,X)$ induced by $ - \circ\Gamma$, and we extend the map $R$-linearly. Since $Y$ is smooth, the composition map is well-defined. Further, $R_{tr}^c$ sends morphisms $g:X\to Z$ to the natural transformation $\eta^{g}:R_{tr}^c[X]\to R_{tr}^c[Z]$ such that for each $Y\in Cor_R(Sm)^{op}$
\[\eta^{g}_{Y}: QDom_R(Y,X)\to  QDom_R(Y,Z)
\]
is induced by the standard pushforward of algebraic cycles along the proper map $g$\begin{footnote}{In \cite{BGN}, $R^c_{tr}$ is referred to as $z_{equi}(-,0)$, which is the notation in \cite{MVW}}\end{footnote}.
\end{proposition}
\begin{proposition}\label{prop:mdgc=mdg}
If $X$ is a proper variety, then for any smooth variety $Y$, we have $QDom_R(Y,X)=Cor_R(Y,X)$. Therefore, $R_{tr}[X]=R_{tr}^c[X]$. 
\end{proposition}
\begin{proof}
    Given a proper variety $X$ and a closed integral subscheme $\Gamma\subset Y\times X$, we have that the projection $\Gamma\to X$ is closed in $X$, making it a proper map. Further, if the projection is quasi-finite and dominant on a connected component of $X$, then it is in fact finite (see \href{https://stacks.math.columbia.edu/tag/02OG}{Tag 02OG}). This implies $QDom_R(Y,X)=Cor_R(Y,X)$. 
\end{proof}
The following proposition is the content of the material following 6.9.3 of \cite{BV}:
\begin{proposition}\label{prop:Mcdg}
Let $\mathcal{V}^{prop}$ denote the category of varieties where we consider only proper morphisms between them. The cdh topology on $\Sch$ will induce a cdh topology on $Cor_R(
\Sch)$, and we restrict our attention to presheaves on $Cor_R(Sm)\subset 
Cor_R(\Sch)$, defining $\mathcal{P}_{tr,cdh}$ to be the DG enhancement of the derived category of cdh sheaves restricted to $Cor_R(Sm)$. There exists a cdh sheafification functor \[C_{Cor_R(Sm)}:(Cor_R(Sm)^{op}-\mathrm{sR-mod})\to \mathcal{P}_{tr}\] along with a dg-$\A^1$-localization functor, $C^{\Delta}:\mathcal{P}_{tr}\to \mathcal{P}_{tr}$, such that if we write $C^{\mathcal{M}}:=C^{\Delta}\circ C_{Cor(Sm)}$, the composition
\[C^{\mathcal{M}}\circ R_{tr}^c:\mathcal{V}^{prop}\to \mathcal{P}_{tr}
\]
is a DG lift of the assignment of the Voevodsky motive with compact support $X\mapsto M^c_{tr}(X)$ in \cite{V}. Define $M^c_{dg}:=C^{\mathcal{M}}\circ R_{tr}^c$ One can extend $M^c_{dg}$ to simplicial proper varieties, i.e. \[sM^c_{dg}:s\mathcal{V}^{prop}\to \mathcal{P}_{tr}\] by normalizing the associated simplicial presheaf with transfers and then applying $C^{\mathcal{M}}$.
\end{proposition}
\begin{definition}
We define $DM_{gm,cdh}^{\eff}(R)$ objects of $Cor_R(Sm)^{op}-\mathrm{R\mathrm{-dgm}}$ quasi-isomorphic to bounded complexes of objects in the essential image of the functor $M^c_{dg}$ of $Cor_R(Sm)^{op}-\mathrm{R\mathrm{-dgm}}$.
\end{definition}
The next proposition shows us that the cdh topology does not change the category of motives:
\begin{proposition}[Proposition 2 of \cite{BV}]\label{prop:cdhtoNis}
    As Nisnevich covers are cdh covers, the `inclusion' functor: $I:DM_{gm,cdh}^{\eff}(R)\to DM_{gm,Nis}^{\eff}(R)$ is a DG functor and quasi-equivalence (i.e. the induced morphism of triangulated categories is an equivalence).
\end{proposition}
We conclude our review by citing an important proposition.
\begin{proposition}[\cite{BV} - Lemma 6.9.1]\label{prop:NMV}
For any scheme $X$ and $U\subset X$:
\[\begin{tikzcd}
M^c_{dg}(X-U)\arrow[r]\arrow[d] &M^c_{dg}(X)\arrow[d]\\ 0\arrow[r]&M^c_{dg}(U)
\end{tikzcd}
\]
is quasi-isomorphic to a homotopy pushout square in $\mathcal{D}_{gm,Nis}^{\eff}$
\end{proposition}
Now, we can build a map of $K$-theory spectra:
\begin{proposition}\label{prop:Mdg}
The functor $M_{dg}:\s\V_*^{\perf}\to \mathrm{C}^{\mathrm{gm}}(R)$ is weakly exact, giving us a map:
\[K(\s\V_*^{\perf})\to K(\mathrm{C}^{\mathrm{gm}}(R))
\]
where $K(\mathrm{C}^{\mathrm{gm}}(R))$ is defined as in Definition \ref{def:kthydg}
\end{proposition}
\begin{proof}
All morphisms in both categories are weak cofibrations by Corollary \ref{prop:ffwcsV} and Proposition \ref{prop:dgffwc}.
It is quick to check that simplicial homotopy equivalences will be global quasi-isomorphisms in $DM^{\eff}_{gm,Nis}(R)$ and therefore weak equivalences in $\mathrm{C}^{\mathrm{gm}}(R)$. As for a hyperenvelope $X_\bullet\to Y_\bullet$, we need to show that the image $R_{tr}[X_\bullet]\to R_{tr}[Y_\bullet]$ in $DM^{\eff}_{gm,Nis}(R)$ is a global quasi isomorphism. Now, $\mathcal{P}_{tr,cdh}(R)$ is defined to be the DG enhancement of the derived category of cdh sheaves restricted to $Cor_R(Sm)$, in which the image of all cdh hypercoverings are global quasi-isomorphisms (see the lemma in 1.1 of \cite{BV}). We see that after applying $C^{\Delta}$ it will remain a global quasi-isomorphism in $DM^{\eff}_{gm,cdh}(R)$, and due to the equivalence on their triangulated categories, will also be a quasi-isomorphism as an object of $DM^{\eff}_{gm,Nis}(R)$ (and therefore a weak equivalence in $\mathrm{C}^{\mathrm{gm}}(R)$). Further, $C^{\mathcal{M}}\circ R_{tr}$ preserves coproducts, as for any varieties $X,Y$ we have that \[0=R_{tr}[0]\to R_{tr}[X]\oplus R_{tr}[Y]\to R_{tr}[X\coprod Y]\] is a 3-term Zariski Mayer-Vietoris complex which is contained in $\mathcal{I}^{Nis}$ (see \cite{BV} 4.4). As a result, we have an isomorphism \[M_{dg}(X\coprod Y)\cong M_{dg}(X)\oplus M_{dg}(Y)\] 
So, $M_{dg}$ will preserve level-wise inclusions onto direct summands, and therefore pushouts along such inclusions, meaning that pushouts along cofibrations are preserved. This further implies that weak equivalences of type (3) in Definition \ref{def:we1} will also be preserved. Next, since products and coproducts are the same in $\mathcal{D}_{gm,Nis}^{\eff}(R)$, we see that coproducts will be preserved when we compose with the Yoneda map $\mathcal{D}_{gm,Nis}^{\eff}(R)\to \mathrm{C}^{\mathrm{gm}}(R)$. Therefore, pushouts along level-wise inclusions will be preserved by $Y\circ M_{dg}$, so we have not only shown that pushouts along cofibrations are preserved, but that morphisms of the form (3) of Definition \ref{def:we1} will also be preserved. As $\mathrm{C}^{\mathrm{gm}}(R)$ inherits a Waldhausen structure from its model structure, the same argument as in the proof of Proposition \ref{prop:sshv} will tell us that morphisms of the form (4) of Definition \ref{def:we1} will also be sent to weak equivalences.
So we have that $w(\s\V_*^{\perf})$ will be sent to weak equivalences in $\mathrm{C}^{\mathrm{gm}}(R)$. The subcategory of $\s\V_*^{\perf}$ whose morphisms are weak equivalences in $\mathrm{C}^{\mathrm{gm}}(R)$ satisfies the two axioms of Definition \ref{def:we} and contains $w(\s\V_*^{\perf})$, so in fact $\textbf{w}(\s\V_*^{\perf})$ is also sent to weak equivalences. 
\end{proof}
\begin{corollary}
Composing the map in Corollary \ref{cor:sVmotivicmeasure} with that of Proposition \ref{prop:Mdg}, we get another derived measure:
\[K(\mathcal{V})\to K(\mathrm{C}^{\mathrm{gm}}(R))
\]
\end{corollary}
Now, we compare this motivic measure to the one on cdh sheaves. The inclusion $Sm\hookrightarrow \Sch$ gives a direct image functor: $f_*:\sShv_{\cdh}(\Sch)_*\hookrightarrow \sShv_{\cdh}(Sm)_*$, where $\sShv_{\cdh}(Sm)_*$ denotes the category of simplicial cdh sheaves that are colimits of representables of $Sm$.
Further, if $R[Sm]$ denotes the universal $R$-linear category associated to $Sm$, i.e. where $Hom_{R[Sm]}(X,Y)=R\l Hom_{Sm}(X,Y)\r$, we have the `free' functor:
\[\sShv_{\cdh}(Sm)_*\to R[Sm]^{op}-smod\]
Lastly, by the remark after (2.1.3) of \cite{BV}, we have a DG functor:
$\Psi:R[Sm]^{op}-smod\to \mathcal{P}_{tr}$ (left adjoint to the `forgetting transfers' functor) fitting into the diagram:
\[
\begin{tikzcd}
&R[Sm]^{op}-mod\arrow[d,"\Psi"]\\
Sm\arrow[ur]\arrow[r,"R_{tr}"]&\mathcal{P}_{tr}
\end{tikzcd}
\]
where the functors out of $Sm$ are Yoneda embeddings. Let us now denote $\Phi:\sShv_{\cdh}(Sm)_*\to \mathcal{P}_{tr}$ to be the composition of the free functor with $\Psi$. Using these functors, we can now construct a map of $K$-theory spectra: 
\begin{theorem}\label{thm:doldkan}
Using the notation of Propositions \ref{prop:Mcdg} and \ref{prop:cdhtoNis}
\[I\circ C^{\mathcal{M}}\circ\Phi \circ f_*:\sShv_{\cdh}(\Sch)_*\to DM_{gm,Nis}^{\eff}(R)\]
is weakly exact between two Waldhausen categories with FFWC, giving us a map of $K$-theory spectra:
\[K(\sShv_{\cdh}(\Sch)_*)\to K(\mathrm{C}^{\mathrm{gm}}(R))\]
\end{theorem}
\begin{proof}
We first check that weak equivalences are preserved. In $\sShv_{\cdh}(\Sch)_*$, they are given as local weak equivalences $\mathcal{F}\to\mathcal{G}$. These remain local weak equivalences under the direct image functor $f_*$. Next, a local equivalence $\mathcal{F}\to\mathcal{G}$ gives a local weak equivalence of sheaves of simplicial abelian groups $R\mathcal{F}\to R\mathcal{G}$, and $\Psi$ will send this to an isomorphism of cohomology sheaves. By the discussion in 1.11 of \cite{BV} this is a locally acyclic map of DG sheaves, which is an isomorphism in the derived category of DG sheaves of modules. Then $I\circ \mathcal{C}^\mathcal{M}$ is a DG functor, so it will preserve this isomorphism on homotopy categories, i.e. it will be a weak equivalence in $\mathrm{C}^{\mathrm{gm}}(R)$. We now need to show that pushouts along cofibrations are sent to \emph{weak} pushouts along cofibrations (which is any pushout in $\mathrm{C}^{\mathrm{gm}}(R)$ as all morphisms are weak cofibrations). As the direct image functor is a right adjoint and is computed as a filtered limit of simplicial sets, it will preserve monomorphisms and finite colimits, i.e. all pushouts along monomorphisms. Next, since the free functor is a left Quillen map, it must also send pushouts along cofibrations to pushouts. Lastly, DG functors preserve finite colimits, which allows us to conclude the second condition for weak exactness. 
\end{proof}

We have the following immediate corollary:
\begin{corollary}
The weakly exact functors in Proposition \ref{prop:sshv}, Proposition \ref{prop:Mdg} and \ref{thm:doldkan} commute, therefore giving us a commutative diagram of spectra:
\[\begin{tikzcd}
K(\s\V_*^{\perf})\arrow[dr]\arrow[r]&K(\sShv_{\cdh}(Sm)_*)\arrow[d]\\
& K(\mathrm{C}^{\mathrm{gm}}(R))
\end{tikzcd}
\]
\end{corollary}

Recall that \cite{BGN} had constructed a more direct measure which we recall here in the terminology of \cite{BV}:
\begin{proposition}\label{BGN}
There is a weak $W$-exact functor  $(F^w,F^!,F_!):\mathcal{V}\to \mathrm{C}^{\mathrm{gm}}(R)$ as follows: on objects, we directly send
\[ X\mapsto M_{dg}^c[X]
\]
On morphisms:
\begin{enumerate}
    \item For a cofibration (closed embedding) $f:X\hookrightarrow Y$, we simply define $F_!(f)=C^{\mathcal{M}}\Big(R_{tr}^c[f]\Big)$, which is well-defined since closed immersions are proper. This construction is covariant, as proper maps define pushforwards. 
    \item For a weak equivalence, i.e. isomorphism $f:X\xrightarrow{\cong}Y$, we just set $F^w(f)=F_!(f)$. 
    \item Finally, for a complement (open embedding) $f:X\xrightarrow{\circ} Y$, we define the `restriction' morphism, or $C^{\mathcal{M}}$ composed with the natural transformation $T_{f}^*:R_{tr}^c[Y]\to R_{tr}^c[X]$ such that on objects $U\in Cor_R(Sm)^{op}$ and correspondences $\Gamma\in QDom_R(U,Y)$ (i.e. not an $R$-linear combination of correspondences) the map is given by the pullback induced by an open immersion $X\times U\to Y\times U$\begin{footnote}{As noted in \cite{BV}, such a pullback can actually be given for any flat and quasi-finite morphism, not just an open immersion.}\end{footnote}:
    \[T_{f}^*:[\Gamma]\mapsto [(f\times id_U)^{*}(\Gamma)]\in QDom_R(U,X)
    \]
    and we then extend $R-$linearly to define the map on all elements of $QDom_R(U,X)$. Recall that $(f\times id_U)^{*}$ is explicitly given by 
    \[[\Gamma]\mapsto [(f^{-1}\times id_U)(\Gamma\cap (f(X)\times U))]=[(f^{-1}\times id_U)(\Gamma)]
    \]
     Here $[\Gamma]$ denotes the closed subvariety associated to the closed subset $\Gamma\subset Y\times U$, i.e. $(\Gamma,\mathcal{O}|_{Y\times U,\Gamma})$, where $\mathcal{O}|_{Y,\Gamma}$ is the structure sheaf localized at the closed subset $\Gamma$. So, we can in turn define \[F^!(f)=C^{\mathcal{M}}\Big(T_{f}^*\Big):M^c_{dg}(Y)\to M^c_{dg}(X)\] 
\end{enumerate}
\end{proposition}
Now, we would like to compare the two different K-theory maps from $K(\mathcal{V}_k)\to K(\mathrm{C}^{\mathrm{gm}}(R))$. However, since the one in Proposition \ref{prop:Mdg} is achieved via the route of pro-objects, we must compare these maps in the category of weakly constant pro-objects of $\mathrm{C}^{\mathrm{gm}}(R)$, i.e. $\Pro^{wc}(\mathrm{C}^{\mathrm{gm}}(R))$. The following proposition follows just as in Theorem \ref{thm:proobj}:
\begin{proposition}
Let $\Pro^{wc}(\mathrm{C}^{\mathrm{gm}}(R))$ denote the subcategory of $\Pro(\mathrm{C}^{\mathrm{gm}}(R))$ consisting of pro-objects such that all the morphisms within the pro-object are weak equivalences. Just as in Propositions \ref{prop:prowald} and \ref{prop:proffwc}, $\Pro^{wc}(\mathrm{C}^{\mathrm{gm}}(R))$ is Waldhausen with FFWC, where a morphism $f:X\to Y$ is a weak equivalence (resp. cofibration) if one can find a representation of $f$ indexed by $Y$ by weak equivalences (resp. cofibrations) in $\mathrm{C}^{\mathrm{gm}}(R)$. 
\end{proposition}
\begin{proposition}
The weakly exact functor of Waldhausen categories in Proposition \ref{prop:Mdg} also gives a weakly exact functor:
\[\Pro^{wc}(\s\V_*^{\perf})\to \Pro^{wc}(\mathrm{C}^{\mathrm{gm}}(R))
\]

\end{proposition}
\begin{proof}
The weakly exact functor in Proposition \ref{prop:Mdg} $\s\V_*^{\perf}\to \mathrm{C}^{\mathrm{gm}}(R)$ extends to a functor \[\Pro(M_{dg}): \Pro^{wc}(\s\V_*^{\perf})\to \Pro^{wc}(\mathrm{C}^{\mathrm{gm}}(R))\] on pro-objects as it preserves weak equivalences. Further if $f:X\to Y$ is a morphism in $\Pro^{wc}(\s\V^{\perf}_*)$ and one can find a representation of $f$ indexed by $Y$ by weak equivalences (resp. cofibrations) in $\s\V_*^{\perf}$, then $\Pro(M_{dg})$ will also send these morphisms to weak equivalences (resp. weak cofibrations) in $\mathrm{C}^{\mathrm{gm}}(R)$ via Proposition \ref{prop:Mdg}, which tells us that $\Pro(M_{dg})$ preserves weak cofibrations and weak equivalences. By the same logic, $\Pro(M_{dg})$ preserves level-wise inclusions onto direct summands just like $M_{dg}$, and therefore sends pushouts along cofibrations (which are inclusions into coproducts) to pushouts along weak cofibrations (i.e. weak pushout squares) concluding the proof.  
\end{proof}
The proof of the next proposition is the same as in Theorem \ref{thm:proobj}:
\begin{proposition}
The inclusion $\mathrm{C}^{\mathrm{gm}}(R)\hookrightarrow \Pro^{wc}(\mathrm{C}^{\mathrm{gm}}(R))$ gives an equivalence:
\[K(\mathrm{C}^{\mathrm{gm}}(R))\xrightarrow{\simeq} K(\Pro^{wc}(\mathrm{C}^{\mathrm{gm}}(R)))\]
\end{proposition}
We can now compare the derived measures. For a variety $X$, if we choose a compactification $\o{X}$ of $X$, and a smooth projective hyperenvelope $\tilde{j}_X$ of $j_X:\o{X}-X\to\o{X}$, we get the following diagram of pre-sheaves with transfers:
\[
\begin{tikzcd}
M^{c}_{dg}(\widetilde{\o{X}-X}_\bullet)\arrow[r,"M^{c}_{dg}(\tilde{j}_X)"]\arrow[d,"\simeq"]& M^{c}_{dg}(\tilde{X}_\bullet)\arrow[r]\arrow[d,"\simeq"]&M^{c}_{dg}(Cone(\tilde{j}_X))=Cone(M^{c}_{dg}[\tilde{j}_X])\arrow[d]\\
M^{c}_{dg}(\o{X}-X)\arrow[r,"M^{c}_{dg}(j_X)"]&M^{c}_{dg}(\o{X})\arrow[r]& M^{c}_{dg}(X)
\end{tikzcd}
\]
Here, we recall Proposition \ref{prop:mdgc=mdg}, namely that $M^c_{dg}=M_{dg}$ for proper varieties. Above, as an abuse of notation, we equate a simplicial complex $X_\bullet$ with its Moore complex $M(X_\bullet)$. By Proposition \ref{prop:NMV}, \[0\to M^{c}_{dg}[\o{X}-X]\to M^{c}_{dg}[\o{X}]\to M^{c}_{dg}[X]\to 0\] is exact. Further, the lemma in Section 1.11 of \cite{BV} tells us that the vertical arrows are also weak equivalences. As a result, the induced maps
\[Cone(M^{c}_{dg}[j_X])\xrightarrow{\simeq}M^{c}_{dg}[X]
\]
and
\[Cone(M^{c}_{dg}[\tilde{j}_X])\xrightarrow{\simeq}Cone(M^{c}_{dg}[j_X])
\]
are weak equivalences. Composing them gives us a collection of weak equivalences, each of which we denote
\[\Phi_{X,\o{X},\tilde{j}_X}: M^{c}_{dg}(Cone(\tilde{j}_X))\xrightarrow{\simeq} M^{c}_{dg}(X)\] Let $G=(G_!,G^!,G_w)$ define the weak $W$-exact functor which gives the motivic measure \[K(\mathcal{V})\to K(\mathrm{C}^{\mathrm{gm}}(R))\to K(\Pro^{wc}(\mathrm{C}^{\mathrm{gm}}(R)))\] where the first map in the composition is defined in Proposition \ref{BGN}, and let $F=(F_!,F^!,F_w)$ define the weak $W$-exact functor which gives the motivic measure \[K(\mathcal{V})\to K(\Pro^{wc}(\s\V_*^{\perf}))\to K(\Pro^{wc}(\mathrm{C}^{\mathrm{gm}}(R)))\] where the first map in the composition is constructed in Theorem \ref{thm:liftmeasure}. We see that for each $X\in\mathcal{V}$, each $\Phi_{X,\o{X},\tilde{j}_X}$ defines a weak equivalence $F(X)\xrightarrow{\simeq} G(X)$, as $G(X)$ is just the constant pro-object $M^c_{dg}(X)$. A priori, it may be that different choices of compactifications/hyperenvelopes represent different weak equivalences between pro-objects. The next lemma shows that this is not the case.
\begin{lemma}
For a fixed $X\in\mathcal{V}$ and choosing any objects $(X,\o{X},\tilde{j}_X)$ in $PH(\mathcal{V})$, each morphism \[\Phi_{X,\o{X},\tilde{j}_X}: M^{c}_{dg}(Cone(\tilde{j}_X))\xrightarrow{\simeq} M^{c}_{dg}(X) \] represents the same morphism of pro-objects \[\Phi_X: F(X)\xrightarrow{\simeq} G(X)\]
which is in fact a weak equivalence.
\end{lemma}
\begin{proof}
Given  $(X,\o{X}_1,\tilde{j}_{X_1}),(X,\o{X}_2,\tilde{j}_{X_2})$ in $PH(\mathcal{V})$, one can check that the morphism $\Phi_{(X,\o{X}_1\times\o{X}_2,\tilde{j}_{X_1,X_2})}$ fits into the diagram:
\[
\begin{tikzcd}
&M^{c}_{dg}[Cone(\tilde{j}_{X_1,X_2})]\arrow[dr]\arrow[dd,"\Phi_{(X,\o{X}_1\times\o{X}_2,\tilde{j}_{X_1,X_2})}"]\arrow[dl]&\\
M^{c}_{dg}[Cone(\tilde{j}_{X_1})]\arrow[dr,"\Phi_{(X,\o{X}_1,\tilde{j}_{X_1})}",swap]&&M^{c}_{dg}[Cone(\tilde{j}_{X_2})]\arrow[dl,"\Phi_{(X,\o{X}_2,\tilde{j}_{X_2})}"]\\
&M^{c}_{dg}[X]&
\end{tikzcd}
\]
This shows that $\Phi_{(X,\o{X}_1,\tilde{j}_{X_1})}$ and $\Phi_{(X,\o{X}_2,\tilde{j}_{X_2})}$ represent the same weak equivalence of pro-objects.
\end{proof}
Therefore, we can label the above weak equivalence $\Phi_X$, and a particular representation of the morphism will be labeled $\Phi_{(X,\o{X},\tilde{j}_{X})}$.

\begin{lemma}\label{lem:nattrans}
The collection $\{\Phi_{X}\}_{X\in ob(\mathcal{V})}$ defines natural weak equivalences:
\[F_!\Rightarrow G_!;\quad F^!\Rightarrow G^!;\quad F_w\Rightarrow G_w
\]
\end{lemma}
\begin{proof}
Given a closed immersion $f:Z\hookrightarrow X$, the morphism $P_f:P_Z\to P_X$ is represented by morphisms of the form \[(Z,\o{Z},\tilde{j}_Z)\to (X,\o{X},\tilde{j}_X)\]
Therefore, for a particular choice of $\o{X},\tilde{j}_X$ and $\tilde{j}_Z$, we have the commutative diagram:
\[\begin{tikzcd}
\widetilde{\o{Z}-Z}_\bullet\arrow[rrr]\arrow[ddd]\arrow[dr,"\tilde{j}_Z"]&&&\widetilde{\o{X}-X}_\bullet\arrow[dl,"\tilde{j}_X",swap]\arrow[ddd]\\ 
&\tilde{Z}_\bullet\arrow[d]\arrow[r]&\tilde{X}_\bullet\arrow[d]&\\
 &\o{Z}\arrow[r]& \o{X}&\\
 \o{Z}-Z\arrow[rrr]\arrow[ur,"j_Z
 "]&&&\o{X}-X\arrow[ul,"j_X
 ",swap]\\
\end{tikzcd}
\]
Applying $R_{tr}^c$ and localizing by $C^{\mathcal{M}}$, we have the square in $\mathrm{C}^{\mathrm{gm}}(R)$:
\[\begin{tikzcd}
R^c_{tr}[Cone(\tilde{j}_Z)]\arrow[r,"F_!"]\arrow[d,"\Phi_{(Z,\o{Z},\tilde{j}_Z)}"]&R^c_{tr}[Cone(\tilde{j}_X)]\arrow[d,"\Phi_{(X,\o{X},\tilde{j}_X)}"]\\
R^c_{tr}[Z]\arrow[r,"G_!"]&R^c_{tr}[X]
\end{tikzcd}
\]
This diagram represents the diagram of pro-objects:
\[
\begin{tikzcd}
F_!(Z)\arrow[r,"F_!(f)"]\arrow[d,"\Phi_Z"]&F_!(X)\arrow[d,"\Phi_X"]\\
G_!(Z)\arrow[r,"G_!(f)"]&G_!(X)
\end{tikzcd}
\]
So, we have a natural transformation $F_!\Rightarrow G_!$ comprised of weak equivalences. In order to construct the natural transformation $F^!\Rightarrow G^!$, we observe that for an open embedding $i:U\to X$, we showed that the associated morphism $P_U\to P_X$ in $PH(\mathcal{V})$ is represented by a pro-object of morphisms of the form $(X, \o{X},\tilde{j}_X)\to (U, \o{X},\tilde{j}_U)$, where $\o{X}$ is a compactification of $X$, $\tilde{j}_U$ is a smooth projective hyperenvelope of $j_U:\o{X}-U\to\o{X}$, $\tilde{j}_X$ is a smooth projective hyperenvelope of $j_X:\o{X}-X\to\o{X}$ and the morphism $\tilde{j}_X\to \tilde{j}_U$ commutes with the inclusion $j_X\to j_U$, i.e. we get the diagram:
\begin{equation}\label{eq:openhyp}
\begin{tikzcd}
\widetilde{\o{X}-U}_\bullet\arrow[ddd]\arrow[dr,"\tilde{j}_U"]&&&\widetilde{\o{X}-X}_\bullet\arrow[dl,"\tilde{j}_X",swap]\arrow[ddd]\arrow[lll]\\ 
&\tilde{X}_\bullet'\arrow[d]&\tilde{X}_\bullet\arrow[d]\arrow[l]&\\
 &\o{X}\arrow[r,equal]& \o{X}&\\
 \o{X}-U\arrow[ur,"j_U
 "]&&&\o{X}-X\arrow[ul,"j_X",swap]\arrow[lll]\\
\end{tikzcd}
\end{equation}
In turn, applying $R^c_{tr}$ and localizing by $C^{\mathcal{M}}$, we have:
\begin{equation}\label{eq:conesquare}\begin{tikzcd}
M_{dg}^c(Cone(\tilde{j}_U))\arrow[d,"\simeq"]&M_{dg}^c(Cone(\tilde{j}_X))\arrow[l]\arrow[d,"\simeq"]\\
M_{dg}^c(Cone(j_U))&M_{dg}^c(Cone(j_X))\arrow[l]
\end{tikzcd}
\end{equation}
where the vertical arrows are quasi-isomorphisms. We just need to show that the following commutative diagram exists:
\begin{equation}\label{eq:targetsquare}
\begin{tikzcd}
M_{dg}^c(Cone(j_U))\arrow[d,"\simeq"]&M_{dg}^c(Cone(j_X))\arrow[l,"F_!(i)",swap]\arrow[d,"\simeq"]\\
M_{dg}^c(U)&M_{dg}^c(X)\arrow[l,"G_!(i)"]
\end{tikzcd}
\end{equation}
Observe that $M_{dg}^c(Cone(j_U))=Cone(M_{dg}^c(j_U))$ is $M_{dg}^c(\o{X})$ in degree $0$ and $M_{dg}^c(\o{X}-U)$ in degree $1$ (and zero in every other degree), and likewise $M_{dg}^c(Cone(j_X))$ is $M_{dg}^c(\o{X})$ in degree $0$ and $M_{dg}^c(\o{X}-X)$ in degree $1$. Further, the vertical morphisms in the above diagram will be represented by morphisms defined by the pullback on algebraic cycles induced by the open embeddings $U\to \o{X}$ and $X\to\o{X}$ in degree $0$, i.e. we have a commutative square:
\begin{equation}\label{eq:degzero}
\begin{tikzcd}
M^c_{dg}(\o{X})\arrow[d]&M^c_{dg}(\o{X})\arrow[l,equal]\arrow[d]\\
M^c_{dg}(U)&M^c_{dg}(X)\arrow[l]
\end{tikzcd}
\end{equation}
since the complexes in the bottom row are only in degree $0$, we see that \eqref{eq:degzero} gives us the degree $0$ part of the maps defining the square \eqref{eq:targetsquare}. The degree one part of \eqref{eq:targetsquare} is simply:
\[
\begin{tikzcd}
    M^c_{dg}(\o{X}-U)\arrow[d]&M^c_{dg}(\o{X}-X)\arrow[l]\arrow[d]\\
    0\arrow[r,equal]&0
\end{tikzcd}
\]
which commutes, and since the bottom row is $0$, the degree one part of \eqref{eq:targetsquare} will trivially commute with the chain maps. So, \eqref{eq:targetsquare} commutes. Composing \eqref{eq:conesquare} with \eqref{eq:targetsquare} gives us:
\[
\begin{tikzcd}
M^c_{dg}(Cone(\tilde{j}_U))\arrow[d,"\Phi_{(U,\o{X},\tilde{j}_U)}"]&M^c_{dg}(Cone(\tilde{j}_X))\arrow[l]\arrow[d,"\Phi_{(X,\o{X},\tilde{j}_X)}"]\\
M^c_{dg}(U)&M^c_{dg}(X)\arrow[l]
\end{tikzcd}
\]
which represents the morphism:
\[
\begin{tikzcd}
F^!(U)\arrow[d,"\Phi_U"]&F^!(X)\arrow[d,"\Phi_X"]\arrow[l,"F_!(i)",swap]\\
G^!(U)& G^!(X)\arrow[l,"G_!(i)"]
\end{tikzcd}
\]
This gives the natural transformation $F^!\Rightarrow G^!$ as desired. $F_w\Rightarrow G_w$ follows from $F_!\Rightarrow G_!$. 
\end{proof}

\begin{corollary}
The natural isomorphisms in Lemma \ref{lem:nattrans} give us the equivalence of the derived motivic measures in Proposition \ref{BGN} and Proposition \ref{prop:Mdg}. 
\end{corollary}
\begin{proof}
Lemma \ref{lem:nattrans} shows that the maps on $K$-theory induced by $F$ and $G$ are homotopy equivalent, i.e $K(F)\simeq K(G)$. One can then check that we therefore have the homotopy commutative diagram:
\[\begin{tikzcd}
&K(\Pro^{wc}(\s\V_*^{\perf}))\arrow[rrr]&&&K(\Pro^{wc}(\mathrm{C}^{\mathrm{gm}}(R)))\\
K(\mathcal{V})\arrow[ddr,swap]\arrow[ur]&K(\s\V_*^{\perf})\arrow[u,swap,"\simeq"]\arrow[dd,"K(M_{dg}\circ N)"]&&\\
\\
&K(\mathrm{C}^{\mathrm{gm}}(R))\arrow[rruruu,swap,"\simeq"]&&&
\end{tikzcd}\]
which verifies the claim.
\end{proof}

\section{From Voevodsky to Chow Motives}
It would be nice to have a map $K(\mathrm{C}^{\mathrm{gm}}(\Z))\to K(\Ch^\perf(\Chow^{\eff}))$ that would allow us to recover the derived Gillet--Soul\'{e} measure in Corollary \ref{cor:derivedGS}. For the rest of the section, we write $\mathrm{C}^{\mathrm{gm}}(\Z)$ as $\mathrm{C}^{\mathrm{gm}}$. The following work is largely based on 6.7.4 of \cite{BV}:
\begin{proposition}\label{thm:weightstructure}
There is a map of K-theory spectra: $K(\mathrm{C}^{\mathrm{gm}})\to K(\Ch^\perf(\Chow^{\eff}))$ 
\end{proposition}
\begin{proof}
Let $DM^{\eff}_{gm,Nis}(\V)$ denote the full subcategory of $DM^{\eff}_{gm,Nis}$ generated by direct summands of Voevodsky motives of smooth projective varieties. Let $\Big(DM^{\eff}_{gm,Nis}(\V)\Big)^{pretr}$ denote the smallest pre-triangulated subcategory of $DM^{\eff}_{gm,Nis}$ containing $DM^{\eff}_{gm,Nis}(\V)$, and let $\Big(DM^{\eff}_{gm,Nis}(\V)\Big)^{\perf}$ denote the full subcategory of $DM^{\eff}_{gm,Nis}$ consisting of objects quasi-isomorphic to objects of $\Big(DM^{\eff}_{gm,Nis}(\V)\Big)^{pretr}$. By the proof of 6.7.4 of \cite{BV}, the fully faithful inclusion:
\[
\Big(DM^{\eff}_{gm,Nis}(\V)\Big)^{\perf}\hookrightarrow DM^{\eff}_{gm,Nis}
\]
induces an equivalence on homotopy categories. This implies that if $\mathrm{C}^{\mathrm{gm}}(\V)$ denotes the presheaf category of $\Big(DM^{\eff}_{gm,Nis}(\V)\Big)^{\perf}$, then the inclusion $\mathrm{C}^{\mathrm{gm}}(\V)\to \mathrm{C}^{\mathrm{gm}}$ is weakly exact, and induces an equivalence of homotopy categories, and reflects weak equivalences. So, by Theorem 1.3 of \cite{BM}, we have that the induced map of K-theory spectra $K(\mathrm{C}^{\mathrm{gm}}(\V))\to K(\mathrm{C}^{\mathrm{gm}})$ is an equivalence. Now, in \cite{BV} the DG-functor sending $DM^{\eff}_{gm,Nis}(\V)$ to its homotopy category is denoted by
\[\epsilon_0:DM^{\eff}_{gm,Nis}(\V)\to \Chow^{\eff}\]
It is the identity on objects and sends a map $f\mapsto H^0(f)\in \Hom_{\Chow^{\eff}}(X,Y)$. This in turn induces a DG functor
\[\epsilon_0^{\perf}:\Big(DM^{\eff}_{gm,Nis}(\V)\Big)^{\perf}\to \Ch^{\perf}(\Chow^{\eff})\]
which is weakly exact and therefore gives a map of K-theory spectra $K(\mathrm{C}^{\mathrm{gm}}(\V))\to K(\Ch^{\perf}(\Chow^{\eff})^{op}-mod)$. So far, we have the map:
\[K(\mathrm{C}^{\mathrm{gm}})\xleftarrow{\simeq} K(\mathrm{C}^{\mathrm{gm}}(\V))\to  K(\Ch^{\perf}(\Chow^{\eff})^{op}-mod)
\]
Finally, viewing $\Ch^{\perf}(\Chow^{\eff})$ as a DG category (as it is a subcategory of complexes over an additive category), the DG-Yoneda embedding $\Ch^{\perf}(\Chow^{\eff})\to \Ch^{\perf}(\Chow)^{op}-\mathrm{mod}$ induces a weak equivalence on homotopy categories \cite{To} and one can quickly check that the rest of the axioms necessary to invoke 1.3 of \cite{BM} are satisfied for the embedding, so the induced map $K(\Ch^{\perf}(\Chow^{\eff}))\to K(\Ch^{\perf}(\Chow^{\eff})^{op}-mod)$ is an equivalence. Therefore, we have the map:
\[K(\mathrm{C}^{\mathrm{gm}})\xleftarrow{\simeq}K(\mathrm{C}^{\mathrm{gm}}(\V))\xrightarrow{K(\epsilon_0^{\perf})}K(\Ch^{\perf}(\Chow^{\eff})^{op}-mod)\xleftarrow{\simeq}K(\Ch^{\perf}(\Chow^{\eff}))
\]
Concluding the proof.
\end{proof}

\begin{theorem}\label{thm:RecoverChow}
The map of K-theory spectra in Proposition \ref{prop:Mdg} (for $R=\Z$) fits into a homotopy commutative diagram:
\[
\begin{tikzcd}
K(\s\V_*^{\perf})\arrow[r]\arrow[dr]&K(\mathrm{C}^{\mathrm{gm}})\arrow[d]\\
&K(\Ch^\perf(Chow^\eff))
\end{tikzcd}
\]
\end{theorem}
\begin{proof}
First observe that the weakly exact functor $\s\V_*^{\perf}\to \mathrm{C}^{\mathrm{gm}}$ factors through $\mathrm{C}^{\mathrm{gm}}(\V)$, and therefore we have the homotopy commutative diagram:
\[
\begin{tikzcd}
K(\s\V_*^{\perf})\arrow[r]\arrow[dr]&K(\mathrm{C}^{\mathrm{gm}})\\
&K(\mathrm{C}^{\mathrm{gm}}(\V))\arrow[u,"\simeq"]
\end{tikzcd}
\]
Then, we have the commutative diagram of weakly exact functors
\[
\begin{tikzcd}
\s\V_*^{\perf}\arrow[r]\arrow[dr]&\Big(DM_{eff,Nis}^{gm}(\V)\Big)^{pretr}\arrow[d,"\epsilon_0^{pretr}"]\arrow[r]&\mathrm{C}^{\mathrm{gm}}(\V)\arrow[d]\\
&\Ch^{\perf}(Chow)\arrow[r]&\Ch^{\perf}(Chow)\mathrm{-dgm}
\end{tikzcd}
\]
By the work done in the proof of Proposition \ref{thm:weightstructure}, the bottom functor induces an equivalence on $K$-theory spectra, so we get the homotopy commutative diagram:
\[
\begin{tikzcd}
K(\s\V_*^{\perf})\arrow[r]\arrow[dr]&K(\mathrm{C}^{\mathrm{gm}})\arrow[d]\\
& K(\Ch^{\perf}(Chow))
\end{tikzcd}
\]
which proves the theorem. 
\end{proof}
\begin{remark}\label{rem:thmheart}
    If we denote $\Chow^{\eff}(R)$ to be the category of Chow motives with coefficients in a commutative ring $R$, i.e. whose objects are the same as $\Chow^{\eff}$ but where
    \[
    \Hom_{\Chow^{\eff}(R)}(X,Y):=\Hom_{\Chow^{\eff}}(X,Y)\otimes R
    \]
    the proof of Proposition \ref{thm:weightstructure} shows us that we have a commutative diagram, for each $R$:
    \[
    \begin{tikzcd}    K(\s\V_*^{\perf})\arrow[r]\arrow[dr]&K(\mathrm{C}^{\mathrm{gm}}(R))\arrow[d]\\
&K(\Ch^\perf(Chow^\eff(R)))
    \end{tikzcd}
    \]
\end{remark}
\section{The Derived Compactly Supported $\A^1$ Euler Characteristic}
We briefly review the work of \cite{Ro} which, with slight modification, provides a map of $S^1$-spectra 
\[K(\sShv_{\cdh}(Sm)_*)\to \mathbf{1}_{\mathrm{SH}_{Nis}(k)}\] 
The work done in the previous sections will promote this to a derived motivic measure
\[
K(\mathcal{V}_k)\to \mathbf{1}_{\mathrm{SH}_{Nis}(k)}
\]
At the level of $\pi_0$, this assigns to a variety its compactly supported $\A^1$ Euler characteristic \[[X]\mapsto \chi_{c,\A^1}(X)\in \GW(k)=\pi_0(\mathbf{1}_{\mathrm{SH}_{Nis}(k)})\]

Let $\sShv_{Nis,\A^1}(Sm)$ denote the category of simplicial \emph{Nisnevich} sheaves on $Sm$ with an $\A^1$-local model structure due to \cite{Jar}. For a review of the Nisnevich topology, see \cite{Vo}. This can be defined as the Bousfield localization of the model structure on $\sShv_{Nis}(Sm)_*$ as defined in Proposition \ref{prop:joyal} with respect to the projections $X\times\A^1\to X$ for all $X\in Sm_k$ (see 2.5 of \cite{MV}). 
\begin{proposition}\label{prop:10.1}
    The inclusion $Sm\hookrightarrow \Sch$ gives us a direct image functor:
    \[f_*:\sShv_{\cdh}(\Sch)_*\to \sShv_{Nis}(Sm)_*\]
    Further, composing $f_*$ with the $\A^1$-localization functor:
    \[L_{\A^1}:\sShv_{Nis}(Sm)_*\to \sShv_{Nis,\A^1}(Sm)_*\] gives an exact functor of Waldhausen categories (of their perfect cofibrant objects).
\end{proposition}
\begin{proof}
    Because all cdh covers are Nisnevich covers, we have a morphism of sites $(Sm)_{Nis}\hookrightarrow (\Sch)_{\cdh}$. This gives a direct image functor 
    \[
    \begin{tikzcd}
    & \sShv_{\cdh}(Sm)_*\arrow[d,"f_*"]\\      Sm\arrow[ur]\arrow[r]&\sShv_{Nis}(Sm)_*
    \end{tikzcd}
    \]
    Note that the direct image will preserve limits, finite colimits (i.e. pushouts along cofibrations) and local weak equivalences. Next, by 2.5 of \cite{MV} the functor $L_{\A^1}$ is a left Bousfield localization, meaning that it preserves weak equivalences and cofibrations (and pushouts along them). So, we get an exact functor of Waldhausen categories as desired.  
\end{proof}

In \cite{Ro}, the $K$-theory of the unstable $\A^1$ category is constructed using a different model structure, which we define below:

\begin{definition}[\cite{Ro} Definition 3.2]
   Let $\textbf{M}(k)$ denote the category $\s\Pre(Sm)_*$ with a Nisnevich $\A^1$-local model structure defined as follows:
   cofibrations are given as 
   \[\{(X\times\partial\Delta^n\hookrightarrow\Delta^n))_+\}_{X\in Sm_k}\]
   The small object argument gives a cofibrant replacement functor $(-)^c$. An object $B$ is fibrant if
   \begin{enumerate}
       \item $B(X)$ is a fibrant simplicial set for all $X\in Sm_k$.
       \item The image of every \emph{Nisnevich elementary distinguished square}, i.e a pullback diagram:
       \[
       \begin{tikzcd}
           V\arrow[r]\arrow[d]&Y\arrow[d]\\
           U\arrow[r]&X
       \end{tikzcd}
       \]
       where $Y\to X$ is \'{e}tale and $Y-V\to X-U$ is an isomorphism, is sent to a homotopy Cartesian square.
       \item  For every $X\in Sm_k$, $B(X\times\A^1\xrightarrow{pr}X)$ is a weak equivalence of simplicial sets.
   \end{enumerate}
   A morphism $\phi:D\to B$ is a weak equivalence if for every fibrant $Z$, we have
   \[
   sSet_{\s\Pre(Sm)_*}(\phi^c,Z):sSet_{\s\Pre(Sm)_*}(B^c,Z)\to sSet_{\s\Pre(Sm)_*}(D^c,Z)
   \]
   is a weak equivalence of simplicial sets. A map is a fibration if it satisfies the right lifting property with respect to acyclic cofibrations. 
\end{definition}
Luckily, the $K$-theory of perfect cofibrant objects associated to this model structure is the same as that associated to the model structure in Proposition \ref{prop:10.1}.
\begin{proposition}\label{prop:10.3}
    $\textbf{M}(k)$ is a cofibrantly generated. Further, there is a left Quillen functor $\sShv_{Nis,\A^1}(Sm)_*\to \textbf{M}(k)$ that induces an equivalence on homotopy categories. Restricting to perfect cofibrant objects on both categories, we have an equivalence of $K$-theory spectra:
    \[
    K(\sShv_{Nis,\A^1}(Sm)_*)\to K(\textbf{M}(k))
    \]
\end{proposition}
\begin{proof}
    First, there is cofibrantly generated Nisnevich-local model structure of Jardine on the category of simplicial presheaves (see remark 7.4 of \cite{D} for further discussion), which is denoted $\s\Pre(Sm)_{*,Jardine}$. Proposition 8.1 of \cite{D} tells us that the known Quillen equivalence given by the inclusion $\sShv_{Nis}(Sm)_*\hookrightarrow \s\Pre(Sm)_{*,Jardine}$ remains a Quillen equivalence after we $\A^1$-localize. Next, Theorem 2.17 of \cite{DRO} proves that the identity functor on simplicial presheaves gives a Quillen equivalence between $\s\Pre(Sm)_{*,Jardine,\A^1}$ and $\textbf{M}(k)$. Since all of these categories have FFWC, we see that we get exact functors between their perfect cofibrant objects that induce isomorphisms on homotopy categories. Since weak equivalences are closed under retracts, we can invoke 1.1 of \cite{BM} to get the equivalences of $K$-theory spectra. 
\end{proof}
\begin{corollary}\label{cor:M(k)motivicmeasure}
    Composing the map $K(\mathcal{V}_k)\to K(\sShv_{\cdh}(\Sch)_*)$ constructed in Proposition \ref{prop:sshv} with Proposition \ref{prop:10.1} and \ref{prop:10.3}, we have a new map of spectra:
    \[
    K(\mathcal{V}_k)\to K(\textbf{M}(k))
    \]
\end{corollary}
Immediately, we can also define the following motivic measure constructed in 5.26 of \cite{Cam} that does not rely on a folklore proof (5.24 of loc. cit):
\begin{proposition}
    Let $k=\C$. By 4.2 of \cite{Ro} there is a map of $K$-theory spectra: $K(\textbf{M}(\C))\to A(*)$, where $A(*)$ is Waldhausen's $A$ theory of a point. Composing this map with the one in Corollary \ref{cor:M(k)motivicmeasure} gives a map
    \[K(\mathcal{V}_\C)\to A(*)\]
    which on path components, will send 
    \[[X]\mapsto [X(\C)]\]
    if $X$ is any smooth projective variety. 
\end{proposition}
Now, we move to the stable setting. The following two propositions are 2.6 and 2.7 of \cite{Ro}, but are written in greater generality than we need for our work:
\begin{proposition}\label{prop:stabilization}
    Let $\mathrm{SH}(k)=\textbf{Spt}_{T}(Sm)$ denote the category of $T$ spectra with $T=S^{0,1}\wedge S^{1,1}$, where $S^1=\Delta^1/\partial\Delta^1$ and $S^{1,1}$ is the mapping cylinder of the inclusion $k\xrightarrow{1} \G_m$. $\mathrm{SH}(k)$ admits a level-wise projective model structure such that the infinite suspension:
    \[\Sigma_{+}^\infty:\textbf{M}(k)\to \mathrm{SH}(k)\]
    induces a map of $K$-theory spectra of perfect cofibrant objects, where in $\mathrm{SH}(k)$, perfect objects are finite homotopy colimits of suspension spectra of smooth varieties. Further $\Sigma_+^\infty$ induces an equivalence of $K$-theory spectra via a d\'{e}vissage argument.
\end{proposition}
    However, this model structure does not respect the monoidal structure on $\mathrm{SH}(k)$, which is important for building the higher trace. As a result, we must actually deal with symmetric spectra (which in this case gives a Quillen equivalent model structure).  
\begin{proposition}\label{prop:symspt}
    Let $SymSpt_{T}(Sm)$ denote the category of symmetric $T$-spectra. $\textbf{SymSpt}_{T}(Sm)$ also admits a level-wise projective model structure, making it a monoidal model category Quillen equivalent to $\mathrm{SH}(k)$ (using the assumption that $k$ is a field) giving a $K$-theory equivalence on the subcategory of cofibrant perfect objects. Further, the infinite suspension:
    \[\Sigma_{+}^\infty:\textbf{M}(k)\to \textbf{SymSpt}_{T}(Sm)\]
    is a left Quillen functor which induces an equivalence on $K$-theory of perfect cofibrant objects, i.e.
    \[K(\textbf{M}(k))\xrightarrow{\simeq} K(\textbf{SymSpt}_{T}(Sm))\]
\end{proposition}
The following is 6.6 of \cite{Ro}:
\begin{proposition}\label{cor:trace}
    There is a trace map:
    \[K(\textbf{SymSpt}_{T}(Sm))\to \End(1_{\mathrm{SH}_k})\]
    where $\End(\mathbf{1}_{\mathrm{SH}_k})$ is the $S^1$ endomorphism spectrum of the unit object in $\mathrm{SH}_k$. The homotopy groups of this $S^1$ spectrum are written as $\pi_{*,0}(1_{\mathrm{SH}(k)})$. At the level of $\pi_{0,0}$, this sends the suspension spectrum of a variety $X$ to its compactly supported $\A^1$ Euler characteristic
    \[[\Sigma_{+}^\infty X]\mapsto \chi_{c,\A^1}(X)\in \GW(k)\]
\end{proposition}

\begin{corollary}
    We have a map $K(\mathcal{V})\to K(\textbf{SymSpt}_T(Sm))$ by composing Corollary \ref{cor:M(k)motivicmeasure} with Proposition \ref{prop:symspt}. Composing this map in turn with Proposition \ref{cor:trace} gives a derived motivic measure:
    \[K(\chi_{c,\A^1}):K(\mathcal{V}_k)\to K(\End(\mathbf{1}_{\mathrm{SH}_k}))\] where $K(\End(\mathbf{1}_{\mathrm{SH}_k}))$ is viewed as an $S^1$ spectrum. In particular, on $K$ groups, we have homomorphisms:
    \[K_*(\mathcal{V})\to \pi_{*,0}(\textbf{1}_{\mathrm{SH}(k)})\] 
    Further, on the level of $\pi_0$, we have, for a smooth variety $X$:
     \[[X]\mapsto \chi_{c,\A^1}(X)\in \GW(k)\]
\end{corollary}

\bibliographystyle{alpha}
\bibliography{refs}
\end{document}